\documentclass[12pt, reqno]{amsart}
\usepackage{amssymb}
\usepackage{amsthm}
\usepackage{amscd}

\usepackage{mathtools}
\mathtoolsset{showonlyrefs}
\numberwithin{equation}{section}

\usepackage{graphicx}
\usepackage{hyperref}

\usepackage{amsmath}
\usepackage{epic,eepic}

\usepackage{latexsym}

\usepackage{pifont}

\theoremstyle{plain}
\newtheorem{lemma}{Lemma}[subsection]
\newtheorem{prop}[lemma]{Proposition}
\newtheorem{thm}[lemma]{Theorem}
\newtheorem{cor}[lemma]{Corollary}
\newtheorem{aplemma}{Lemma~A.\hspace{-1.5mm}}
\newtheorem{approp}{Proposition~A.\hspace{-1.5mm}}
\newtheorem{apthm}{Theorem~A.\hspace{-1.5mm}}
\newtheorem{apcor}{Corollary~A.\hspace{-1.5mm}}
\newtheorem{intthm}{Theorem}

\usepackage[all]{xy}

\theoremstyle{definition}

\newtheorem{rema}[lemma]{Remark}

\newtheorem{remb}{Remark}

\newtheorem{defi}[lemma]{Definition}
\newtheorem{exa}[lemma]{Example}
\newtheorem{aprem}{Remark~A.\hspace{-1.5mm}}
\newtheorem{apdefi}{Definition~A.\hspace{-1.5mm}}
\newcommand{\bde}{\begin{defi}}
\newcommand{\ede}{\end{defi}\vspace{1mm}}
\newcommand{\ble}{\begin{lemma}}
\newcommand{\ele}{\end{lemma}}
\newcommand{\bpr}{\begin{prop}}
\newcommand{\epr}{\end{prop}}
\newcommand{\bt}{\begin{thm}}
\newcommand{\et}{\end{thm}}
\newcommand{\bco}{\begin{cor}}
\newcommand{\eco}{\end{cor}}
\newcommand{\bre}{\begin{rema}}
\newcommand{\ere}{\end{rema}}
\newcommand{\brea}{\begin{rema}}
\newcommand{\erea}{\end{rema}\vspace{1mm}}
\newcommand{\breb}{\begin{remb}}
\newcommand{\ereb}{\end{remb}\vspace{1mm}}
\newcommand{\bex}{\begin{exa}}
\newcommand{\eex}{\end{exa}}
\newcommand{\bpf}{\begin{proof}}
\newcommand{\epf}{\end{proof}\vspace{1mm}}

\newcommand{\bade}{\begin{apdefi}}
\newcommand{\eade}{\end{apdefi}}
\newcommand{\bale}{\begin{aplemma}}
\newcommand{\eale}{\end{aplemma}}
\newcommand{\bapr}{\begin{approp}}
\newcommand{\eapr}{\end{approp}}
\newcommand{\bat}{\begin{apthm}}
\newcommand{\eat}{\end{apthm}}
\newcommand{\baco}{\begin{apcor}}
\newcommand{\eaco}{\end{apcor}}
\newcommand{\bare}{\begin{aprem}}
\newcommand{\eare}{\end{aprem}}

\newcommand{\be}{\begin{enumerate}}
\newcommand{\ee}{\end{enumerate}}
\newcommand{\bcd}{\[\begin{CD}}
\newcommand{\ecd}{\end{CD}\]}
\newcommand{\bit}{\begin{itemize}}
\newcommand{\eit}{\end{itemize}}
\newcommand{\bq}{\begin{quote}}
\newcommand{\eq}{\end{quote}}
\newcommand{\ba}{\begin{array}}
\newcommand{\ea}{\end{array}}
\newcommand{\mcA}{\mathcal{A}}

\newcommand{\mcD}{\mathcal{D}}
\newcommand{\mcE}{\mathcal{E}}
\newcommand{\mcF}{\mathcal{F}}

\newcommand{\mcH}{\mathcal{H}}

\newcommand{\mcK}{\mathcal{K}}
\newcommand{\mcL}{\mathcal{L}}
\newcommand{\mcM}{\mathcal{M}}
\newcommand{\mcN}{\mathcal{N}}
\newcommand{\mcO}{\mathcal{O}}
\newcommand{\mcP}{\mathcal{P}}

\newcommand{\mcS}{\mathcal{S}}
\newcommand{\mcT}{\mathcal{T}}
\newcommand{\mcU}{\mathcal{U}}
\newcommand{\mcV}{\mathcal{V}}
\newcommand{\mcW}{\mathcal{W}}

\newcommand{\mbA}{\mathbb{A}}

\newcommand{\mbC}{\mathbb{C}}

\newcommand{\mbG}{\mathbb{G}}
\newcommand{\mbH}{\mathbb{H}}

\newcommand{\mbN}{\mathbb{N}}

\newcommand{\mbP}{\mathbb{P}}

\newcommand{\mbR}{\mathbb{R}}

\newcommand{\mbZ}{\mathbb{Z}}

\newcommand{\migi}{\rightarrow}
\newcommand{\longmigi}{\longrightarrow}

\newcommand{\isom}{\stackrel{\sim}{\migi}}

\newcommand{\migiincl}{\hookrightarrow}

\newcommand{\migisurj}{\twoheadrightarrow}

\newcommand{\Y}{X}

\newcommand{\PP}{\mr{PGL}_{n+1}^\infty}
\newcommand{\QQ}{\mr{PGL}_{n+1}^\mbA}
\newcommand{\RR}{\mr{PGL}_{n+1}^{\mbA, \infty}}

\newcommand{\mr}{\mathrm}
\newcommand{\hidden}[1]{\,}

\newcommand{\N}{N}
\newcommand{\M}{m}
\newcommand{\A}{\mbA}

\newcommand{\Proj}[2]{\mr{Proj}_{#1#2}}
\newcommand{\Aff}[2]{\mr{Aff}_{#1#2}}

\newcommand{\DID}[2]{\mr{ID}^{^\mr{Zzz...}}_{#1#2}}
\newcommand{\DAID}[2]{\mr{AID}^{^\mr{Zzz...}}_{#1#2}}

\newcommand{\eDID}[2]{\overline{\mr{ID}}^{^{_\mr{Zzz...}}}_{#1#2}}
\newcommand{\eDAID}[2]{\overline{\mr{AID}}^{^{_\mr{Zzz...}}}_{#1#2}}

\newcommand{\SSP}{\vspace{0mm}}
\newcommand{\LSP}{\vspace{0mm}}

\begin{document}

\title[Projective and affine structures  in positive characteristic I]
{Projective and affine structures  in positive characteristic I: 
Chern class formulas and Characterizations of projective spaces}
\author{Yasuhiro Wakabayashi}
\date{}
\markboth{Yasuhiro Wakabayashi}{}
\maketitle
\footnotetext{Y. Wakabayashi: 
Graduate School of Information Science and Technology, The  University of Osaka, Suita, Osaka 565-0871, Japan;}
\footnotetext{e-mail: {\tt wakabayashi@ist.osaka-u.ac.jp};}
\begin{abstract}
This paper aims to develop a theory of projective and affine structures on higher-dimensional varieties in positive characteristic. This theory deals with Frobenius-projective and Frobenius-affine structures, which have been previously investigated  in the case where the underlying space is a curve. We first provide a description of such structures in terms of Berthelot's higher-level differential operators. That description leads us to obtain a positive characteristic version of Gunning's formulas, which give necessary conditions on Chern classes for the existence of Frobenius-projective and Frobenius-affine structures, respectively. Finally, we establish some characterizations of projective spaces using Frobenius-projective structures.

\end{abstract}
\setcounter{tocdepth}{2}
\tableofcontents 

\setcounter{section}{0}
\section{Introduction}
\SSP

This paper forms the first part of a series in which we 
develop a theory of projective and affine structures on higher-dimensional varieties in positive characteristic\footnote{This paper is one of two main parts of our unpublished paper available at:  {\tt https://arxiv.org/pdf/2011.04846.pdf} (a follow up paper, ~\cite{Wak16}, is in preparation).}. 
This theory deals with  {\it Frobenius-projective structures}  and {\it Frobenius-affine structures}, 
which are a specific type of Frobenius-Ehresmann structure introduced in ~\cite{Wak12} and
 has been extensively  investigated  
 in the case of curves under the correspondence with $\mr{PGL}_2$-opers.

We expect that the geometry of Frobenius-projective and Frobenius-affine structures will expand in a meaningful way  by comparing with what has been shown in the one-dimensional case and with the theory of projective and affine structures on complex manifolds.
As the first step in expanding the theory, 
we prove in this paper  some basic properties of them.

\LSP
\subsection{Review of  complex projective and affine structures} \label{S90hj2}
Here,  let us briefly review what is a projective (resp., an affine) structure defined on a complex manifold.

A {\it    projective} (resp., {\it an affine}) {\it structure} on a complex manifold of dimension $n>0$ is a maximal system of local coordinates in the analytic topology modeled on the complex projective space $\mbP_\mbC^n$ (resp., the complex affine space $\mbA_\mbC^n$) such  that 
on any two overlapping coordinate patches, the change of coordinates may be described   as a projective transformation on $\mbP_\mbC^n$ (resp., an affine transformation on $\mbA^n_\mbC$).
On a complex manifold equipped with a projective  (resp., an affine) structure, there is  a {\it   projective} (resp., an {\it affine}) {\it geometry},  in the spirit of F. Klein's Erlangen program, that locally agrees with the geometry of $\mbP^n_\mbC$ (resp., $\mbA^n_\mbC$).

Projective and affine structures on complex manifolds  
  arise in many areas of mathematics and are still actively being investigated.
For example,  such additional structures on Riemann surfaces (i.e., complex manifolds of dimension $1$)
play a major role in understanding the framework of uniformization theorem.
An important consequence of the uniformization theorem for Riemann surfaces is that any closed Riemann surface 
 is isomorphic either to the complex projective line  $\mbP^1_\mbC$, or to a quotient of the complex affine line $\mbA_\mbC^1  \left(\subseteq \mbP_\mbC^1 \right)$ by a discrete group of translations (i.e., a  $1$-dimensional complex torus), or to a quotient of the $1$-dimensional complex hyperbolic space $\mbH^1_\mbC  \left(\subseteq \mbP_\mbC^1 \right)$ by a torsion-free discrete subgroup of 
$\mr{SU}(1, 1)\cong \mr{SL}_2(\mbR)$.
This implies that by collecting various local inverses of a universal covering map (from $\mbP^1_\mbC$, $\mbA^1_\mbC$,  or $\mbH^1_\mbC$), 
   we have a canonical projective structure on each  closed Riemann surface.
Also, it is shown that $1$-dimensional complex tori are the only closed Riemann surfaces admitting an affine structure.

    In the case of higher dimensions, there are nontrivial obstructions for the existence of a projective structure
because not all  complex manifolds have  such a  structure.
This fact can be verified, e.g., by a very useful formula providing a necessary condition on   
  Chen classes  for the existence of a projective structure,
  which  was  proved by  R. C. Gunning (cf. ~\cite[Part II, Section 8, Theorem 5]{Gun}).

Moreover, 
S. Kobayashi and T. Ochiai (cf.  ~\cite{KO1}, ~\cite{KO2}) studied and classified compact  complex K\"{a}hler(-Einstein) manifolds  admitting 
a projective structure (or more generally, a projective connection, which is  an infinitesimal version of a projective structure).
The standard examples of  them
are as follows:
\begin{itemize}
\item
the $n$-dimensional projective space $\mbP^n_\mbC$;
\item
all \'{e}tale quotients of $n$-dimensional complex tori;
\item
all compact quotients of the $n$-dimensional complex hyperbolic space $\mbH^n_\mbC$ 
by torsion-free discrete subgroup of $\mr{SU}(n,1)  \left(\subseteq \mr{PGL}_{n+1}(\mbC) \right)$.
\end{itemize}
Each manifold in the second admits an affine structure.

In the case of dimension $1$, 
the respective  examples in this list correspond to the closed Riemann surfaces of genus $0$, $1$, and $>1$ in order. 
A result by S. Kobayashi and T. Ochiai  (cf. ~\cite[Corollary 5.3]{KO1} and Remark preceding that corollary) asserts that
the above  standard examples are
 the only compact  K\"{a}hler-Einstein manifolds admitting a projective connection.
In dimension $2$, a compact complex surface admits a projective structure if and only if 
it is one of the standard examples in the case of $n=2$  (cf. ~\cite[Section 1,  Main Theorem]{KO1}).
On the other hand, it was proved (cf. ~\cite[Theorem 7.1]{JR1}) that there is one more  example in dimension $3$, namely \'{e}tale quotients of smooth modular families of false elliptic curves parametrized by a Shimura curve (cf. ~\cite{JR1}).
See, e.g., ~\cite{Ye}, ~\cite{Kl}, and ~\cite{Du2}, for other related studies
in higher dimensions.

\LSP
\subsection{Previous works related to Frobenius-projective and Frobenius-affine structures} \label{S902}
Then, is it possible to develop the geometry in positive characteristic  based on the same viewpoint?
One thing to keep in mind is that due to its analytic formulation, the definitions of  projective and affine structure cannot be adopted, in positive characteristic, as they are.
Indeed, 
if  we had  defined them  in a naive fashion, they would be nothing  other than  the trivial examples of varieties having such structures. 
By taking this unfortunate fact into account, we need to deal with  appropriate replacements.

The positive characteristic analogues discussed in this paper are  
what we call {\it  $F^\N$-projective structures} and {\it $F^\N$-affine structures}, where $\N$ denotes either a positive integer  or $\N=\infty$.
Both notions were  introduced  by  Y. Hoshi (cf. ~\cite{Hos1}, ~\cite{Hos2})
in the case where $\N<\infty$ and the underlying variety is a curve.
Their higher-dimensional extensions are our  central focus, and we 
 attempt to understand what facts, such as differences from the case of complex manifolds,  can be shown about them.

Note that we can find   previous studies of $F^\N$-projective and $F^\N$-affine structures on curves since these are equivalent to   dormant $\mr{PGL}_2$-opers and dormant (generic) Miura $\mr{PSL}_2$-opers, respectively.
Various properties of dormant $\mr{PGL}_2$-opers and their moduli space were discussed 
in the context of $p$-adic Teichm\"{u}ller theory (cf. ~\cite{Mzk2}).
As discussed  in that study,
an understanding about  their moduli
 provides some important perspectives, including the relationship with 
 the Gromov-Witten theory of certain Quot schemes (cf. ~\cite{JP}, ~\cite{J},  ~\cite{Wak20}, ~\cite{Wak10},  ~\cite{Wak5}, ~\cite{Wak33})  and combinatorics of edge-colored graphs, as well as of convex rational polytopes (cf. ~\cite{LO}, ~\cite{Wak2}, ~\cite{Wak8}).
It gives us  an effective way to solve  the counting  problem of dormant $\mr{PGL}_2$-opers (cf. ~\cite{Wak}, ~\cite{Wak6}).
Some of these facts are generalized to 
the theory of dormant $G$-opers
for a semisimple  algebraic group $G$   (cf. e.g., ~\cite{JRXY}, ~\cite{JP},  ~\cite{Wak6},  ~\cite{Wak7},  ~\cite{Wak5}).

On the other hand, by a result of ~\cite{Wak7},  dormant Miura $\mr{PGL}_2$-opers  correspond bijectively to {\it Tango structures}.
A Tango structure is a certain line bundle on an algebraic curve and brings various sorts of pathological phenomena in positive characteristic.
As an application of the theory   of dormant Miura $\mr{PGL}_2$-opers  developed in ~\cite{Wak7},
we achieve a detailed understanding of the moduli space of algebraic surfaces in positive characteristic violating  the Kodaira vanishing theorem. Also, the main result of ~\cite{Wak9} describes dormant Miura $\mr{PGL}_2$-opers, or equivalently Tango structures, by means of   
  solutions to
 the Bethe ansatz equations  for Gaudin model.

\vspace{2mm}
\begin{center}
\begin{picture}(500,150)

\put(45, 65){char.\,$p>0$ analogue}

\put(272, 60){higher level/dim. analogue}

\put(25,120){\fbox{$\begin{matrix}
\text{{\bf projective} (resp., {\bf affine}) {\bf str's}}
\\
\text{on real/complex manifolds} 
\end{matrix}$}}

\put(255,120){\fbox{$\begin{matrix}
\text{{\bf dormant $\mr{PGL}_2$-opers}}
\\
\text{(resp., {\bf gen.\,Miura $\mr{PGL}_2$-opers})} 
 \\
\text{on alg.\,curves in char.\,$p>0$} 
\end{matrix}$}}

\put(110,15){\fbox{$\begin{matrix}
\text{{\bf $F^N$-projective} (resp., {\bf $F^N$-affine}) {\bf str's}}
 \\
\text{on alg.\,varieties  in char.\,$p>0$} 
\end{matrix}$}}

\put(155,95){\vector(2,-3){35}}
\put(145,95){\vector(2,-3){35}}

\put(270,95){\vector(-2,-3){35}}
\put(280,95){\vector(-2,-3){35}}

\end{picture}
\end{center}

\LSP
\subsection{Results in this paper}
  \label{S903}
In what follows, we shall describe some (relatively important)  results of this paper.
Let $p$ be a prime,  $k$  an algebraically closed field of characteristic $p$,  
and $X$ a smooth projective  variety over $k$ of dimension $n>0$.
For a positive integer $\N$, denote by $X^{(\N)}$ the $\N$-th relative Frobenius twist of $X$ (cf. Section \ref{S1}), i.e., the base-change of $X$ via the $\N$-th iterate of the Frobenius automorphism of $k$.
Denote also by $(\mr{PGL}_{n+1})^{(\N)}_X$ the sheaf  on $X^{(\N)}$ represented by  $\mr{PGL}_{n+1}$, which is considered as a sheaf on $X$ via the underlying homeomorphism of the $\N$-th relative Frobenius morphism $F_{X/k}^{(\N)} : X \migi X^{(\N)}$  (cf. Section \ref{S1}). 

We define a {\it Frobenius-projective} 
{\it structure of level $\N$}, or an {\it $F^\N$-projective structure} for short,  
 on $X$   as a maximal  collection 
  of 
 \'{e}tale coordinate  charts on $X$
 forming a torsor modeled on $(\mr{PGL}_{n+1})^{(\N)}_X$  (cf. Definition \ref{D0188} for the precise definition). 
The local sections of $(\mr{PGL}_{n+1})^{(\N)}_X$ 
play the same  role as projective transformations, appearing  in the definition of a projective structure.
 Moreover, an {\it $F^\infty$-projective structure} means a compatible collection of $F^\N$-projective structures  for various $\N$. 
In a similar manner,  one  can define the  ``affine'' versions, i.e., the notions of an {\it $F^\N$-affine structure}  for various $\N \in \mbZ_{>0} \sqcup \{ \infty \}$.

Thus, we obtain the set
\begin{align}
\mr{Proj}_{\N} \ \left(\text{resp.,}  \  \mr{Aff}_\N \right)
\end{align}
(cf. \eqref{Eg61}, \eqref{Ew011})
of all $F^\N$-projective (resp., $F^\N$-affine)  structures on $X$.
Notice, moreover (cf. \eqref{Eg62}) that we can associate an $F^\N$-projective structure to any $F^\N$-affine structure.

We first consider a description of $F^\N$-projective and $F^\N$-affine structures 
in terms of Berthelot's sheaf of differential operators $\mcD_X^{(\N-1)}$ of level $\N-1$.
 (The case of  infinite level, i.e., the sheaf $\mcD_X^{(\infty)}$, is due to Grothendieck).
By an  {\it  indigenous} (resp.,  {\it affine-indigenous}) {\it $\mcD_{X}^{(\N-1)}$-module}, we mean a certain $\mcD_X^{(\N -1)}$-module  equipped with a filtration satisfying a strict form  of Griffiths transversality (cf. Definitions \ref{D01234} and  \ref{D01235}).
Also, such a  $\mcD_X^{(\N -1)}$-module is called {\it dormant} if
the  $p$-$(\N -1)$-curvature (i.e., the higher-level generalization of $p$-curvature)
  associated to  it vanishes identically 
 (cf. Definition \ref{Dbb030}).
If the underlying space $X$ has dimension $1$, 
then dormant indigenous $\mcD_X^{(\N -1)}$-bundles are nothing but dormant $\mr{PGL}_2^{(\N)}$-opers in the sense of ~\cite{Wak18}  and ~\cite{Wak19}.
Regarding the general case,  we can  provide a bijective correspondences,  described  below.

\SSP
\begin{intthm}[cf. Theorem \ref{P01111} and Proposition \ref{C06p0gh} for the full statement] \label{TGG}
Let $\N \in \mbZ_{> 0} \sqcup \{ \infty \}$.
\begin{itemize}
 \item[(i)]
 Denote by $\overline{\mr{ID}}_{\N -1}^{^{_\mr{Zzz...}}}\!$ 
 the  set of equivalence classes (relative to the equivalence relation determined by tensoring with a line bundle in a certain way) of dormant indigenous $\mcD_X^{(\N -1)}$-modules (cf.  \eqref{Ef1113}).
 Suppose that  $p \nmid (n+1)$.
 Then, there exists a canonical bijection of sets
 \begin{align}
 \overline{\zeta}^{\diamondsuit \Rightarrow \heartsuit}_\N : \overline{\mr{ID}}_{\N -1}^{^{_\mr{Zzz...}}} \isom \mr{Proj}_\N.
 \end{align}
 \item[(ii)]
 Denote by $\overline{\mr{AID}}_{\N -1}^{^{_\mr{Zzz...}}}\!$ the  set of equivalence classes (relative to the equivalence relation determined by tensoring with a line bundle in a certain way) of dormant  affine-indigenous $\mcD_X^{(\N -1)}$-modules (cf.  \eqref{Ef1113}).
Then, there exists a canonical bijection of sets
  \begin{align}
 \overline{\xi}^{\diamondsuit \Rightarrow \heartsuit}_\N : \overline{\mr{AID}}_{\N -1}^{^{_\mr{Zzz...}}} \isom \mr{Aff}_\N.
 \end{align}
 \end{itemize}
 \end{intthm}
\SSP

\begin{rema}[Extension to infinite level]
Before describing the next  theorem,
we shall comment on the advantage of introducing the various structures of infinite level, not only those of finite levels.
A key ingredient in the study of $F^\infty$-projective structures is the Tannakian fundamental group associated to the category of $\mcD_X^{(\infty)}$-modules, or equivalently the  category of $F$-divided sheaves (cf. ~\cite[Section 2.2]{dS}), which we call the {\it stratified fundamental group} and denote by $\pi_1^\mr{str}(X)$.
That is to say, the category of representations of $\pi_1^\mr{str}(X)$ classifies  $\mcD_X^{(\infty)}$-modules.
The dormant indigenous $\mcD_X^{(\infty)}$-module corresponding, via the bijection $\overline{\zeta}_\infty^{\diamondsuit \Rightarrow \heartsuit}$ obtained in  the above theorem, to each $F^\infty$-projective structure specifies a $\mr{PGL}_{n+1}$-representation of
$\pi_1^\mr{str}(X)$ up to conjugation;
it may be thought of as the monodromy representation of this $F^\infty$-projective structure, and enables us to obtain a much better understanding of
related things from  the viewpoint of algebraic groups.
\end{rema}
\SSP

Next, we recall Gunning's formula, that  provides  necessary conditions on Chern classes  for the existence of a projective (resp., an affine) structure on complex manifolds  (cf. ~\cite[Part II, Section 8, Theorem 5]{Gun}).
This   formula implies that
if a K\"{a}hler manifold  carries a projective (resp., an affine) structure, then the Chern numbers are proportional to those of the projective  (resp.,  affine) space of the same dimension.
As an application of the above theorem,
we can prove 
a version for $F^\N$-projective (resp., $F^\N$-affine) structures
 by examining 
 the Chern classes of the corresponding 
  dormant indigenous  (resp.,  dormant affine-indigenous) $\mcD_X^{(\N-1)}$-module; this result is very helpful in characterizing higher-dimensional varieties admitting an $F^\N$-projective (resp., $F^\N$-affine) structure.

\SSP
\begin{intthm}[cf. Theorems \ref{T016} and \ref{T016a} for the full statement] \label{TB}
 Let $\N \in \mbZ_{>0}\sqcup \{ \infty \}$ and let $X$ be a smooth projective variety over $k$  of dimension $n \left(>0\right)$.
We denote by $c_l^\mr{crys} (X)$ (for $l >0$) the $l$-th  crystalline Chern class of $X$.
\begin{itemize}
\item[(i)]
Suppose that $p \nmid (n+1)$ and $X$ admits an $F^\N$-projective structure.
Then,   for each $l > 0$, the  equality 
\begin{align}
c_l^\mr{crys} (X) = \frac{1}{(n+1)^l} \cdot \binom{n+1}{l} \cdot c_1^\mr{crys} (X)^l 
\end{align}
holds in
the $2l$-th crystalline cohomology $H^{2l}_{\mr{crys}}(X/W_\N)$ of $X$ over $W_\N$ (:= the ring of Witt vectors over $k$ of length $\N$).
\item[(ii)]
Suppose that $X$ admits an $F^\N$-affine structure.
Then, 
 for each $l > 0$, the  equality 
\begin{align}
c_l^\mr{crys} (X)  =0 
\end{align}
holds in
$H^{2l}_{\mr{crys}}(X/W_\N)$.
\end{itemize}

 \end{intthm}
\SSP

The most typical examples in our discussion are projective and affine spaces.
They have 
   globally defined coordinate charts taking values in themselves,  and specify   the trivial $F^\N$-projective and $F^\N$-affine structures,  respectively.
 The projective space $\mbP^n$ (with $p \nmid (n+1)$) has only this trivial  $F^\N$-projective  structure,
but the situation is different for the affine space $\mbA^n$, i.e., there are various nontrivial $F^\N$-affine structures on $\mbA^n$.
This fact is  an exotic  phenomenon in positive characteristic
 because  the complex affine space has only the trivial  affine structure.

In the case of $\N = \infty$,
 we prove some characterizations of projective spaces in terms of $F^\infty$-projective structures.
These characterizations can be regarded as the positive characteristic versions of ~\cite[Theorem 4.4]{KO1} and ~\cite[Theorem 4.1]{JR1} (or the main theorem in ~\cite{Ye}).
As investigated in ~\cite{Wak16}, this result will contribute to (the positive characteristic version of)  the classification problem, starting with S. Kobayashi and T. Ochiai, of varieties admitting projective or affine structures (cf. ~\cite{KO1}, ~\cite{KO2}).

\SSP
\begin{intthm} [cf. Proposition \ref{P0104}, Theorems   \ref{T01022},  \ref{T044}, and Remark \ref{R1k3401}] \label{TC}
Let $n$ be a positive integer.  Then, the following assertions hold:
\begin{itemize}
\item[(i)]
The projective space $\mbP^n$ admits exactly one $F^\infty$-projective structure, but  admits  no $F^\infty$-affine structures.
Moreover, the affine space $\mbA^n$ admits an $F^\infty$-affine structure. 
In particular, it admits  an $F^\infty$-projective structure.
\item[(ii)]
Let $X$ be a smooth projective  variety over $k$ of dimension $n$.
Suppose that  $X$  admits  an $F^\infty$-projective structure and moreover
one of the following conditions is fulfilled:
\begin{itemize}
\item
 $X$ contains a rational curve;
\item
the \'{e}tale fundamental group $\pi_1^{\text{\'{e}t}} (X)$ of $X$ is trivial;
\item
the stratified fundamental group $\pi_1^\mr{str} (X)$ of $X$ is trivial.
\end{itemize}
Then, $X$ is isomorphic to the projective space $\mbP^n$.
In particular, $\mbP^n$ is the only Fano variety admitting an $F^\infty$-projective structure.

\end{itemize}
 \end{intthm}

\LSP
\subsection{Notation and Conventions} \label{S9025}

Throughout this paper, we fix
 a prime number $p$ and  an algebraically closed field $k$ of characteristic $p$.
(Our discussion will work  even if  $k$ is replaced with  a  separably closed field of characteristic $p$.)



By a {\it variety (over $k$)}, we mean 
 a connected integral  scheme of finite type over $k$.
Moreover,  a {\it curve}
 means a variety of dimension $1$.
Unless stated otherwise, we will always be working  over $k$; for example,  products  of varieties will be taken over $k$, i.e., $X_1 \times X_2 := X_1 \times_k X_2$.

We fix a positive integer $n$ and a smooth variety $X$ over $k$ of dimension $n$.
Let $\Omega_X$ (resp., $\mcT_X$) denote  the sheaf of $1$-forms (resp., the sheaf of vector fields) on $X$ relative to $k$, i.e., $\Omega_{X} := \Omega_{X/k}$ (resp., $\mcT_X := \mcT_{X/k} = \Omega_{X/k}^\vee$).
Since $X$ is smooth,
  the canonical bundle   $\omega_X$ of $X$ (over $k$) can be constructed  by setting $\omega_X := \mr{det}(\Omega_X)  \left(= \bigwedge^{\mr{dim}(X)} \Omega_X \right)$.

Finally, denote by
$\mbP^n$ (resp., $\mbA^n$) the $n$-dimensional projective (resp., affine) space over $k$.
The $k$-rational points of $\mbP^n$ correspond bijectively to the ratios $[a_0: a_1: \cdots : a_n]$ (with $a_0, \cdots, a_n \in k$ and $(a_0, \cdots, a_n)\neq (0, \cdots, 0)$).
We identify  $\mbA^n$ with the open subscheme of $\mbP^n$ consisting of the elements $\left\{[1: a_1: \cdots : a_n] \, | \, a_1, \cdots, a_n \in k \right\}$.

\vspace{10mm}
\section{$F^\N$-projective and $F^\N$-affine structures} \label{S0002}
\SSP

In this section, we introduce 
 Frobenius-projective  and Frobenius-affine structures
on smooth projective varieties of arbitrary dimension (cf. Definition \ref{D0188}), including those of infinite level (cf. Definition \ref{D0050}).
These are particular cases of the Frobenius-Ehresmann structures in ~\cite{Wak12}.
Also, the previous study for the one-dimensional  case  can be found in  ~\cite{Hos1} and ~\cite{Hos2}.

\LSP
\subsection{Frobenius twists and relative Frobenius morphisms} \label{S1}

We shall  write $f : X \migi \mr{Spec}(k)$ for the structure morphism of $X/k$ and 
 $F_{\Y}$ for
the absolute Frobenius endomorphism of $\Y$.
For each positive integer $\N$, the  {\bf $\N$-th Frobenius twist of $\Y$ over $k$} is, by definition,  the base-change $\Y^{(\N)}$ ($:= \Y \times_{k, F_k^\N} k$) of $\Y$ via
 the $\N$-th iterate $F_k^\N$ of $F_k  \left(= F_{\mr{Spec}(k)} \right) : \mr{Spec}(k) \migi \mr{Spec}(k)$.
Denote by $f^{(\N)} : \Y^{(\N)} \migi \mr{Spec}(k)$ the structure morphism of 
$\Y^{(\N)}$.
The {\bf $\N$-th relative Frobenius morphism of $\Y$ over $k$} is  the unique morphism $F_{\Y/k}^{(\N)} : \Y \migi \Y^{(\N)}$ over $k$ that makes the following diagram commute:
\begin{align}
\vcenter{\xymatrix@C=36pt@R=36pt{
\Y \ar@/^10pt/[rrrrd]^{F_\Y^\N}\ar@/_10pt/[ddrr]_{f} \ar[rrd]_{  \ \ \ \  \  \   F_{\Y/k}^{(\N)}} & & &   &   \\
& & \Y^{(\N)}  \ar[rr]_{W_{X/k}}
 \ar[d]^-{f^{(\N)}}  \ar@{}[rrd]|{\Box}  &  &  \Y \ar[d]^-{f} \\
&  & \mr{Spec}(k) \ar[rr]_{F_k^\N} & &  \mr{Spec}(k).
}}
\end{align}
For convenience, we often use the notation $X^{(0)}$  to denote $X$ itself.

Given a connected smooth algebraic group $G$ over $k$,  
 we denote by 
$G_{\Y}$ the \'{e}tale sheaf of groups on $\Y$ represented by $G$.
Moreover, given  a positive integer $\N$, we shall  write
\begin{align} \label{Eg60}
G_{\Y}^{(\N)} := (F^{(\N)}_{\Y/k})^{-1} (G_{\Y^{(\N)}})  \left(\subseteq G_{\Y} \right).
\end{align}
These sheaves constitute, via pull-back by  various relative  Frobenius morphisms,
  a sequence of inclusions
\begin{align} \label{E0078}
G_\Y \supseteq G_\Y^{(1)} \supseteq G_\Y^{(2)} \supseteq \cdots \supseteq G_\Y^{(\N)} \supseteq \cdots.
\end{align}

\LSP
\subsection{$F^\N$-projective and $F^\N$-affine structures of level $\N < \infty$} \label{Sdd1}

Denote by 
$\mr{GL}_{n+1}$ (resp., $\mr{PGL}_{n+1}$) the general (resp., projective) linear group over $k$ of rank $n+1$.  We shall use the notation  $\overline{(-)}$ to denote  the natural quotient $\mr{GL}_{n+1} \migisurj \mr{PGL}_{n+1}$.
 The group $\mr{PGL}_{n+1}$ can be identified with the automorphism group  of $\mbP^n$ in such a way that $\overline{A}([a_0: \cdots: a_n]) :=  [a_0: \cdots :a_n] \overline{{^t A}}$ for each $A \in \mr{GL}_{n+1}$ and $[a_0: \cdots: a_n]\in \mbP^n$.
Also, denote by 
\begin{align}
\QQ
\end{align}
the subgroup of $\mr{PGL}_{n+1}$ consisting of automorphisms $h$ of $\mbP^n$ with $h (\mbA^n) = \mbA^n$.
That is to say, each element of $\QQ$ may be expressed as $\overline{\begin{pmatrix} a & {\bf 0}\\ {^t {\bf a}}& A\end{pmatrix}}$ for some $a \in \mbG_m$,  ${\bf a} \in \mbA^n$, and  $A \in \mr{GL}_n$.

Let us write
\begin{align} \label{Ep120}
\mcP^{\text{\'{e}t}}_X  \  \left(\text{resp.,} \ \mcA^{\text{\'{e}t}}_X\right)
\end{align}
for the 
set-valued \'{e}tale sheaf
on $X$ that assigns, to each \'{e}tale scheme  $U$ over  $X$, the set of \'{e}tale $k$-morphisms
  from $U$ to $\mbP^n$ (resp., $\mbA^n$).
We shall write $\mcP^{\text{\'{e}t}}:= \mcP^{\text{\'{e}t}}_X$ (resp., $\mcA^{\text{\'{e}t}}:=\mcA^{\text{\'{e}t}}_X$)
if there is no fear of confusion.
The graph of
a local section $\phi : U \migi \mbP^n$ (resp., $\phi : U \migi \mbA^n$) of 
this sheaf,
 which will be denoted by 
 $\Gamma_\phi : U \migi U \times \mbP^n$ (resp., $\Gamma_\phi : U \migi U \times \mbA^n$), defines 
 a local section of 
   the trivial $\mbP^n$-bundle $X \times \mbP^n \xrightarrow{\mr{pr}_1} X$ (resp., the trivial $\mbA^n$-bundle $X \times \mbA^n \xrightarrow{\mr{pr}_1} X$) on $X$.

For each positive integer $\N$, the sheaf $\mcP^{\text{\'{e}t}}$ has a left $(\mr{PGL}_{n+1})_X^{(\N)}$-action defined as follows.
Let $U$ be an open subscheme of $X$, $\phi : U \migi \mbP^n$ an element of  $\mcP^{\text{\'{e}t}}(U)$, and
$\overline{A}$ an element of $(\mr{PGL}_{n+1})_X^{(\N)}(U)  \left(= \mr{PGL}_{n+1} (U^{(\N)}) \subseteq  \mr{PGL}_{n+1} (U) \right)$.
Then,  
 one verifies that the composite 
\begin{align} \label{Eg69}
\overline{A} (\phi) : U \xrightarrow{\Gamma_\phi} U \times \mbP^n \xrightarrow{\overline{A}} U \times \mbP^n \xrightarrow{\mr{pr}_2} \mbP^n
\end{align}
belongs to $\mcP^{\text{\'{e}t}}(U)$ (cf. ~\cite[Lemma 2.2.1]{Wak12}).
The assignment $(\overline{A}, \phi) \mapsto \overline{A}  (\phi)$ defines a $(\mr{PGL}_{n+1})_X^{(\N)}$-action on $\mcP^{\text{\'{e}t}}$.
These actions for various $\N$ are compatible with the sequence of inclusions \eqref{E0078} in the case where $G = \mr{PGL}_{n+1}$.
In a similar manner, one can define a  $(\QQ)_X^{(\N)}$-action on  the sheaf $\mcA^{\text{\'{e}t}}$  that is compatible with 
the $(\mr{PGL}_{n+1})_X^{(\N)}$-action on $\mcP^{\text{\'{e}t}}$ via
the natural inclusions $\mcA^{\text{\'{e}t}} \migiincl  \mcP^{\text{\'{e}t}}$ and $(\QQ)^{(\N)}_X \migiincl (\mr{PGL}_{n+1})_X^{(\N)}$.

The following definition generalizes \cite[Definition 2.1]{Hos1} and  ~\cite[Definition 2.1]{Hos2}.

\SSP
\bde
 \label{D0188}
 Let $\mcS^\heartsuit$  be a subsheaf of $\mcP^{\text{\'{e}t}}$  (resp., $\mcA^{\text{\'{e}t}}$) and  $\N$ a positive integer.
Suppose that  $\mcS^\heartsuit$ is closed under the left $(\mr{PGL}_{n+1})_X^{(\N)}$-action on $\mcP^{\text{\'{e}t}}$ (resp., the left $(\QQ)_X^{(\N)}$-action on $\mcA^{\text{\'{e}t}}$).
 Then, 
 we  say that $\mcS^\heartsuit$ is a {\bf Frobenius-projective} (resp.,  {\bf Frobenius-affine}) {\bf  structure of level $\N$} on $X$, or simply an {\bf $F^\N$-projective}  (resp., {\bf $F^\N$-affine}) {\bf  structure} on $X$,
 if
 it forms a $(\mr{PGL}_{n+1})_X^{(\N)}$-torsor  (resp., a $(\QQ)_X^{(\N)}$-torsor) on $X$ with respect to the resulting $(\mr{PGL}_{n+1})_X^{(\N)}$-action (resp., $(\QQ)_X^{(\N)}$-action)  on $\mcS^\heartsuit$.
That is to say, 
we use the term ``{\it $F^\N$-projective} (resp., {\it $F^\N$-affine}) {\it structure}" to refer to
$F^\N$-$(\mr{PGL}_{n+1}, \mbP^n)$-structures  (resp., $F^\N$-$(\QQ, \mbA^n)$-structures)
 in the sense of  ~\cite[Definition 2.2.2]{Wak12}.
     \ede
\SSP

Let $\N$ be a positive integer.
 Denote by
\begin{align} \label{Eg61}
\Proj{\N}{, X}  \ 
 \left(\text{resp.,} \ 
 \Aff{\N}{, X}
  \right), \ \text{or simply} \ 
  \Proj{\N}{}  \ 
 \left(\text{resp.,} \ 
 \Aff{\N}{}
   \right),
\end{align}
the set of $F^\N$-projective  (resp., $F^\N$-affine) structures on $X$.
If $\mcS^\heartsuit$ is  an $F^\N$-affine structure  on $X$, 
then the smallest subsheaf of $\mcP^{\text{\'{e}t}}$ that contains $\mcS^\heartsuit  \left(\subseteq \mcA^{\text{\'{e}t}} \subseteq \mcP^{\text{\'{e}t}} \right)$ and is closed under the $(\mr{PGL}_{n+1})_X^{(\N)}$-action forms an $F^\N$-projective structure $\iota_\N^\heartsuit (\mcS^\heartsuit)$ on $X$.
The resulting assignment $\mcS^\heartsuit \mapsto \iota_\N^\heartsuit (\mcS^\heartsuit)$ gives a map of sets
 \begin{align} \label{Eg62}
 \iota^\heartsuit_\N :
  \Aff{\N}{, X} \migi \Proj{\N}{, X}.
 \end{align}

\LSP
\subsection{$F^\N$-projective and $F^\N$-affine structures of level $\N = \infty$} \label{Sjkl1}

Let $N'$ be another positive integer with $\N \leq \N'$ and 
 $\mcS^\heartsuit$  an $F^{\N'}$-projective (resp.,  $F^{\N'}$-affine) structure on $X$.
Denote by 
\begin{align}
\mcS^\heartsuit |^{\langle \N \rangle}
\end{align}
  the smallest  subsheaf of $\mcP^{\text{\'{e}t}}$ (resp., $\mcA^{\text{\'{e}t}}$) which contains $\mcS^\heartsuit$ and  is closed under the  $(\mr{PGL}_{n+1})_X^{(\N)}$-action (resp., the $(\QQ)_X^{(\N)}$-action).
One verifies that $\mcS^\heartsuit |^{\langle \N \rangle}$ forms an $F^{\N}$-projective (resp.,   $F^{\N}$-affine)  structure on $X$.
We  refer to $\mcS^\heartsuit |^{\langle \N \rangle}$ as the {\bf $\N$-th truncation} of $\mcS^\heartsuit$.
The resulting assignments $\mcS^\heartsuit\mapsto \mcS^\heartsuit |^{\langle \N \rangle}$ for various positive integers $\N$, $\N'$ with $\N'>\N$ give  a projective system of sets
\begin{align} \label{E33421}
\cdots \migi 
\Proj{\N}{}
 \migi  \cdots\migi
 \Proj{3}{}
    \migi 
   \Proj{2}{}
     \migi 
   \Proj{1}{}
     \\
\left(\text{resp.,} \ \cdots \migi 
\Aff{\N}{}
 \migi  \cdots\migi 
 \Aff{3}{}
 \migi
 \Aff{2}{}
   \migi 
   \Aff{1}{}
   \right).
    \notag
\end{align}

\SSP
\bde \label{D0050} 
 A {\bf Frobenius-projective} (resp.,  {\bf Frobenius-affine}) {\bf structure of level $\infty$} on $X$, or simply an {\bf $F^\infty$-projective} (resp., {\bf $F^\infty$-affine})  {\bf  structure} on $X$, is a collection 
 \begin{align}
 \mcS_\infty^\heartsuit := \{\mcS^\heartsuit_\N \}_{\N \in \mbZ_{>0}},
 \end{align}
  where each $\mcS^\heartsuit_\N$ ($\N \in \mbZ_{>0}$) is an $F^\N$-projective  (resp., $F^\N$-affine) structure with
  $\mcS_{\N+1}^\heartsuit |^{\langle \N \rangle} = \mcS_\N^\heartsuit$.
  That is, an $F^\infty$-projective (resp., $F^\infty$-affine) structure is the same as an 
$F^\infty$-$(\mr{PGL}_{n+1}, \mbP^n)$-structures  (resp., $F^\infty$-$(\QQ, \mbA^n)$-structures)
 in the sense of  ~\cite[Definition 2.3.1]{Wak12}.
    \ede
\SSP

Denote by
\begin{align} \label{Ew011}
\Proj{\infty}{, X}
  \left(\text{resp.,} \ 
  \Aff{\infty}{, X}
   \right), \text{or simply} \ 
   \Proj{\infty}{}
  \left(\text{resp.,} \ 
  \Aff{\infty}{}
   \right),
\end{align}
the set of $F^\infty$-projective (resp., $F^\infty$-affine) structures on $X$.
Then, 
$\Proj{\infty}{}$
 (resp., 
 $\Aff{\infty}{}$)
  may be naturally identified with the limit of the projective system \eqref{E33421}, i.e.,
\begin{align} \label{Ew010}
\varprojlim_{\N \in \mbZ_{>0}}
\Proj{\N}{}
  = 
  \Proj{\infty}{} \ 
    \left(\text{resp.,}   \ \varprojlim_{\N \in \mbZ_{>0}} 
    \Aff{\N}{}
     = 
     \Aff{\infty}{}
     \right).
\end{align}
The collection of maps $\{\iota^\heartsuit_\N \}_{\N \in \mbZ_{>0}}$ (cf. \eqref{Eg62}) is compatible with
the projective systems   $\{ 
\Proj{\N}{}
\}_{\N \in  \mbZ_{>0}}$ and  $\{ 
\Aff{\N}{}
\}_{\N \in \mbZ_{>0}}$.
This collection  induces
a map  of sets
\begin{align} \label{Ew001}
\iota^\heartsuit_\infty :
\Aff{\infty}{} 
 \migi 
\Proj{\infty}{}.
\end{align}

\vspace{10mm}
\section{Dormant (affine-)indigenous $\mcD^{(\N-1)}$-modules} \label{S0210}\SSP

This section starts by 
recalling  Berthelot's sheaf  of differential operators of finite level, and then,
defines the notion of  a dormant (affine-)indigenous $\mcD^{(\N-1)}$-module (cf. Definitions \ref{D01234} and \ref{D01235}).
As proved in the next section,  dormant (affine-)indigenous $\mcD^{(\N-1)}$-modules will be used to describe $F^\N$-projective/affine structures in terms of $\mcD$-modules.

\LSP
\subsection{Differential operators of level $\M$} \label{SS018b}

 For each $\M \in \mbZ_{\geq 0} \sqcup \{ \infty \}$ and $l \in \mbZ_{\geq 0}$, denote by
$\mcD_{X, l}^{(\M)}$
 the sheaf of {\it differential operators of level $\M$ and order at most $l$} (hence $\mcD_{X, 0}^{(\M)} = \mcO_X$), and write
\begin{align}
\mcD^{(\M)}_X := \bigcup_{l \geq 0} \mcD_{X, l}^{(\M)}
\end{align}
(cf. ~\cite{Be1}, ~\cite{Be2}).
This sheaf acts on $\mcO_X$ and  forms a sheaf of (noncommutative) $k$-algebras.
We shall write ${^L \mcD}_{X}^{(\M)}$ (resp., ${^R \mcD}_{X}^{(\M)}$) for the sheaf $\mcD_X^{(\M)}$ endowed with a structure of $\mcO_X$-module  arising from  left (resp., right) multiplication by sections in $\mcD_{X, 0}^{(\M)}  \left(=\mcO_X \right)$.

In what follows, let us explain the local description of these structures.
Given an integer $r$, we write
 $q_r^{}$ for  the unique  integer satisfying that $r = q_r^{} \cdot p^\M + l$ with $0 \leq  l < p^\M$, where 
 $q_r^{\M} := 0$ if $\M= \infty$.
For nonnegative integers  $a$, $a'$, and $a''$ with $a = a' + a''$, we set
\begin{align}
\left\{ \begin{matrix} a \\ a' \end{matrix} \right\} := \frac{q_{a}!}{ q_{a'}!  \cdot q_{a''}!}, 
\hspace{10mm}
\left\langle\begin{matrix} a  \\ a' \end{matrix} \right\rangle := \left(\begin{matrix} a  \\ a' \end{matrix} \right)\cdot \left\{\begin{matrix} a  \\ a' \end{matrix} \right\}^{-1}.
\end{align}
We also use the usual conventions on multi-indices;
if two elements $\underline{a} := (a_1, \cdots, a_n)$,  $\underline{a}' := (a'_1, \cdots, a'_n)$ in $\mbZ^n_{\geq 0}$ satisfies $\underline{a}' \leq \underline{a}$, that is, if $a'_i \leq a_i$ for every $i=1, \cdots, n$, then we define
\begin{align}
\underline{a} !:= \prod_{i=1}^n a_i!, 
\hspace{5mm}
\left(\begin{matrix}\underline{a} \\ \underline{a}' \end{matrix} \right) := \prod_{i=1}^n \left(\begin{matrix} a_i \\ a'_i\end{matrix} \right),
\hspace{5mm}
\left\{\begin{matrix}\underline{a} \\ \underline{a}' \end{matrix} \right\} := \prod_{i=1}^n \left\{\begin{matrix} a_i \\ a'_i\end{matrix} \right\},
\hspace{5mm}
\left\langle\begin{matrix}\underline{a} \\ \underline{a}' \end{matrix} \right\rangle := \prod_{i=1}^n \left\langle\begin{matrix} a_i \\ a'_i\end{matrix} \right\rangle.
\end{align}

We shall fix a  local coordinate system $\underline{t} := (t_1, \cdots, t_n)$ 
of $X$ defined on an open subscheme $U \left( \subseteq X \right)$.
Given an $n$-tuple of nonnegative integers $\underline{h} := (h_1, \cdots, h_n) \in \mbZ^n_{\geq 0}$, we write
$\underline{t}^{\underline{h}} := \prod_{i=1}^n t_i^{h_i}  \left(\in \Gamma (U, \mcO_U) \right)$.
Then,  the $\mcO_U$-module ${^L \mcD}_{X}^{(\M)} |_U$ is freely generated by   
 the  symbols $\underline{\partial}^{\langle \underline{r}\rangle} := \partial_1^{\langle r_1\rangle} \cdots \partial_n^{\langle r_n\rangle}$ for various $ \underline{r} := (r_1, \cdots, r_n) \in \mbN^n$, and
the action of $\mcD_X^{(\M)}$ on $\mcO_X$ is determined by
 \begin{align} \label{aaaa}
 \underline{\partial}^{\langle \underline{r} \rangle} (\underline{t}^{\underline{h}}) = q^{}_{\underline{r}}! \cdot \binom{\underline{h}}{\underline{r}} \cdot  \underline{t}^{\underline{h}-\underline{r}}  \left(:= \prod_{i=1}^n q^{}_{r_i}! \cdot  \binom{h_i}{r_i} \cdot t_i^{h_i -r_i} \right).
 \end{align}
 Moreover,   the multiplication in $\mcD_X^{(\M)}$ is given by 
\begin{align}
\underline{\partial}^{\langle \underline{l} \rangle}\cdot
 \underline{\partial}^{\langle \underline{r} \rangle}
  = \left\langle\begin{matrix} \underline{l} + \underline{r} \\ \underline{l} \end{matrix} \right\rangle  \cdot  \underline{\partial}^{\langle \underline{l} + \underline{r} \rangle},
\hspace{10mm}
\underline{\partial}^{\langle \underline{r}\rangle}
\cdot  f = \sum_{\underline{i} \leq \underline{r}}   \left\{\begin{matrix} \underline{r}  \\ \underline{i} \end{matrix} \right\}  \cdot \underline{\partial}^{\langle \underline{i}\rangle} (f) \cdot  \underline{\partial}^{\langle \underline{r} - \underline{i}\rangle}
\end{align}
for any $f \in \Gamma (U, \mcO_X)$.

Given a nonnegative integer  $r < p^{m+1}$ and $i=1, \cdots, n$,
we shall write 
 $\partial^{[r]}_i := \frac{1}{q_r !} \cdot \partial^{\langle r \rangle}_i$, so that, with this notation, 
 the formulas for multiplication in $\mcD^{(\M)}_X$ become
 \begin{align}
\underline{\partial}^{[ \underline{l} ]}\cdot
 \underline{\partial}^{[ \underline{r} ]}
  = \left(\begin{matrix} \underline{l} + \underline{r} \\ \underline{l} \end{matrix} \right)  \cdot  \underline{\partial}^{[ \underline{l} + \underline{r} ]},
\hspace{10mm}
\underline{\partial}^{[ \underline{r}]}
\cdot  f = \sum_{\underline{i} \leq \underline{r}}    \underline{\partial}^{[ \underline{i}]} (f) \cdot  \underline{\partial}^{[\underline{r} - \underline{i}]}.
\end{align}
 Then,  $\mcD_X^{(\M)}$ is locally generated as an $\mcO_X$-algebra  by the sections $\partial_i^{[p^l]}  \left(= \partial_i^{\langle p^l\rangle}\right)$ for various $l \leq \M$.
$\mcD_X^{(0)}$ is, in particular,  locally generated by $\partial_1  \left(:= \partial_1^{[p^0]} \right), \cdots, \partial_n  \left(:= \partial_n^{[ p^0 ]} \right)$ and  $\mcD_X^{(\infty)}$ is Grothendieck's sheaf of differential operators (cf. ~\cite[Section 16.8]{EGA4}).
For each pair of nonnegative integers $(\M_1, \M_2)$ with $\M_1 \leq \M_2$,
there exists  a natural $k$-linear morphism $\mcD_X^{(\M_1)} \migi \mcD_X^{(\M_2)}$  given by $\underline{\partial}^{[\underline{r}]} \mapsto \underline{\partial}^{[\underline{r}]}$, which restricts to an isomorphism $\mcD_{X, l}^{(\M_1)} \isom \mcD_{X, l}^{(\M_2)}$ for each $l < p^{\M_1+1}$.
The sheaves $\mcD_X^{(\M)}$  ($m \geq 0$) and the morphisms $\mcD_X^{(\M_1)} \migi \mcD_X^{(\M_2)}$ form 
 an inductive system with  
\begin{align} \label{E1004}
\varinjlim_{\M \geq 0} \mcD_X^{(\M)} \cong  \mcD_X^{(\infty)}.
\end{align}

By  a {\bf (left) $\mcD_X^{(\M)}$-module}, we shall mean  a pair $(\mcV, \nabla)$ consisting of 
a vector bundle (i.e., a locally free coherent sheaf of finite rank)
$\mcV$ on $X$ and an $\mcO_X$-linear morphism $\nabla : {^L \mcD}_X^{(\M)} \migi \mcE nd_k (\mcV)$ of sheaves of $k$-algebras, i.e., a left $\mcD_{X}^{(\M)}$-action $\nabla$ on $\mcV$ extending its $\mcO_X$-module structure, where $\mcE nd_k (\mcV)$ denotes the sheaf of locally defined $k$-linear endomorphisms of $\mcV$ endowed with a structure of $\mcO_X$-module given by left multiplication. 
An {\bf invertible $\mcD_X^{(\M)}$-module} is defined to be a $\mcD_X^{(\M)}$-module $(\mcL, \nabla_\mcL)$ such that $\mcL$ is a line bundle.

We see that giving  a $\mcD_X^{(0)}$-module amounts to giving
a flat bundle on $X$ in the sense of ~\cite[Definition 4.5]{Wak5}  (i.e., a vector bundle on $X$ together with  a flat connection).
As for a general $\M$,
giving a $\mcD_X^{(\M)}$-module structure on an $\mcO_X$-module is equivalent to giving an $\M$-PD stratification on that sheaf (cf. ~\cite[Proposition 2.3.2]{Be1}).

Given two $\mcD_X^{(\M)}$-modules $(\mcV_i, \nabla_{i})$ ($i=1,2$), we define  a {\bf morphism} from  $(\mcV_1, \nabla_{1})$ to $(\mcV_2, \nabla_{2})$ as an $\mcO_X$-linear morphism $\mcV_1 \migi \mcV_2$ that is compatible with  the respective $\mcD_X^{(m)}$-actions  $\nabla_{1}$, $\nabla_{2}$.
Denote by $\nabla_{\mcO_X}^{\mr{triv}(m)}$ the $\mcD_X^{(m)}$-action on $\mcO_X$  given in \eqref{aaaa}.
Thus, we obtain the trivial  (invertible) $\mcD_{X}^{(m)}$-module
\begin{align} \label{Eh1063}
(\mcO_X, \nabla_{\mcO_X}^{\mr{triv}(m)}).
\end{align}

Given a vector bundle $\mcV$ on $X$,
  we always suppose that the tensor product  $\mcD_{X, l}^{(\M)} \otimes \mcV := {^R \mcD}_{X, l}^{(\M)} \otimes \mcV$ (resp., $\mcV \otimes \mcD_{X, l}^{(\M)}:= \mcV \otimes {^L \mcD}_{X, l}^{(\M)}$) is endowed with 
a structure of $\mcO_X$-module
 arising  from
 the left (resp., right) $\mcO_X$-module structure on $\mcD_{X, l}^{(\M)}$.

Let $\M'$ be another element of $\mbZ_{\geq 0} \sqcup \{\infty\}$  with $\M \leq \M'$ and let $(\mcV, \nabla)$ be a $\mcD_X^{(\M')}$-module.
Then, the composite 
$\nabla |^{(\M)} : {^L \mcD}_{X}^{(\M)} \migi {^l \mcD}_{X}^{(\M')} \xrightarrow{\nabla} \mcE nd_k (\mcV)$
 specifies a $\mcD_X^{(\M)}$-module structure  on $\mcV$. That is to say, we obtain a $\mcD_X^{(\M)}$-module 
\begin{align}\label{Ew23467}
(\mcV, \nabla |^{(\M)}).
\end{align}
We  shall refer to $(\mcV, \nabla |^{(\M)})$ as the {\bf $\M$-th truncation} of $(\mcV, \nabla)$.

If $(\mcV_i, \nabla_{i})$ ($i=1,2$) are $\mcD_{X}^{(\M)}$-modules, 
then
the tensor product $\mcV_1 \otimes \mcV_2$ has a canonical  $\mcD_X^{(\M)}$-module structure 
\begin{align} \label{Ep1}
\nabla_{1} \otimes \nabla_{2} : {^L \mcD}_X^{(m)} \migi \mcE nd_k (\mcV_1 \otimes \mcV_2)
\end{align} 
determined, under the local description discussed above,  by
\begin{align}
\partial^{[r]}_i (v_1 \otimes v_2) = \sum_{j=0}^r \partial_i^{[j]} (v_1) \otimes \partial^{[r-j]}_i (v_2)
\end{align} 
for  any $r <p^{m+1}$, $i=1, 2, \cdots, n$ and 
any local sections $v_1 \in \mcV_1$, $v_2 \in \mcV_2$.
In a similar way, we can define the dual $(\mcV^\vee, \nabla^\vee)$ of a $\mcD_X^{(\M)}$-module $(\mcV, \nabla)$ (cf. ~\cite[Corollary 2.3.3]{Be1} for both constructions).

\LSP
\subsection{Kodaira-Spencer map associated to a $\mcD^{(\M)}$-module} \label{SS018c}

 
\SSP
\bde  \label{D0199}
 Let
 $(\mcV, \nabla)$ be a $\mcD_X^{(\M)}$-module and $\mcU$ a subbundle of $\mcV$.
  Then,  the {\bf Kodaira-Spencer map} associated to the triple  $(\mcV, \nabla, \mcU)$  is defined as  
  the $\mcO_X$-linear morphism
  \begin{align} \label{Eh1010}
  \mr{KS}_{(\mcV, \nabla, \mcU)} : \mcT_X \migi \mcU^\vee \otimes (\mcV/\mcU)  \left(= \mcH om_{\mcO_X} (\mcU, \mcV/\mcU) \right)
  \end{align}
  induced, in the natural manner, from 
   the $\mcO_X$-linear composite
  \begin{align} \label{Eh1011}
  \mr{KS}'_{(\mcV, \nabla, \mcU)} : 
  \mcU \xrightarrow{\mr{inclusion}} \mcV \xrightarrow{\breve{\nabla}^{(0)}} \Omega_X \otimes \mcV \xrightarrow{\mr{quotient}} \Omega_{X} \otimes (\mcV/ \mcU),
  \end{align}
   where $\breve{\nabla} |^{(0)}$ denotes the flat connection on $\mcV$
 corresponding to the $\mcD_X^{(0)}$-action $\nabla |^{(0)}$.
 Notice that if  $\mr{KS}_{(\mcV, \nabla, \mcU)}$ is an isomorphism and $\mcU$ is a line bundle, then (since $\mcT_X$ forms a rank $n$ vector bundle) $\mcV$ must be a vector bundle of rank $n+1$.
  \ede
\SSP

\begin{rema}[Various Kodaira-Spencer maps] \label{R13401} 
In this remark, we compare two notions of Kodaira-Spencer map defined
above and in ~\cite[Definition 4.1.1]{Wak12}.

Let $[\infty]$ be the point $[1:0:0:\cdots :0]$ in $\mbP^n$, and 
 write
 \begin{align} \label{Ep101}
 \PP
 \end{align}
  for the subgroup of $\mr{PGL}_{n+1}$ consisting of elements  fixing $[\infty]$.
 For a $\mr{PGL}_{n+1}$-bundle $\mcE$ on  $X$, one verifies that
 giving a $\PP$-reduction of $\mcE$ is equivalent to giving a global section $X \migi \mbP^n_\mcE$ of the $\mbP^n$-bundle $\mbP^n_\mcE$ corresponding to $\mcE$. 
Also, write
$\mr{GL}_{n+1}^\infty := \mr{GL}_{n+1} \times_{\mr{PGL}_{n+1}} \PP  \left(\subseteq \mr{GL}_{n+1} \right)$, and
 \begin{align} \label{Ep104}
 \RR:= 
 \PP
  \cap 
 \QQ
    \left(\subseteq \mr{PGL}_{n+1} \right).
 \end{align}
The assignment  given by  $\overline{A} \mapsto \overline{A} ([\infty])$
for any $\overline{A} \in \mr{PGL}_{n+1}$ defines an isomorphism
\begin{align} \label{Ep105}
\mr{PGL}_{n+1} / \PP \isom \mbP^n  \ \left(\text{resp.,} \ \QQ / \RR \isom \mbA^n\right).
\end{align}

Next,  let $(\mcV, \nabla)$ be a $\mcD_X^{(0)}$-module such that $\mcV$ is of rank $n+1$, and 
  $\mcU$ a line subbundle of $\mcV$.
 Denote by $(\mcE, \nabla_\mcE)$  the flat $\mr{GL}_{n+1}$-bundle on $X$  (in the sense of ~\cite[Definition 1.28]{Wak5} or ~\cite[Section 7.1]{Wak12})
 corresponding to $(\mcV, \nabla)$.
 Then $\mcU$ specifies a $\mr{GL}_{n+1}^\infty$-reduction $\mcE_\mr{red}$ on $\mcE$.
Also, denote by $(\overline{\mcE}, \nabla_{\overline{\mcE}})$ the flat $\mr{PGL}_{n+1}$-bundle induced by $(\mcE, \nabla_\mcE)$, and by $\overline{\mcE}_\mr{red}$ the $\PP$-reduction on $\overline{\mcE}$ determined by $\mcE_\mr{red}$.

Observe that there exists a canonical isomorphism $\mcA d_{\mcE} \isom \mcE nd_{\mcO_X}(\mcV)$ (where $\mcA d_{(-)}$ denotes the adjoint vector bundle associated to a principal bundle $(-)$) that  restricts to an isomorphism from $\mcA d_{\mcE_\mr{red}}$ to the $\mcO_X$-submodule of $\mcE nd_{\mcO_X}(\mcV)$ consisting of locally defined $\mcO_X$-linear endomorphisms $f$ with $f (\mcU) \subseteq \mcU$.
Hence, we have a commutative diagram consisting  of isomorphisms
 \begin{align}  \label{Ew3232}
\vcenter{\xymatrix@C=36pt@R=36pt{
\mcA d_{\overline{\mcE}}/ \mcA d_{\overline{\mcE}_\mr{red}}  \ar[d]_-{\wr} & \mcA d_{\mcE}/ \mcA d_{\mcE_\mr{red}} \ar[l]_-{\sim} \ar[r]^-{\sim} \ar[d]^-{\wr}&\mcH om_{\mcO_X} (\mcU, \mcV/\mcU)
\\
\widetilde{\mcT}_{\overline{\mcE}}/\widetilde{\mcT}_{\overline{\mcE}_\mr{red}}
&\widetilde{\mcT}_{\mcE}/\widetilde{\mcT}_{\mcE_\mr{red}}
\ar[l]^-{\sim}&
}}
\end{align}
such that  the both sides of vertical arrows are the isomorphism displayed preceding 
~\cite[Example 4.1.3]{Wak12} in the case where the pair ``$(\mcE_G, \mcE_H)$" is taken to be $(\mcE, \mcE_\mr{red})$ and  $(\overline{\mcE}, \overline{\mcE}_\mr{red})$, respectively.
By  this diagram,    both  the Kodaira-Spencer maps $\mr{KS}_{(\mcE, \nabla_\mcE, \mcE_\mr{red})}$ and $\mr{KS}_{(\overline{\mcE}, \nabla_{\overline{\mcE}}, \overline{\mcE}_\mr{red})}$ defined in the manner  of ~\cite[Definition 4.1.1]{Wak12}  can be identified with 
  $\mcO_X$-linear morphisms $\mcT_X \migi \mcH om_{\mcO_X} (\mcU, \mcV/\mcU)$.
Under this identification, the morphisms  $\mr{KS}_{(\mcE, \nabla_\mcE, \mcE_\mr{red})}$, $\mr{KS}_{(\overline{\mcE}, \nabla_{\overline{\mcE}}, \overline{\mcE}_\mr{red})}$  
   coincide with $\mr{KS}_{(\mcV, \nabla, \mcN)}$ at least up to composition with an automorphism of the vector bundle $\mcH om_{\mcO_X}(\mcU, \mcV/\mcU)$.
 In particular, $\mr{KS}_{(\mcV, \nabla, \mcN)}$ is an isomorphism if and only if  $\mr{KS}_{(\mcE, \nabla_\mcE, \mcE_\mr{red})}$, or equivalently  $\mr{KS}_{(\overline{\mcE}, \nabla_{\overline{\mcE}}, \overline{\mcE}_\mr{red})}$,  is an isomorphism.
 \end{rema}

\LSP
\subsection{Indigenous and affine-indigenous  $\mcD^{(\M)}$-modules} \label{SS018}

We shall introduce   
indigenous and affine-indigenous $\mcD^{(\M)}_X$-modules defined on a variety of arbitrary positive  dimension.

\SSP
\bde  \label{D01234}
 \begin{itemize}
 \item[(i)]
 An  {\bf indigenous $\mcD_X^{(\M)}$-module}
   is a triple
 \begin{align} \label{Eh99}
 \mcV^\diamondsuit := (\mcV, \nabla, \mcN)
 \end{align}
 consisting of  a $\mcD_{X}^{(\M)}$-module $(\mcV, \nabla)$ and a line subbundle $\mcN$ of $\mcV$ such that
  the Kodaira-Spencer map $\mr{KS}_{(\mcV, \nabla, \mcN)} : \mcT_X \migi \mcH om_{\mcO_X} (\mcN, \mcV/\mcN)$ is an isomorphism.

  \item[(ii)]
 Let $\mcV_i^\diamondsuit := (\mcV_i, \nabla_{i}, \mcN_i)$ ($i=1,2$) be  two indigenous
 $\mcD_X^{(\M)}$-modules.
 Then, an {\bf isomorphism} from $\mcV_1^\diamondsuit$ to $\mcV_2^\diamondsuit$
 is defined as an isomorphism  of vector bundles
 $\eta : \mcV_1 \isom \mcV_2$  that is compatible with the respective $\mcD_{X}^{(\M)}$-actions $\nabla_{1}$, $\nabla_{2}$ and satisfies the equality
   $\eta (\mcN_1) = \mcN_2$.
\end{itemize}
  \ede

\SSP
\bde  \label{D01235}
 \begin{itemize}
 \item[(i)]
 An {\bf affine-indigenous  $\mcD_X^{(\M)}$-module}
  is a quadruple 
 \begin{align}
{^\A \mcV}^\diamondsuit :=  (\mcV, \nabla, \mcN, \delta)
 \end{align}
consisting of an indigenous $\mcD_X^{(\M)}$-module  $(\mcV, \nabla, \mcN)$ 
  and  
a left inverse morphism $\delta : \mcV \migisurj \mcN$ to 
the natural inclusion $\mcN \migiincl \mcV$  (i.e., the composite $\mcN \migiincl \mcV \xrightarrow{\delta} \mcN$ coincides with the identity morphism of $\mcN$)
such that $\mr{Ker}(\delta)  \left(\subseteq \mcV \right)$ is closed under the $\mcD_X^{(\M)}$-action $\nabla$.

 \item[(ii)]
 Let ${^\A\mcV}_i^\diamondsuit := (\mcV_i, \nabla_{i}, \mcN_i, \delta_i)$ ($i=1,2$) be 
  two affine-indigenous $\mcD_X^{(\M)}$-modules.
Then, an {\bf isomorphism} from ${^\A \mcV}_1^\diamondsuit$ to ${^\A \mcV}_2^\diamondsuit$
 is defined as an isomorphism of 
 indigenous $\mcD_X^{(\M)}$-modules $\eta : (\mcV_1, \nabla_{1}, \mcN_1) \isom (\mcV_2, \nabla_{2}, \mcN_2)$   that is compatible  with $\delta_1$ and $\delta_2$, i.e., satisfies the equality $\delta_2 \circ \eta = \eta |_{\mcN_1} \circ \delta_1$.
 \end{itemize}
  \ede
\SSP

We shall prove several  properties on  indigenous $\mcD_X^{(\M)}$-modules.

\SSP
\bpr \label{P36h9}
Let $(\mcV, \nabla, \mcN)$ be a triple consisting of a $\mcD_X^{(\M)}$-module $(\mcV, \nabla)$ and a line subbundle $\mcN$ of $\mcV$.
Then, the triple $(\mcV, \nabla, \mcN)$ forms an indigenous
$\mcD_X^{(\M)}$-module
 if and only if the following $\mcO_X$-linear  composite is an isomorphism:
 \begin{align} \label{Eh101}
\mr{KS}^\circledast_{(\mcV, \nabla, \mcN)} :  \mcD_{X,  1}^{(\M)} \otimes \mcN \xrightarrow{\mr{inclusion}} \mcD_{X}^{(\M)} \otimes \mcV \xrightarrow{\nabla}\mcV.
 \end{align}
  \epr
\begin{proof}
Recall that the natural morphism  $\mcD_{X, 1}^{(0)} \rightarrow \mcD_{X, 1}^{(\M)}$ is an isomorphism.
Under the identification $\mcD_{X, 1}^{(0)} \otimes \mcN = \mcD_{X, 1}^{(\M)} \otimes \mcN$ arising from  this isomorphism, 
$\mr{KS}^\circledast_{(\mcV, \nabla, \mcN)}$ coincides with
the composite
\begin{align} \label{Eq33gr}
 \mcD_{X,  1}^{(0)} \otimes \mcN \xrightarrow{\mr{inclusion}} \mcD_{X}^{(\M)} \otimes \mcV \xrightarrow{\nabla |^{(0)}}\mcV.
\end{align}
This composite restricts to the identity morphism of $\mcN$ via the inclusions $ \mcN \left(= \mcD_{X, 0}^{(0)}  \otimes \mcN \right)$
$\hookrightarrow \mcD_{X, 1}^{(0)} \otimes \mcN$ and $\mcN \hookrightarrow \mcV$, so we obtain the following morphism of short exact sequences:
\begin{align} \label{Ejwe4}
\vcenter{\xymatrix@C=21pt@R=36pt{
0 \ar[r] & \mcN \left(= \mcD_{X, 0}^{(0)} \otimes \mcN \right)\ar[r] \ar[d]^-{\mr{id}_\mcN} & \mcD_{X, 1}^{(0)} \otimes \mcN \ar[r] \ar[d]^-{\eqref{Eq33gr}} & \mcT_X \otimes \mcN  \left(= ( \mcD_{X, 1}^{(0)}/\mcD_{X, 0}^{(0)}) \otimes \mcN \right)\ar[r] \ar[d] & 0
\\
0 \ar[r] & \mcN \ar[r] & \mcV \ar[r] & \mcV/\mcN \ar[r] & 0. 
}}
\end{align}
By this diagram,  we see that
$\mr{KS}^\circledast_{(\mcV, \nabla, \mcN)}$  is an isomorphism if and only if the right-hand vertical arrow $\mcT_X \otimes \mcN \rightarrow \mcV/\mcN$ is an isomorphism.
Hence, the assertion follows from the fact that this morphism $\mcT_X \otimes \mcN \rightarrow \mcV/\mcN$ corresponds to $\mr{KS}_{(\mcV, \nabla, \mcN)}$ via a canonical bijection $\mr{Hom}_{\mcO_X} (\mcT_X \otimes \mcN, \mcV/\mcN) \cong  \mr{Hom} (\mcT_X, \mcN^\vee \otimes (\mcV/\mcN))$. 
\end{proof}

\SSP
\bpr \label{P3fg6h9}  
 Let $\mcV^\diamondsuit_i := (\mcV_i, \nabla_{i}, \mcN_i)$ ($i=1,2$) be indigenous $\mcD_X^{(\M)}$-modules.
 Then, the following assertions hold:
\begin{itemize}
\item[(i)] 
Denote by $\mr{Isom}(\mcV_1^\diamondsuit, \mcV_2^\diamondsuit)$ (resp., $\mr{Isom}(\mcN_1, \mcN_2)$) the set of isomorphisms of indigenous
$\mcD_X^{(\M)}$-modules
  (resp., $\mcO_X$-modules)  from $\mcV_1^\diamondsuit$ (resp., $\mcN_1$) to $\mcV_2^\diamondsuit$ (resp., $\mcN_2$). 
Then, the assignment $\eta \mapsto \eta |_{\mcN_1}$ gives an injection of sets
\begin{align} \label{Eh1120}
\mr{Isom}(\mcV_1^\diamondsuit, \mcV_2^\diamondsuit) \migiincl 
\mr{Isom}(\mcN_1, \mcN_2).
\end{align}
\item[(ii)]
Suppose that $X$ is proper over $k$, or more generally, the equality $\Gamma (X, \mcO_X) =k$ holds.
Suppose further  that 
$\mcV_1 = \mcV_2$ and $\mcN_1 = \mcN_2$.
Then, $\nabla_{1} = \nabla_{2}$ if and only if $\mcV_1^\diamondsuit \cong  \mcV_2^\diamondsuit$.
Moreover, there exist canonical identifications 
\begin{align}
\mr{Isom}(\mcV_1^\diamondsuit, \mcV_2^\diamondsuit) =k^\times = \mr{Isom}(\mcN_1, \mcN_2)
\end{align}
 and the map  \eqref{Eh1120} coincides with the identity map of $k^\times$ under these identifications.
\end{itemize}
  \epr
\begin{proof}
First, we shall consider assertion (i).
Let $\eta$
   be an isomorphism $\mcV_1^\diamondsuit \isom \mcV_2^\diamondsuit$.
Since $\eta$ preserves the respective $\mcD_X^{(\M)}$-actions $\nabla_{1}$ and $\nabla_{2}$,
the following square diagram is commutative:
\begin{align} \label{Eh456}
\vcenter{\xymatrix@C=56pt@R=36pt{
\mcD_{X, 1}^{(\M)} \otimes \mcN_1\ar[r]^-{\mr{id} \otimes (\eta |_{\mcN_1})}_-{\sim} \ar[d]_-{\mr{KS}^\circledast_{\mcV_1^\diamondsuit}}& \mcD_{X, 1}^{(\M)} \otimes \mcN_2 \ar[d]^-{\mr{KS}^\circledast_{\mcV_2^\diamondsuit}}
\\
\mcV_1 \ar[r]^-{\sim}_{\eta} &\mcV_2.
}}
\end{align}
By Proposition \ref{P36h9}, both the right-hand and left-hand vertical arrows are isomorphisms.
Hence, the commutativity of this diagram implies that 
$\eta$ is uniquely determined by its  restriction $\eta |_{\mcN_1}$ to $\mcN_1$.
This completes the proof of assertion (i).

Next, we shall consider the former assertion of (ii).
The ``only if'' part is clear.
To prove the ``if'' part, we suppose that there exists an isomorphism $\eta : \mcV_1^\diamondsuit \isom \mcV_2^\diamondsuit$.
Because of the assumption on $X$, i.e, $\Gamma (X, \mcO_X) =k$,
the automorphism $\eta |_{\mcN}$ of the line bundle $\mcN := \mcN_1 = \mcN_2$ may be given by multiplication by some $a \in k^\times$. 
By the commutativity of the  diagram \eqref{Eh456}, 
$\eta$ turns out to be  given by multiplication by $a$.
But,  such an  automorphism cannot  preserve  the respective $\mcD_X^{(\M)}$-actions $\nabla_{1}$, $\nabla_{2}$ unless $\nabla_{1} = \nabla_{2}$.
This proves the former assertion of  (ii).

 Finally, the latter assertion of (ii) follows immediately from the argument just discussed.
\end{proof}
\SSP

\begin{rema}[Determinant] \label{R1f3401}
Note that the sheaf $\mcD_{X, 1}^{(\M)}$ forms a rank $n+1$ vector bundle   obtained as an extension of $\mcD_{X, 1}^{(\M)}/ \mcD_{X, 0}^{(\M)}$ $\left( = \mcT_X\right)$  by $\mcD_{X, 0}^{(\M)}   \left(= \mcO_X\right)$.
Hence, we obtain natural isomorphisms
\begin{align} \label{Ep30}
\mr{det} (\mcD_{X,1}^{(\M)}) \isom \mr{det}(\mcD_{X,0}^{(\M)}) \otimes \mr{det}(\mcD_{X, 1}^{(\M)}/ \mcD_{X, 0}^{(\M)})\isom \mr{det} (\mcT_X) \isom \omega_X^\vee.
\end{align} 
If  $\mcV^\diamondsuit := (\mcV, \nabla, \mcN)$ is an indigenous 
$\mcD_X^{(\M)}$-module, then
the isomorphism $\mr{KS}^\circledast_{(\mcV, \nabla, \mcN)} : \mcD_{X, 1}^{(\M)}\otimes \mcN \isom \mcV$ (cf. \eqref{Eh101}) induces, via taking determinants,   the following composite isomorphism
\begin{align} \label{Eh401}
\mr{det}(\mcV) & \isom \mr{det}(\mcD_{X, 1}^{(\M)}\otimes \mcN)  \isom
\mr{det}(\mcD_{X, 1}^{(\M)}) \otimes \mcN^{\otimes (n+1)}  \isom \omega^\vee_X \otimes \mcN^{\otimes (n+1)}.
\end{align}
In particular,  the determinant $\mr{det}(\mcV)$ is isomorphic to the line bundle
$\omega^\vee_X \otimes \mcN^{\otimes (n+1)}$. 

Also, let ${^\mbA \mcV}^\diamondsuit := (\mcV, \nabla, \mcN, \delta)$ be an affine-indigenous $\mcD_X^{(\M)}$-module.
By passing through $\mr{KS}^\circledast_{(\mcV, \nabla, \mcN)}$,
we can identify the inclusion $\mcN \migiincl \mcV$ with the natural inclusion $\mcN  \left(= \mcD_{X, 0}^{(\M)} \otimes \mcN \right)$ $\migiincl \mcD_{X, 1}^{(\M)} \otimes \mcN$.
This implies that the surjection $\delta : \mcV \migisurj \mcN$ allows us to identify $\mr{Ker}(\delta)$ with $(\mcD_{X, 1}^{(\M)} \otimes \mcN) / (\mcD_{X, 0}^{(\M)} \otimes \mcN)$.
Hence, we have the following composite isomorphism
\begin{align} \label{Eh405}
\mr{Ker}(\delta) \isom  (\mcD_{X, 1}^{(\M)} \otimes \mcN) / (\mcD_{X, 0}^{(\M)} \otimes \mcN) \isom  (\mcD_{X, 1}^{(\M)}/\mcD_{X, 0}^{(\M)}) \otimes \mcN \isom \mcT_X \otimes \mcN,
\end{align}
which induces an isomorphism
\begin{align} \label{Eh406}
\mr{det}(\mr{Ker}(\delta)) \isom \left(\mr{det}(\mcT_X \otimes \mcN) = \right)  \omega_X^\vee \otimes \mcN^{\otimes n}.
\end{align}
 \end{rema}

\LSP
\subsection{Twisting by invertible $\mcD^{(\M)}$-modules} \label{SS018a}

Let $\mcV^\diamondsuit := (\mcV, \nabla, \mcN)$
 (resp., ${^\A\mcV}^\diamondsuit := (\mcV, \nabla, \mcN, \delta)$) be an indigenous  (resp., affine-indigenous)
 $\mcD_X^{(\M)}$-module
   and $(\mcL, \nabla_\mcL)$ an invertible $\mcD_X^{(\M)}$-module.
Since $\mr{KS}_{(\mcV, \nabla, \mcN)}$ is an isomorphism, the Kodaira-Spencer map $\mr{KS}_{(\mcV \otimes \mcL, \nabla\otimes \nabla_\mcL, \mcN \otimes \mcL)}$ associated to the triple 
$(\mcV \otimes \mcL, \nabla \otimes \nabla_\mcL, \mcN \otimes \mcL)$   (cf. \eqref{Ep1} for the definition of $\nabla \otimes \nabla_\mcL$) is an isomorphism.
In particular, the collection
\begin{align} \label{Ep2}
\mcV^\diamondsuit_{\otimes (\mcL, \nabla_\mcL)} := (\mcV \otimes \mcL, \nabla_\mcV \otimes \nabla, \mcN \otimes \mcL) \hspace{17mm}\\
\left(\text{resp.,} \  {^\A \mcV}^\diamondsuit_{\otimes (\mcL, \nabla_\mcL)} := (\mcV \otimes \mcL, \nabla \otimes \nabla_\mcL, \mcN \otimes \mcL, \delta \otimes \mr{id}_{\mcL})\right) \notag
\end{align}
forms an indigenous (resp., affine-indigenous) $\mcD_X^{(\M)}$-module.

Regarding the underlying vector bundles of  indigenous and affine-indigenous $\mcD_X^{(\M)}$-modules, we  prove the following assertion.

\SSP
\ble \label{L067} 
 For each $i=1,2$, suppose that we are given $\M_i \in \mbZ_{\geq 0} \sqcup \{ \infty \}$ and  an indigenous $\mcD_X^{(\M_i)}$-module 
 $\mcV^\diamondsuit_i := (\mcV_i, \nabla_{i}, \mcN_i)$ 
 (resp., an affine-indigenous  $\mcD_X^{(\M_i)}$-module
 ${^\A \mcV}^\diamondsuit_i := (\mcV_i, \nabla_{i}, \mcN_i, \delta_i)$).
 In the non-affine case, 
  we suppose further that $p \nmid (n+1)$.
Then, there exist a line bundle $\mcL$ on $X$ and an $\mcO_X$-linear  isomorphism
$\eta : \mcV_1 \isom \mcV_2 \otimes \mcL$ with $\eta (\mcN_1) = \mcN_2 \otimes \mcL$.
(We can choose $\mcL$ as the line bundle $\mcN_1 \otimes \mcN_2^\vee$.)
 \ele
\begin{proof}
By replacing $\mcV^\diamondsuit_i$ (resp., ${^\A \mcV}^\diamondsuit_i$) with 
$(\mcV_i, \nabla_{i}|^{(0)}, \mcN_i)$ (resp., $(\mcV_i, \nabla_{i}|^{(0)}, \mcN_i, \delta_i)$), we may assume that $\M_1 = \M_2 =0$.

First, we shall consider the assertion for the  non-affine case.
It follows from Proposition \ref{P36h9} that, for each $i=1,2$,  the vector bundle $\mcV_i$  may be identified with  $\mcD_{X, 1}^{(0)} \otimes \mcN_i$
 via  
$\mr{KS}^\circledast_{(\mcV_i, \nabla_{i}, \mcN_i)}$. 
Denote by $\breve{\nabla}_{i}^{\mr{det}}$ 
 the connection on $\omega_X^\vee \otimes \mcN_i^{\otimes (n+1)}$  induced, after taking determinants,  from $\nabla_{i}$  via
 the isomorphism $\mr{det}(\mcV_i) \isom \omega_X^\vee \otimes \mcN_i^{\otimes (n+1)}$ (cf. \eqref{Eh401}).
  Then, tensoring $\breve{\nabla}_{1}^{\mr{det}}$ with the dual of $\breve{\nabla}_{2}^{\mr{det}}$ yields a connection $\breve{\nabla}$ on the line bundle $(\mcN_1 \otimes \mcN_2^\vee)^{\otimes (n+1)}  \left(\cong (\omega_X^\vee \otimes \mcN_i^{\otimes (n+1)}) \otimes (\omega_X^\vee \otimes \mcN_i^{\otimes (n+1)})^\vee \right)$.
Because of the assumption $p \nmid n+1$, one can find uniquely  a  connection $\breve{\nabla}^{\frac{1}{n+1}}$ on $\mcL := \mcN_1 \otimes \mcN_2^\vee$ such that the $(n+1)$-st tensor product
of $\breve{\nabla}^{\frac{1}{n+1}}$ coincides with $\breve{\nabla}$.
In fact, if $\breve{\nabla}$ may be described locally as $\breve{\nabla} = d + A$ for some $A \in \Gamma (X, \Omega_X)$ under a local trivialization $\mcL \cong \mcO_X$, then we can write $\breve{\nabla}^{\frac{1}{n+1}}= d + \frac{1}{n+1}A$.
The tensor product 
 of $\breve{\nabla}^{\frac{1}{n+1}}$ and  $\nabla_{2}$ (more precisely, the connection corresponding to $\nabla_{2}$) defines  a  $\mcD_X^{(0)}$-action $\nabla_{\mcV \otimes \mcL, 2}$ on 
$\mcV_2 \otimes \mcL$.
Since $\mcV^\diamondsuit_2$ forms an indigenous bundle,
the Kodaira-Spencer map 
\begin{align} \label{Eh2037}
\eta:=\mr{KS}^\circledast_{(\mcV_2 \otimes \mcL, \nabla_{\mcV \otimes \mcL, 2}, \mcN_2 \otimes \mcL)} :  \mcD_{X, 1}^{(0)} \otimes (\mcN_2 \otimes \mcL) \migi \mcV_2 \otimes \mcL
\end{align}
(cf. \eqref{Eh101})
 is an isomorphism satisfying that $\eta (\mcD_{X, 0}^{(0)} \otimes (\mcN_2 \otimes \mcL)) = \mcN_2 \otimes \mcL$.
Using  the natural identifications
\begin{align}
\mcV_1   \left(= \mcD_{X, 1}^{(0)} \otimes \mcN_1\right) \cong   \mcD_{X, 1}^{(0)} \otimes (\mcN_2 \otimes \mcL), \hspace{5mm} \mcN_1   \left(= \mcD_{X, 0}^{(0)} \otimes \mcN_1\right) \cong   \mcD_{X, 0}^{(0)} \otimes (\mcN_2 \otimes \mcL),
\end{align}
 $\eta$ specifies the required isomorphism.
This completes the proof of the non-affine case.

Next, we shall consider the assertion for the affine case.
For each $i =1,2$, denote by $\nabla_{\mcN, i}$ the $\mcD_X^{(0)}$-action on $\mcN_i$ induced from $\nabla_{i}$ via the surjection $\delta_i$.
Also, denote by $\nabla_\mcL$ the  structure of $\mcD_X^{(0)}$-action  on $\mcL := \mcN_1 \otimes \mcN_2^\vee$ defined as the tensor product of $\nabla_{\mcN,1}$ and the dual of $\nabla_{\mcN,2}$.
Just as \eqref{Eh2037},  the Kodaira-Spencer map associated to 
the underlying  indigenous $\mcD_X^{(0)}$-module of ${^\A \mcV}^\diamondsuit_{\otimes (\mcL, \nabla_\mcL)}$ specifies an isomorphism $\eta: \mcV_1 \isom \mcV_2 \otimes \mcL$ with $\eta (\mcN_1) = \mcN_2 \otimes \mcL$.
This completes the proof of assertion (ii).
\end{proof}

\LSP
\subsection{$p$-curvature and $F$-divided sheaves} \label{SS058}

Suppose that $\M \neq \infty$.
Denote by $\mcK_{X}^{(\M)}$ the kernel of  $\nabla_{\mcO_X}^{\mr{triv}(\M)} : {^L\mcD}_X^{(\M)} \migi \mcE nd_k (\mcO_X)$ (cf. \eqref{Eh1063}).
If $(\mcV, \nabla)$ is a $\mcD_X^{(\M)}$-module, then the composite
\begin{align}
{^p \psi}_{(\mcV, \nabla)}^{(\M)} : \mcK_X^{(\M)} \migiincl {^L \mcD}_X^{(\M)} \xrightarrow{\nabla} \mcE nd_k (\mcV)
\end{align}
is called the {\bf $p$-$\M$-curvature} of $(\mcV, \nabla)$ (cf. ~\cite[Definition 3.1.1]{LSQ}).

Given  a $\mcD_{X}^{(\M)}$-module $(\mcV, \nabla)$,
we shall write 
\begin{align} \label{Ew90}
\mcS ol(\nabla)
\end{align}
for the subsheaf of $\mcV$ 
on which $\mcD_{X, +}^{(\M)}$ acts as zero, where $\mcD_{X, +}^{(\M)}$ denotes the kernel of the canonical projection $\mcD_X^{\M} \migi \mcO_X$.
We shall refer to $\mcS ol(\nabla)$ as
 the {\bf sheaf of horizontal sections} of $(\mcV, \nabla)$.
This sheaf  may be thought of as an $\mcO_{X^{(\M +1)}}$-module via the underlying homeomorphism of $F_{X/k}^{(\M +1)}$.

Let $\mcU$ be  a vector bundle on $X^{(\M +1)}$.
Then,  its pull-back $F_{X/k}^{(\M +1)*}(\mcU)$  via $F_{X/k}^{(\M +1)}$ admits naturally a  ${\mcD}_X^{(\M)}$-action
\begin{align} \label{Ew9090}
\nabla^{\mr{can}(\M)}_{\mcU} : {^L \mcD}_X^{(\M)}  \migi \mcE nd_k (F_{X/k}^{(\M +1)*}(\mcU)), 
\end{align}
for which 
the pair  $(F_{X/k}^{(\M +1)*}(\mcU), \nabla^{\mr{can}(\M)}_{\mcU})$ forms a $\mcD_X^{(\M)}$-module with vanishing $p$-$\M$-curvature.
 The ${\mcD}_X^{(\M)}$-action  $\nabla^{\mr{can}(\M)}_{\mcU}$ may be given, under the local description  of  $\mcD_X^{(\M)}$ introduced in Section \ref{SS018b}, by 
 $\underline{\partial}^{[\underline{r}]} (x \otimes v) := \underline{\partial}^{[\underline{r}]} (x) \otimes v$
  for each local section $x \otimes v \in \mcO_X \otimes_{(F_{X/k}^{(\M +1)})^{-1}(\mcO_{X^{(\M +1)}})} (F_{X/k}^{(\M +1)})^{-1}(\mcU)  \left(=  F_{X/k}^{(\M +1)*}(\mcU)\right)$.

According to  ~\cite[Corollary 3.2.4]{LSQ}, 
 the assignments $\mcU \mapsto (F_{X/k}^{(\M +1)*}(\mcU), \nabla^{\mr{can}(\M)}_{\mcU})$ and $(\mcV, \nabla) \mapsto \mcS ol(\nabla)$ determine 
an equivalence of categories
\begin{align} \label{E3002}
\begin{pmatrix}
\text{the category of} \\
\text{vector bundles on $X^{(\M +1)}$} 
\end{pmatrix}
\isom \begin{pmatrix}
\text{the category of $\mcD_X^{(\M)}$-modules} \\
\text{with vanishing $p$-$\M$-curvature} \\
\end{pmatrix}.
\end{align}

Recall (cf. ~\cite[Section 2.2]{dS}) that
an {\bf $F$-divided sheaf} on $X$ is a collection
\begin{align}
\{ (\mcE_\M, \alpha_\M)\}_{\M\in \mbZ_{\geq  0}},
\end{align}
where each  pair $(\mcE_\M, \alpha_\M)$  consists of a vector bundle  $\mcE_\M$ on $X^{(\M)}$ 
 and
  an $\mcO_{X^{(\M)}}$-linear isomorphism $\alpha_\M : F^{(1)*}_{X^{(\M)}/k}(\mcE_{\M+1}) \isom \mcE_\M$.
(Note that $F$-divided sheaves are the same as {\it flat bundles} in the sense of ~\cite[Definition 1.1]{G}.)
$F$-divided sheaves form a category, in which the morphisms from $\{ (\mcE_\M, \alpha_\M) \}_\M$ to $\{ (\mcE'_\M, \alpha'_\M) \}_{\M}$ are the collections of $\mcO_{X^{(\M)}}$-linear morphisms $\mcE_\M\migi \mcE'_\M$ which commute with the $\alpha_\M$'s and $\alpha'_\M$'s in the obvious sense.

For an  $F$-divided sheaf $\{ (\mcE_\M, \alpha_\M) \}_{\M}$ on $X$,
$\mcE_0$ admits a canonical   
 $\mcD_X^{(\infty)}$-action.
In fact, the $\mcD_X^{(\M)}$-action $\nabla_{\mcE_{\M+1}}^{\mr{can}(\M)}$ (cf. \eqref{Ew9090}) on $F_{X/k}^{(\M+1) *}(\mcE_{\M+1})$ induces 
a $\mcD_{X}^{(\M)}$-action 
on $\mcE_0$ (with vanishing $p$-$\M$-curvature) via the composite isomorphism 
\begin{align}
F^{(\M)*}_{X/k}(\alpha_{\M}) \circ \cdots\circ F^{(1)*}_{X/k}(\alpha_1)\circ\alpha_0 : F_{X/k}^{(\M+1) *}(\mcE_{m+1}) \isom \mcE_0.
\end{align}
The collection of these actions for various $\M$ yields, via  \eqref{E1004},  a $\mcD_X^{(\infty)}$-action
\begin{align}
\nabla^\mr{can (\infty)}_{\mcE_0} : {^L \mcD}_X^{(\infty)} \migi \mcE nd_k (\mcE_0)
\end{align}
 on $\mcE_0$.
By Katz's theorem appearing in ~\cite[Theorem 1.3]{G}, the assignment $\{ (\mcE_\M, \alpha_\M) \}_{\M} \mapsto (\mcE_0, \nabla^{\mr{can}(\infty)}_{\mcE_0})$ determines an equivalence of categories:  
\begin{align} \label{E3001}
\begin{pmatrix}
\text{the category of} \\
\text{$F$-divided sheaves on $X$} 
\end{pmatrix}
\isom \begin{pmatrix}
\text{the category of } \\
\text{$\mcD_X^{(\infty)}$-modules} \\
\end{pmatrix}.
\end{align}
Moreover, for each $\N \geq 1$, the following diagram of categories is $2$-commutative:
\begin{align}  \label{Ew3216}
\begin{CD}
\begin{pmatrix}
\text{the category of} \\
\text{$F$-divided sheaves on $X$} 
\end{pmatrix}
@>\eqref{E3001}> \sim > \begin{pmatrix}
\text{the category of } \\
\text{$\mcD_X^{(\infty)}$-modules} \\
\end{pmatrix}
\\
@V \{ (\mcE_\M, \alpha_\M)\}_\M \mapsto \mcE_{\N}VV @VV (\mcV, \nabla) \mapsto (\mcV, \nabla |^{(\N-1)}) V
\\
\begin{pmatrix}
\text{the category of} \\
\text{vector bundles on $X^{(\N)}$} 
\end{pmatrix} @> \sim >\eqref{E3002} > \begin{pmatrix}
\text{the category of $\mcD_X^{(\N-1)}$-modules} \\
\text{with vanishing $p$-$(\N -1)$-curvature} 
\end{pmatrix}.
\end{CD}
\end{align}
For convenience, we will say that any $\mcD_X^{(\infty)}$-module {\it has vanishing $p$-$\infty$-curvature}.

\LSP
\subsection{Dormant indigenous  and affine-indigenous $\mcD^{(\M)}$-modules} \label{SS258}

\SSP
\bde  \label{Dbb030}
 Let 
  $\mcV^\diamondsuit := (\mcV, \nabla, \mcN)$ (resp., ${^\mbA \mcV}^\diamondsuit := (\mcV, \nabla, \mcN, \delta)$)  be an indigenous  (resp.,  affine-indigenous) 
 $\mcD_X^{(\M)}$-module.
 Then, we shall say that $\mcV^\diamondsuit$ (resp., ${^\mbA \mcV}^\diamondsuit$) is  {\bf dormant}  if the $p$-$\M$-curvature of $(\mcV, \nabla)$ vanishes identically.
 In the case of infinite level,
  we shall say that any  indigenous  (resp., affine-indigenous) 
  $\mcD_X^{(\infty)}$-module
  is {\bf dormant}, for convenience.
  \ede
\SSP

The following proposition will be useful when we
try to understand which varieties have an indigenous or affine-indigenous  $\mcD_X^{(\infty)}$-module.

\SSP
\bpr \label{P3kjff7} 
Let $\mcN$ be a line bundle on $X$.
Then, the following assertions hold:
\begin{itemize}
\item[(i)]
Suppose that there exists a dormant indigenous
$\mcD_X^{(\infty)}$-module
 of the form $\mcV^\diamondsuit := (\mcV, \nabla, \mcN)$.
Then, the order of the line bundle $\omega^\vee_X \otimes \mcN^{\otimes (n+1)}$
in the Neron-Severi group of $X$ is finite of order prime to $p$.
In particular, $\omega^\vee_X \otimes \mcN^{\otimes (n+1)}$
 is numerically trivial.
\item[(ii)]
Suppose that there exists a dormant affine-indigenous
$\mcD_X^{(\infty)}$-module
  of the form ${^\mbA\mcV}^\diamondsuit := (\mcV, \nabla, \mcN, \delta)$.
 Then, both $\omega_X$ and $\mcN$ are numerically trivial.
 \end{itemize}
 \epr
\begin{proof}
 Assertion (i)  follows from ~\cite[Theorem 1.8]{G}.
 In fact,
 by passing to 
 the isomorphism $\mr{det}(\mcV) \isom \omega^\vee \otimes \mcN^{\otimes (n+1)}$ (cf. \eqref{Eh401}),
  we can obtain a structure of $\mcD_X^{(\infty)}$-action $\nabla_{\omega^\vee \otimes \mcN^{\otimes (n+1)}}$ on the line bundle $\omega^\vee \otimes \mcN^{\otimes (n+1)}$ corresponding to  the determinant of $\nabla_\mcV$.

Next, we shall consider assertion (ii).
Since $\mr{Ker}(\delta)$ is closed under the $\mcD_X^{(\infty)}$-action  
$\nabla$,
$\mcN$ has   a structure of $\mcD_{X}^{(\infty)}$-action $\nabla_{\mcN}$ induced by $\nabla$ via the quotient $\delta : \mcV \migisurj \mcN$.
Let $\nabla_{\mcN^{\otimes (n+1)}}$ denote the $\mcD_{X}^{(\infty)}$-action on $\mcN^{\otimes (n+1)}$ defined as the $(n+1)$-st tensor product of  $\nabla_{\mcN}$.
Then, tensoring $\nabla_{\mcN^{\otimes (n+1)}}$ and the dual of  $\nabla_{\omega^\vee \otimes \mcN^{\otimes (n+1)}}$ obtained in the proof of (i)
yields a $\mcD_X^{(\infty)}$-action $\nabla_{\omega_X}$ on $\omega_X$.
Thus, the proof of the assertion is completed by applying again   ~\cite[Theorem 1.8]{G} to both $(\mcN, \nabla_\mcN)$ and $(\omega_X, \nabla_{\omega_X})$.
\end{proof}
\SSP

For each $\M \in \mbZ_{\geq 0} \sqcup \{\infty\}$,  denote by
\begin{align} \label{Eq22290}
\DID{\M}{, X}
  \ \left(\text{resp.,} \ 
 \DAID{\M}{, X}
    \right), \ \text{or simply} \
    \DID{\M}{}
  \ \left(\text{resp.,} \ 
 \DAID{\M}{}
    \right),
\end{align}
the set of isomorphism classes of dormant indigenous  (resp., dormant affine-indigenous) $\mcD_X^{(\M)}$-modules.

\vspace{10mm}
\section{Chern class formula for $F^\N$-projective/affine structures} \label{Sg02gg10}\SSP

In this section, we
establish  a bijective correspondence between $F^\N$-projective  (resp., $F^\N$-affine) structures and 
dormant  indigenous  (resp., dormant  affine-indigenous) $\mcD_X^{(\N-1)}$-modules.
This correspondence
yields   necessary conditions on Chern classes for the existence of  $F^\N$-projective or $F^\N$-affine structures, respectively  (cf. Theorem \ref{T016}).

\LSP
\subsection{$F^\N$-projective structures associated  to dormant indigenous $\mcD_X^{(\N -1)}$-modules}
 \label{SS0DDF7}
 Let $\N$ be a positive integer and $\mcV^\diamondsuit := (\mcV, \nabla, \mcN)$  a dormant indigenous
$\mcD_X^{(\N-1)}$-module.

Since $(\mcV, \nabla)$ has vanishing $p$-$(\N-1)$-curvature,
the sheaf of horizontal sections $\mcS ol (\nabla)$ 
forms a rank $(n+1)$ vector bundle
 on $X^{(\N)}$ and the   morphism $F_{X/k}^{(\N)*}(\mcS ol (\nabla)) \migi \mcV$ extending  the $\mcO_{X^{(\N)}}$-linear  inclusion $\mcS ol (\nabla) \migiincl \mcV$
  is an isomorphism (cf. \eqref{E3002}).
  If $\mcE^{}$ (resp.,  $\mcE^\nabla$) denotes the $\mr{PGL}_{n+1}$-bundle on $X$ (resp., $X^{(\N)}$) induced by $\mcV$ (resp., $\mcS ol (\nabla)$)
   via projectivization,
  then this isomorphism yields an isomorphism of $\mr{PGL}_{n+1}$-bundles   $\tau : F_{X/k}^{(\N)*}(\mcE^{\nabla}) \isom \mcE$.
According to ~\cite[Proposition A.7.3]{Wak12}, 
there exists a canonical   $(\N -1)$-PD stratification   $\varepsilon^\blacklozenge$ (with vanishing $p$-$(\N -1)$-curvature) on the Frobenius pull-back $F_{X/k}^{(\N)*}(\mcE^\nabla)$;
we use the same notation ``$\varepsilon^\blacklozenge$" to denote the corresponding  $(\N -1)$-PD stratification on $\mcE$ via $\tau$. 

On the other hand, the line subbundle $\mcN  \left(\subseteq \mcV \right)$ determines  a $\PP$-reduction $\mcE_\mr{red}$ of $\mcE$.  
Since $\mr{KS}_{\mcV^\diamondsuit}$ is an isomorphism,  the Kodaira-Spencer map associated to $(\mcE, \varepsilon^\blacklozenge, \mcE_\mr{red})$
 (in the sense of ~\cite[Definition 4.1.1]{Wak12}) turns out to be an isomorphism  (cf.
Remark \ref{R13401}).
As a consequence,  the pair 
$(\mcE_\mr{red}, \varepsilon^\blacklozenge)$
 forms 
 a dormant indigenous $(\mr{PGL}_{n+1}, \mr{PGL}_{n+1}^\infty)$-bundle of level $\N$ (cf. ~\cite[Definitions 4.5.1 and 4.5.2]{Wak12}).
 It follows from ~\cite[Theorem B]{Wak12} that this pair corresponds to a unique  $F^\N$-projective structure  on $X$ (cf. ~\cite[Example 2.3.3, (i)]{Wak12}), which we denote by 
 \begin{align} \label{EQrri}
 \mcV^{\diamondsuit \Rightarrow \heartsuit}.
 \end{align}

Next, 
suppose further that we are given a left inverse morphism $\delta: \mcV \migi \mcN$ to the inclusion $\mcN \migiincl \mcV$ for which the quadruple 
 ${^\mbA \mcV}^\diamondsuit := (\mcV, \nabla, \mcN, \delta)$ forms a dormant  affine-indigenous
$\mcD_X^{(\N-1)}$-module.
The $\mcD_X^{(\N-1)}$-action $\nabla$ restricts to 
a $\mcD_X^{(\N-1)}$-action $\nabla_{\mr{Ker}(\delta)}$ on $\mr{Ker}(\delta)$ with vanishing $p$-$(\N-1)$-curvature.
Then, the sheaf of horizontal sections $\mcS ol(\nabla_{\mr{Ker}(\delta)})$ 
 determines a $\QQ$-reduction $\mcE^{\nabla \mbA}$ of  $\mcE^{\nabla}$ via projectivization, by which $\varepsilon^{\blacklozenge}$ is reduced to
 an $(\N -1)$-PD stratification $\varepsilon^{\mbA\blacklozenge}$  on $\mcE^\mbA := F_{X/k}^{(\N)*}(\mcE^{\nabla \mbA})$.
Since $\RR = \PP \times_{\mr{PGL}_{n+1}}\QQ$,  the intersection of the two reductions of $\mcE$, i.e.,  
$\mcE^\mbA$
 and $\mcE_\mr{red}$,  
specifies a $\RR$-reduction $\mcE_\mr{red}^{\mbA}$.
The  pair  $(\mcE_\mr{red}^{\mbA},  \varepsilon^{\mbA \blacklozenge})$ forms a dormant indigenous $(\mr{PGL}_{n+1}^\mbA, \mr{PGL}_{n+1}^{\mbA, \infty})$-bundle of level $\N$.
By ~\cite[Theorem B]{Wak12} again,
it corresponds to a unique $F^\N$-affine structure on $X$ (cf. ~\cite[Example 2.3.3, (ii)]{Wak12}), which  will be denoted by
\begin{align}
{^\mbA}\mcV^{\diamondsuit \Rightarrow \heartsuit}.
\end{align}

The resulting  assignment $\mcV^\diamondsuit \mapsto \mcV^{\diamondsuit \Rightarrow \heartsuit}$ (resp., ${^\A\mcV}^\diamondsuit \mapsto {^\A\mcV}^{\diamondsuit \Rightarrow \heartsuit}$) 
determines a map of sets
\begin{align} \label{Eq2304}
 \zeta_\N^{\diamondsuit \Rightarrow \heartsuit} : 
 \DID{\N -1}{}
 \migi 
 \Proj{\N}{}
  \ \left(\text{resp.,} \ 
  \xi_\N^{\diamondsuit \Rightarrow \heartsuit} : 
  \DAID{\N -1}{}
   \migi 
   \Aff{\N}{}
    \right),
 \end{align}
which 
 is compatible with the formation of truncation to lower levels.
Hence, by considering this assignment for every $\N \in \mbZ_{>0}$, 
we also obtain an $F^\infty$-projective (resp., $F^\infty$-affine) structure $\mcV^{\diamondsuit \Rightarrow \heartsuit}$ (resp., ${^\A\mcV}^{\diamondsuit \Rightarrow \heartsuit}$) associated to  a dormant indigenous (resp., dormant affine-indigenous) $\mcD_X^{(\infty)}$-module $\mcV^\diamondsuit$ (resp., ${^\A}\mcV^\diamondsuit$).

In any case of  $\N \in \mbZ_{>0} \sqcup \{\infty\}$, we shall refer to 
 $\mcV^{\diamondsuit \Rightarrow \heartsuit}$ (resp., ${^\A\mcV}^{\diamondsuit \Rightarrow \heartsuit}$)
 as the {\bf $F^\N$-projective} (resp., {\bf $F^\N$-affine}) {\bf structure associated to  $\mcV^\diamondsuit$} (resp., {\bf ${^\A}\mcV^\diamondsuit$}).

\LSP
\subsection{Correspondence  via projectivization I}
 \label{SS0ffw6h7}

We shall introduce certain  equivalence relations in 
$\DID{\N -1}{}$
and
$\DAID{\N -1}{}$.

\SSP
\bde  \label{D0f30}
 For each $i=1,2$,
 let $\mcV^\diamondsuit_i := (\mcV_i, \nabla_{i}, \mcN_i)$ (resp., ${^\A \mcV}^\diamondsuit_i := (\mcV_i, \nabla_{i}, \mcN_i, \delta_i)$)  be a dormant indigenous (resp., dormant affine-indigenous) 
 $\mcD_X^{(\N-1)}$-module.
 We shall say that
 these
   are {\bf $\mbG_m$-equivalent} if there exists an invertible  $\mcD_X^{(\N-1)}$-module $(\mcL, \nabla_\mcL)$  with vanishing $p$-$(\N-1)$-curvature 
   satisfying 
  $\mcV^\diamondsuit_{1 \otimes (\mcL, \nabla_\mcL)} \cong \mcV_2^\diamondsuit$ 
 (resp., ${^\mbA }\mcV^\diamondsuit_{1 \otimes (\mcL, \nabla_\mcL)} \cong {^\A \mcV}_2^\diamondsuit$) (cf. \eqref{Ep2}).
  \ede
\SSP

The relation of  {\it being $\mbG_m$-equivalent}
in fact
 defines an equivalence relation in 
$\DID{\N -1}{}$
 (resp., 
 $\DAID{\N -1}{}$).
Write
  \begin{align} \label{Ef1113}
  \eDID{\N -1}{, X}
 \ \left(\text{resp.,} \ 
 \eDAID{\N -1}{, X}
   \right), \ \text{or simply} \ 
    \eDID{\N -1}{}
 \ \left(\text{resp.,} \ 
 \eDAID{\N -1}{}
   \right),
\end{align}
for the quotient set of 
$\DID{\N -1}{}$
 (resp., 
 $\DAID{\N -1}{}$)
  by this equivalence relation.
If  $\mcV^\diamondsuit$  (resp.,  ${^\A \mcV}^\diamondsuit$) is a dormant indigenous  (resp., dormant affine-indigenous) 
$\mcD_X^{(\N-1)}$-module,
 then
we denote by $[\mcV^\diamondsuit]$
 (resp., $[{^\A\mcV}^\diamondsuit]$)
  the image of $\mcV^\diamondsuit$ in
  $\eDID{\N -1}{}$
    (resp., the image of ${^\A\mcV}^\diamondsuit$ in 
    $\eDAID{\N -1}{}$).
If $(\mcL, \nabla_\mcL)$ is an invertible  $\mcD_{X}^{(\N-1)}$-module with vanishing $p$-$(\N-1)$-curvature,
then, since $\mcV^{\diamondsuit \Rightarrow \heartsuit}$ (resp., ${^\mbA} \mcV^{\diamondsuit \Rightarrow \heartsuit}$) is constructed via   change  of structure group by the quotient $\mr{GL}_{n+1} \twoheadrightarrow \mr{PGL}_{n+1}$,    we have  
$\zeta_\N^{\diamondsuit \Rightarrow \heartsuit}(\mcV^{\diamondsuit}) = \zeta_\N^{\diamondsuit \Rightarrow \heartsuit}(\mcV^{\diamondsuit}_{\otimes (\mcL, \nabla_\mcL)})$
 (resp., 
 $\xi_\N^{\diamondsuit \Rightarrow \heartsuit}({^\A}\mcV^{\diamondsuit}) = \xi_\N^{\diamondsuit \Rightarrow \heartsuit}({^\A}\mcV^{\diamondsuit}_{\otimes (\mcL, \nabla_\mcL)})$) (cf. \eqref{Ep2}).
It follows  that
the 
assignment
  $[\mcV^\diamondsuit] \mapsto \zeta_\N^{\diamondsuit \Rightarrow \heartsuit} (\mcV^\diamondsuit)$
 (resp., 
  $[{^\A}\mcV^\diamondsuit] \mapsto \xi_\N^{\diamondsuit \Rightarrow \heartsuit} ({^\A}\mcV^\diamondsuit)$)
  determines a  well-defined map
\begin{align} \label{Ep34}
 \overline{\zeta}_\N^{\diamondsuit \Rightarrow \heartsuit} : 
 \eDID{\N -1}{}\!
\migi \Proj{\N}{}
\ \left(\text{resp.,} \ 
\overline{\xi}_\N^{\diamondsuit \Rightarrow \heartsuit} : 
\eDAID{\N -1}{}\!
\migi 
\Aff{\N}{}\right).
\end{align}

\LSP
\subsection{Correspondence via projectivization II}
 \label{SS06h7}

In this subsection, we prove   the bijectivities of  $\overline{\zeta}_\N^{\diamondsuit \Rightarrow \heartsuit}$ 
and $\overline{\xi}_\N^{\diamondsuit \Rightarrow \heartsuit}$ for $\N \neq \infty$ (cf. Theorem \ref{P01111}).

Until the end of Lemma \ref{L020}, {\it we assume that $X$ is quasi-projective over $k$}.
For each smooth affine algebraic group $G$ over $k$,
we denote by $\check{H}^1_{\text{\'{e}t}} (X, G)$ the first
\v{C}ech  cohomology set of $X$ with coefficients in $G$ with respect to the \'{e}tale topology.
This pointed set may be identified with the set of isomorphism classes of $G$-bundles on $X$, where the distinguished point is represented by  the trivial $G$-bundle.
Since we assumed that $X$ is quasi-projective,  there exists 
a coboundary map 
\begin{align}
\gamma : \check{H}^1_{\text{\'{e}t}}(X, \mr{PGL}_{n+1}) \migi H^2_{\text{\'{e}t}} (X, \mbG_m)
\end{align}
   induced by  the short exact sequence
 \begin{align} \label{Ew22345}
 1 \longmigi \mbG_m \longmigi \mr{GL}_{n+1} \longmigi \mr{PGL}_{n+1} \longmigi 1
 \end{align}
(cf. ~\cite[Chapter III, Theorem 2.17, and Chapter IV, Theorem 2.5]{Mil}).
Given a $\mr{PGL}_{n+1}$-bundle 
 $\mcE$ on $X$, we shall denote by $[\mcE]$ the element of $\check{H}^1_{\text{\'{e}t}}(X, \mr{PGL}_{n+1})$ 
  representing $\mcE$.

\SSP
\ble \label{L019} 
 The coboundary map 
 \begin{align}
 \check{H}^1_{\text{\'{e}t}} (X, \PP) \migi H^2_{\text{\'{e}t}} (X, \mbG_m)
  \end{align}
induced
  by 
the short exact sequence
\begin{align} \label{E4001}
 1 \longmigi \mbG_m \longmigi \mr{GL}_{n+1}^\infty \longmigi \PP \longmigi 1
 \end{align}
 is the zero morphism.
 \ele
\begin{proof}
Let
$\pi_1 : \mr{GL}_{n+1}^\infty \migi \mbG_m$,  $\pi_2 : \mr{GL}_{n+1}^\infty \migi \mr{GL}_n$,
and $\pi_3 : \PP \migi \mr{GL}_n$
be the homomorphisms given, respectively,  by 
\begin{align}
\pi_1 (\begin{pmatrix} a & {\bf a} \\ {\bf 0} & A \end{pmatrix}) = a, \hspace{10mm} \pi_2 (\begin{pmatrix} a & {\bf a} \\ {\bf 0} & A \end{pmatrix}) = A,  \hspace{10mm} \pi_3 (\overline{\begin{pmatrix} a & {\bf a} \\ {\bf 0} & A \end{pmatrix}}) = a^{-1} \cdot A
\end{align}
for any $a \in \mbG_m$, ${\bf a} \in \mbA^n$, and $A \in \mr{GL}_n$.
These homomorphisms fit into the following morphism of short exact sequences: 
\begin{align} \label{E0010}
\begin{CD}
1 @>>> \mbG_m@> \varDelta>> \mr{GL}_{n+1}^\infty  @>\mr{projection}>> \PP @>>>1
\\
@. @VV\mr{id}V @VV \pi_1 \times \pi_2 V @VV\pi_3 V @.
\\
1 @>>> \mbG_m @>>a \mapsto (a,  \varDelta (a))> \mbG_m \times \mr{GL}_{n}@>>(a, A) \mapsto a^{-1} \cdot A> \mr{GL}_n @>>>1,
\end{CD}
\end{align}
where the upper horizontal sequence is \eqref{E4001} and $\varDelta$'s denote the diagonal embeddings.
It induces a commutative square diagram
\begin{align} \label{E0011}
\begin{CD}
\check{H}^1_{\text{\'{e}t}} (X, \PP) @>\gamma_1>> H^2_{\text{\'{e}t}} (X, \mbG_m) 
\\
@V \check{H}^1_{\text{\'{e}t}} (\pi_3)VV @VV \mr{id} V
\\
\check{H}^1_{\text{\'{e}t}} (X, \mr{GL}_n) @>>\gamma_2> H^2_{\text{\'{e}t}}(X, \mbG_m).
\end{CD}
\end{align}
Observe here that the lower horizontal sequence in \eqref{E0010} splits by taking the first  projection $\mbG_m \times \mr{GL}_n \migisurj \mbG_m$ as a split surjection.
Hence, $\gamma_2 = 0$, which implies  $\gamma_1  \left(= \gamma_2 \circ \check{H}^1_{\text{\'{e}t}} (\pi_3) \right)=0$.
This completes the proof of the lemma.
\end{proof}

\SSP
\ble \label{L020} 
Let $\N$ be a positive integer and $\mcE^\nabla$ a $\mr{PGL}_{n+1}$-bundle
  on $X^{(\N)}$ with $\gamma ([F_{X/k}^{(\N)*}(\mcE^\nabla)]) =0$.
Suppose further that $p \nmid (n+1)$.
Then, there exists a $\mr{GL}_{n+1}$-bundle 
 on $X^{(\N)}$ inducing 
$\mcE^\nabla$ via projectivization.
  \ele
\begin{proof}
First, we shall prove  the following claim: 
\begin{quotation}
{\it For each positive integer $\N$,  the endomorphism $f_{\mbG_m}$ of $\check{H}^2_{\text{\'{e}t}} (X, \mbG_m)$ induced  by the $\N$-th relative Frobenius morphism  $F_{\mbG_m/k}^{(N)}$ of $\mbG_m$
    restricts to a bijective endomap of  the subset  $\mr{Im}(\gamma) \left(\subseteq \check{H}^2_{\text{\'{e}t}} (X, \mbG_m) \right)$.}
\end{quotation}
Let us consider the following morphism of short exact sequences:
\begin{align} \label{E4007}
\begin{CD}
1 @>>> \mu_{n+1} @>>> \mr{SL}_{n+1} @>>> \mr{PSL}_{n+1} @>>>1
\\
@. @VV\mr{inclusion} V @VV \mr{inclusion} V @VV\mr{inclusion} V @.
\\
1 @>>> \mbG_m @>>> \mr{GL}_{n+1} @>>> \mr{PGL}_{n+1} @>>>1,
\end{CD}
\end{align}
where $\mu_{n+1}$ denotes the group of $(n+1)$-st roots of unity.
The right-hand vertical arrow $\mr{PSL}_{n+1} \migi \mr{PGL}_{n+1}$ is an isomorphism of algebraic groups  because of the assumption $p \nmid n+1$.
This implies that $\gamma$ can be  decomposed as
\begin{align}
 \check{H}^1_{\text{\'{e}t}}(X, \mr{PGL}_{n+1}) \left(\cong  \check{H}^1_{\text{\'{e}t}}(X, \mr{PSL}_{n+1}) \right) \xrightarrow{\gamma_1}  H^2_{\text{\'{e}t}}(X, \mu_{n+1}) \xrightarrow{\gamma_2}  H^2_{\text{\'{e}t}} (X, \mbG_m),
\end{align}
where $\gamma_1$ denotes the coboundary map induced by the upper horizontal arrow in \eqref{E4007} and $\gamma_2$ denotes the map induced by the natural inclusion  $\mu_{n+1}\migiincl \mbG_m$.
In particular, the inclusion relation  $\mr{Im}(\gamma) \subseteq \mr{Im}(\gamma_2)$ holds.
Since $H^2_{\text{\'{e}t}}(X, \mu_{n+1})$  is finite (cf. ~\cite[Chapter III, Corollary 2.10, and Chapter VI,  Corollary 5.5]{Mil}), 
 the set $\mr{Im}(\gamma_2)$, as well as $\mr{Im}(\gamma)$,  turns out to be finite.
Next, note that the diagram \eqref{E4007} is compatible with the $N$-th relative Frobenius morphisms  on the various algebraic groups appearing there.
It follows that $f_{\mbG_m}(\mr{Im}(\gamma)) \subseteq \mr{Im}(\gamma)$.
Moreover, if $f_{\mu_{n+1}}$ denotes the endomorphism of  $H^2_{\text{\'{e}t}}(X, \mu_{n+1})$ induced by the $\N$-th relative Frobenius morphism  $F^{(\N)}_{\mu_{n+1}/k}$ of $\mu_{n+1}$, then
  the surjection $H^2_{\text{\'{e}t}}(X, \mu_{n+1}) \migisurj \mr{Im}(\gamma_2)$ is compatible with 
  $f_{\mu_{n+1}}$ and $f_{\mbG_m} |_{\mr{Im}(\gamma_2)}$.
Since $F^{(\N)}_{\mu_{n+1}/k}$   is an isomorphism,   $f_{\mu_{n+1}}$  is bijective.
This implies from the finiteness of $\mr{Im}(\gamma_2)$  that $f_{\mbG_m} |_{\mr{Im}(\gamma_2)}$ is bijective.
Thus, the endomap $f_{\mbG_m} |_{\mr{Im}(\gamma)}$ of $\mr{Im}(\gamma)$ is bijective, so
this completes the proof of the claim.

Now, let us go back to the   proof of  the lemma.
Let $\mcE^\nabla$ be a $\mr{PGL}_{n+1}$-bundle on $X^{(\N)}$ with $\gamma ([F_{X/k}^{(\N)*}(\mcE^\nabla)]) =0$.
If $\mcE^{\nabla F}$ denotes the  $\mr{PGL}_{n+1}$-bundle on $X$ corresponding to $\mcE^\nabla$ via $W_{X/k} : X^{(\N)}\isom X$,
then $f_{\mbG_m} (\gamma ([\mcE^\nabla])) = \gamma ([F_{X/k}^{(\N)*}(\mcE^\nabla)]) =0$.
It follows from the above claim that 
the equality $\gamma (\mcE^{\nabla F}) =0$ holds.
 Since the sequence
 \begin{align}
 \check{H}^1_{\text{\'{e}t}}(X, \mr{GL}_{n+1}) \longmigi 
 \check{H}^1_{\text{\'{e}t}}(X, \mr{PGL}_{n+1}) \stackrel{\gamma}{\longmigi}
 H^2_{\text{\'{e}t}}(X, \mbG_m)
 \end{align}
  induced by \eqref{Ew22345} is exact (cf. ~\cite[Chapter IV, Theorem 2.5, Step 3]{Mil}), 
 there exists a $\mr{GL}_{n+1}$-bundle  $\mcF^F$ 
 on $X$ inducing  $\mcE^{\nabla F}$  via projectivization.
Hence, the $\mr{GL}_{n+1}$-bundle on $X^{(\N)}$ corresponding to this  $\mr{GL}_{n+1}$-bundle via $W_{X/k}$ satisfies  the  required conditions.
This completes the proof of the assertion.
\end{proof}
\SSP

\bt \label{P01111} 
Suppose that $X$ is quasi-projective over $k$.
 Then, the following assertions hold:
\begin{itemize}
\item[(i)]
Suppose that $p \nmid (n+1)$.
 Then, for each $\N \in \mbZ_{>0}$,  the map $\overline{\zeta}_\N^{\diamondsuit \Rightarrow \heartsuit}$
    is 
  bijective.
  \item[(ii)]
 For each $\N \in \mbZ_{>0}$,  the map $\overline{\xi}_\N^{\diamondsuit \Rightarrow \heartsuit}$
    is bijective.
   \end{itemize}
  \et
\begin{proof}
 First, we shall consider  the surjectivity of  $\overline{\zeta}_\N^{\diamondsuit \Rightarrow \heartsuit}$  asserted in  (i).
 Let $\mcS^\heartsuit$ be an $F^\N$-projective structure on $X$.
 Since $(\mr{PGL}_{n+1})_X^{(\N)}$ is affine, the {\it right} $(\mr{PGL}_{n+1})_X^{(\N)}$-torsor defined as the inverse of $\mcS^\heartsuit$ specifies a  $\mr{PGL}_{n+1}$-bundle $\mcE^\nabla$ on $X^{(\N)}$ via the underlying homeomorphism of $F_{X/k}^{(\N)}$.
 The graphs $\Gamma_\phi : U \migi U \times \mbP^n$ associated to various local sections $\phi : U \migi \mbP^n$ in $\mcS^\heartsuit$ may be glued together to construct a $\mr{PGL}_{n+1}^{\infty}$-reduction $\mcE_\mr{red}$ of $F_{X/k}^{(\N)*}(\mcE^\nabla)$.
This means  that  the element $[F_{X/k}^{(\N)*}(\mcE^\nabla)]$ of $\check{H}^1_{\text{\'{e}t}} (X, \mr{PGL}_{n+1})$ represented by  $F_{X/k}^{(\N)*}(\mcE^\nabla)$
    comes from $\check{H}^1_{\text{\'{e}t}} (X, \PP)$.
   Hence, Lemma \ref{L019} implies 
    $\gamma ([F_{X/k}^{(\N)*}(\mcE^\nabla)])  =0$. 
 It follows from Lemma \ref{L020} that there exists a vector bundle 
 $\mcV^\nabla$ on $X^{(\N)}$ inducing  $\mcE^\nabla$  via projectivization.
Write $\mcV := F^{(\N)*}_{X/k}(\mcV^\nabla)$.
According to ~\cite[Theorem B]{Wak12} and the discussion in Remark \ref{R13401},   
 the triple 
$\mcV^{\diamondsuit} := (\mcV, \nabla_{\mcV^\nabla}^{\mr{can}(\N-1)}, \mcN)$ forms a dormant indigenous 
$\mcD_X^{(\N-1)}$-module, where $\mcN$ denotes the line bundle determined by the $\PP$-reduction $\mcE_\mr{red}$.
One verifies that 
$\mcV^{\diamondsuit \Rightarrow \heartsuit} \cong \mcS^\heartsuit$ by construction, and 
this completes the proof of the surjectivity of
$\overline{\zeta}_\N^{\diamondsuit \Rightarrow \heartsuit}$.

Next, we shall prove the injectivity of  $\overline{\zeta}_\N^{\diamondsuit \Rightarrow \heartsuit}$.
Let $\mcV_i^\diamondsuit := (\mcV_i, \nabla_{i}, \mcN_i)$ ($i=1,2$) be  dormant indigenous
$\mcD_X^{(\N-1)}$-modules, and
 suppose that there exists  an isomorphism
$\alpha : \mcV_1^{\diamondsuit \Rightarrow \heartsuit} \isom \mcV_2^{\diamondsuit \Rightarrow \heartsuit}$.
Then, $\alpha$ induces an isomorphism $\mcE^\nabla_1 \isom \mcE_2^\nabla$ between the $\mr{PGL}_{n+1}$-bundles  associated to  $\mcS ol (\nabla_{1})$ and $\mcS ol (\nabla_{2})$, respectively.
It follows (cf. ~\cite[Chapter II, Exercise 7.9]{Har}) that we can find a line bundle $\mcL$ on $X^{(\N)}$ and an isomorphism 
$\mcS ol (\nabla_{1}) \isom  \mcL \otimes \mcS ol (\nabla_{2})$ inducing this isomorphism of $\mr{PGL}_{n+1}$-bundles.
In order to complete the proof, we can assume (after  twisting $\mcV_2^\diamondsuit$ by $(F_{X/k}^{(\N)*}(\mcL), \nabla^{\mr{can}(\N-1)}_{\mcL})$) that there exists an $\mcO_{X^{(\N)}}$-linear  isomorphism $\alpha^\nabla_\mcV : \mcS ol (\nabla_{1}) \isom \mcS ol (\nabla_{2})$ inducing $\alpha$ via projectivization.
Let $\alpha_\mcV : (\mcV_1, \nabla_{1}) \isom (\mcV_2, \nabla_{2})$ denote  the isomorphism of $\mcD_{X}^{(\N-1)}$-modules corresponding to $\alpha_\mcV^\nabla$ via the equivalence of categories \eqref{E3002}.
Since  $\alpha$
 preserves the respective $\PP$-reductions
 in $\mcE_1^\nabla$ and $\mcE_2^\nabla$ (defined as in the above discussion),
   the equality $\alpha_\mcV (\mcN_1) = \mcN_2$ holds.
Thus, $\alpha_\mcV$ specifies an isomorphism $\mcV_1^\diamondsuit \isom \mcV_2^\diamondsuit$ of indigenous $\mcD_X^{(\N-1)}$-modules.
This completes the proof of the injectivity of $\overline{\zeta}_\N^{\diamondsuit \Rightarrow \heartsuit}$, and consequently, completes the proof of  assertion (i).

Next, we shall prove assertion (ii).
 Let us consider 
 the surjectivity of $\overline{\xi}_\N^{\diamondsuit \Rightarrow \heartsuit}$.
Write $\mr{GL}_{n+1}^\A := \mr{GL}_{n+1} \times_{\mr{PGL}_{n+1}} \QQ \left(\subseteq  \mr{GL}_{n+1}\right)$.
This algebraic group fits into the following natural short exact sequence:
\begin{align}
1 \longmigi \mbG_m \longmigi \mr{GL}_{n+1}^\A \longmigi \mr{PGL}_{n+1}^\A \longmigi 1.
\end{align}
The homomorphism $\mr{PGL}_{n+1}^\A \migi \mr{GL}_{n+1}^\A$ given by $\overline{\begin{pmatrix} a & {\bf 0} \\ {^t{\bf a}} & A\end{pmatrix}} \mapsto \begin{pmatrix} 1 & {\bf 0} \\ a^{-1} \cdot {^t{\bf a}} & a^{-1} \cdot A\end{pmatrix}$ specifies a split injection $\sigma$ of this sequence with $\sigma (\RR) \subseteq \mr{GL}_{n+1}^{\A, \infty}$, where $\mr{GL}_{n+1}^{\A, \infty} := \mr{GL}_{n+1} \times_{\mr{PGL}_{n+1}} \RR$.
Now, let
$\mcS^\heartsuit$ be an $F^\N$-affine structure on $X$.
The associated $\QQ$-bundle $\mcE^\nabla$ (constructed as in the case of $F^\N$-projective structures)  determines, via 
 change of structure group by $\sigma$, a vector bundle $\mcV^\nabla$ on $X^{(\N)}$ together with a surjection $\delta^\nabla : \mcV^\nabla \migisurj \mcN^\nabla$ onto a line bundle $\mcN^\nabla$.
 Here, we shall write $\mcV := F^{(\N)*}_{X/k}(\mcV^\nabla)$ and $\delta := F^{(\N)*}_{X/k}(\delta^\nabla)$.
The  $\mr{GL}_{n+1}^{\A, \infty}$-reduction  $\mcE_\mr{red} \times^{\RR, \sigma} \mr{GL}_{n+1}^{\A, \infty}$ of $F^{(\N)*}_{X/k}(\mcE^\nabla) \times^{\QQ, \sigma} \mr{GL}_{n+1}^\A$ corresponds to a line subbundle $\mcN$ of $\mcV$ such that the composite $\mcN \migiincl \mcV \xrightarrow{\delta} F^{(\N)*}_{X/k}(\mcN^\nabla)$ is an isomorphism.
By passing to this composite, we identify $\mcN$ with $F^{(\N)*}_{X/k}(\mcN^\nabla)$ and regard $\delta$ as a left inverse  of the natural inclusion $\mcN \migi \mcV$.
It is verified  that the quadruple ${^\A \mcV}^\diamondsuit := (\mcV, \nabla_{\mcV^\nabla}^{\mr{can}(\N-1)}, \mcN, \delta)$ forms a dormant affine-indigenous
$\mcD_X^{(\N-1)}$-modules
 with ${^\A \mcV}^{\diamondsuit \Rightarrow \heartsuit} \cong \mcS^\heartsuit$.
This implies the desired surjectivity of $\overline{\xi}_\N^{\diamondsuit \Rightarrow \heartsuit}$.

The injectivity of $\overline{\xi}_\N^{\diamondsuit \Rightarrow \heartsuit}$ follows from an argument similar to the proof of the injectivity of ${\overline{\zeta}}_\N^{\diamondsuit \Rightarrow \heartsuit}$ discussed above.
This completes the proof of the theorem.
\end{proof}
\SSP

\begin{rema}[$F^\N\text{-}(\mr{PGL}_{n+1}, \mbP^n)$-structures\,=\,$F^\N$-projective structures] \label{Remmm}  
In ~\cite[Example 2.3.3]{Wak12}, we introduced   the notion of an $F^\N\text{-}(\mr{PGL}_{n+1}, \mbP^n)$-structure (resp., an $F^\N$-$(\mr{Aff}_n, \mbA^n)$-structure) as  a specific type of Frobenius-Ehresmann structure.
At first glance, this notion differs from  that of an $F^\N$-projective  (resp., $F^\N$-affine) structure:
 The former is defined as an {\it \'{e}tale} sheaf of locally defined \'{e}tale
morphisms to $\mbP^n$, whereas the latter is defined as a {\it Zariski} sheaf.
However, Theorem \ref{P01111} implies that, under certain assumptions,  
these notions  are  actually equivalent.
Under this equivalence,  the map $\overline{\zeta}^{\diamondsuit \Rightarrow \heartsuit}_\N$ (resp., $\overline{\xi}^{\diamondsuit \Rightarrow \heartsuit}_\N$) coincides with the inverse of the bijection ``$\zeta_\N^{\heartsuit \Rightarrow \spadesuit}$" appearing   in ~\cite[Theorem B]{Wak12}.
 \end{rema}

\LSP
\subsection{Rigidification by theta characteristics} \label{SShjk06g7}

We  consider a rigidification of dormant indigenous  and affine-indigenous $\mcD_X^{(\N-1)}$-modules constructed  by fixing their determinants by means of a (generalized) theta characteristic. 

By  a  {\bf theta characteristic} of $X$, we shall mean a line bundle $\varTheta$ on $X$ together with an isomorphism $\varTheta^{\otimes (n+1)} \isom \omega_X$.
This notion is  well-known and has been  studied in the case where  the underlying variety $X$ is a curve (cf. ~\cite{Mum}).
(Note that there always exists a theta characteristic of any smooth curve.)
In the following,   theta characteristics will be generalized 
in terms of higher-level $\mcD_X$-modules.

Let us fix an element $\N$ of $\mbZ_{>0}\sqcup \{ \infty \}$.

\SSP
\bde \label{D00ggh50}
 An {\bf $F^\N$-theta characteristic}
   of $X$ is a pair
 \begin{align}
 \vartheta := (\varTheta, \nabla_\vartheta),
 \end{align}
 where $\varTheta$ denotes a line bundle on $X$ and $\nabla_\vartheta$
  denotes a $\mcD_X^{(\N-1)}$-module structure  on $\omega_X^\vee \otimes \varTheta^{\otimes (n+1)}$ 
   with vanishing $p$-$(\N-1)$-curvature.
    \ede
\SSP

\begin{rema}[Truncation] \label{Rfghjf3401}  
Let $\N'$ be an element of $\mbZ_{>0}\sqcup \{ \infty \}$ with $\N' \geq \N$, and suppose that we are given an $F^{\N'}$-theta characteristic $\vartheta := (\varTheta, \nabla_\vartheta)$ of $X$.
Then, the $\N$-th truncation  $(\varTheta, \nabla_\vartheta |^{(\N )})$ (cf. \eqref{Ew23467}) forms an $F^\N$-theta characteristic of $X$.
In particular, if $X$ admits an $F^\infty$-theta characteristic, then there exists an $F^\N$-theta characteristic of $X$ for any positive integer $\N$.
By abuse of notation, we also use the notation  $\vartheta$ to denote  $(\varTheta, \nabla_\vartheta  |^{(\N )})$. 
 \end{rema}
\SSP

\begin{rema}[Theta characteristic] \label{Rfg1f3401} 
Suppose that $X$ admits a theta characteristic $\varTheta$.
 Then, 
 the pair 
 $(\varTheta, \nabla^{\mr{triv}(\N-1)}_{\mcO_X})$ (cf. \eqref{Eh1063}) forms  an $F^\N$-theta characteristic  on $X$, where  $\omega_X^\vee \otimes \varTheta^{\otimes (n+1)}$   is identified with $\mcO_X$ by means of the fixed isomorphism $\varTheta^{\otimes (n+1)} \isom \omega_X$.
 In this way, each theta characteristic specifies an $F^\N$-theta characteristic. 

 More generally, if $\omega_X$ is $(n+1)$-divisible in  the Neron-Severi group of $X$, then one can verify that $X$ admits an $F^\infty$-theta characteristic because of  the equivalence of categories \eqref{E3001} and the $p$-divisibility of $\mr{Pic}^0 (X)$, i.e., the identity component of the Picard scheme.
 \end{rema}
\SSP

\begin{rema}[Explicit construction] \label{Rfkk3401} 
If $p \nmid (n+1)$, then
 $X$ admits an $F^\N$-theta characteristic for any  positive integer $\N$ even when $X$ has no theta characteristics. 
In fact,  the assumption $p \nmid (n+1)$ implies that
there exists a pair of integers $(a, b)$ with $a (n+1)-1 = bp^\N$.
Then, we obtain a sequence of natural isomorphisms
\begin{align}
\omega_X^\vee \otimes (\omega_X^{\otimes a})^{\otimes (n+1)}  \isom \omega_X^{a (n+1)-1}  \left(= \omega_X^{\otimes bp^\N} \right) \isom F_{X/k}^{(\N)*}(\omega_{X^{(\N)}}^{\otimes b}).
\end{align}
Let $\nabla_{a,b}$ denote the $\mcD_{X}^{(\N-1)}$-action on $\omega_X^\vee \otimes (\omega_X^{\otimes a})^{\otimes (n+1)}$ corresponding, via this composite isomorphism, to  $\nabla^{\mr{can} (\N-1)}_{\omega_{X^{(\N)}}^{\otimes b}}$ (cf. \eqref{Ew9090}).
 The resulting pair
 \begin{align}
 (\omega_X^{\otimes a}, \nabla_{a,b})
 \end{align}
 forms an $F^\N$-theta characteristic of $X$, as desired.
  \end{rema}
\SSP

Let us fix an $F^\N$-theta  characteristic $\vartheta := (\varTheta, \nabla_\vartheta)$  
of $X$.

\SSP
\bde \label{Defgie}
    A dormant indigenous
    $\mcD_X^{(\N-1)}$-module $\mcV^\diamondsuit := (\mcV, \nabla, \mcN)$   is said to be {\bf $\vartheta$-rigidified}
 if it satisfies   
the two conditions (a) and (b) 
 described as follows:
\begin{itemize}
\item[(a)]
$\mcV = \mcD_{X, 1}^{(\N-1)}\otimes \varTheta$ and $\mcN = \varTheta  \left(= \mcD_{X, 0}^{(\N-1)}\otimes \varTheta \subseteq \mcD_{X, 1}^{(\N-1)}\otimes \varTheta  \right)$;
\item[(b)]
The determinant of $\nabla$ corresponds to $\nabla_\vartheta$ via the  isomorphism $\mr{det}(\mcV) \isom  \omega_X^\vee \otimes \mcN^{\otimes (n+1)}$ $\left(= \omega_X^\vee \otimes \varTheta^{\otimes (n+1)} \right)$
displayed  in  \eqref{Eh401}.
\end{itemize}
\ede

We shall write  
\begin{align} \label{Eh709}
{^\dagger}\DID{\N -1}{, \vartheta}
\end{align}
for the subset of 
$\DID{\N -1}{}$\!
  consisting of    isomorphism classes of $\vartheta$-rigidified dormant indigenous
    $\mcD_X^{(\N-1)}$-modules.
If we take  an $F^\N$-theta 
 characteristic of the form $(\varTheta, \nabla_{\mcO_X}^{\mr{triv}(\N-1)})$ for some theta characteristic $\varTheta$ (cf. Remark \ref{Rfg1f3401}),
   then, for simplicity, we write 
\begin{align}
{^\dagger}\DID{\N -1}{, \varTheta}
 := 
 {^\dagger}\DID{\N -1}{, (\varTheta, \nabla_{\mcO_X}^{\mr{triv}(\N-1)})}.
\end{align}

As for the ``affine" version, we make the following definition.

\SSP
\bde \label{Defei82}
A dormant affine-indigenous
    $\mcD_X^{(\N-1)}$-modules  ${^\A\mcV}^\diamondsuit := (\mcV, \nabla, \mcN, \delta)$ 
    is said to be {\bf rigidified}
    if it satisfies the two conditions (c) and (d) 
 described as follows:
 \begin{itemize}
\item[(c)]
$\mcV = \mcD_{X, 1}^{(\N-1)}\otimes \mcO_X$ and $\mcN = \mcO_X  \left(= \mcD_{X, 0}^{(\N-1)}\otimes \mcO_X \subseteq \mcD_{X, 1}^{(\N-1)}\otimes \mcO_X  \right)$;
\item[(d)]
The $\mcD_{X}^{(\N-1)}$-action on $\mcO_X$ induced from $\nabla$ via the surjection $\delta : \mcV \migisurj \mcN \left(= \mcO_X \right)$ coincides with $\nabla_{\mcO_X}^{\mr{triv}(\N-1)}$.
\end{itemize}
\ede
\SSP

We shall write
\begin{align} \label{Eh5001}
{^\dagger}\DAID{\N -1}{}
\end{align}
for the subset of
$\DAID{\N -1}{}$
  consisting of    isomorphism classes of
 rigidified   dormant affine-indigenous
    $\mcD_X^{(\N-1)}$-modules.

By restricting  
 $\zeta_\N^{\diamondsuit \Rightarrow \heartsuit}$ 
 and
  $\xi_\N^{\diamondsuit \Rightarrow \heartsuit}$, we obtain  maps
\begin{align} \label{Ep3400}
 {^\dagger}\zeta_{\N, \vartheta}^{\diamondsuit \Rightarrow \heartsuit} : {^\dagger}\DID{\N -1}{, \vartheta}\! 
\migi 
\Proj{\N}{}
\hspace{5mm} \text{and} \hspace{5mm}
{^\dagger}\xi_{\N}^{\diamondsuit \Rightarrow \heartsuit} : 
{^\dagger}\DAID{\N -1}{}\! \migi \Aff{\N}{},
\end{align}
respectively.

\SSP
\bpr \label{P06p0hgh} 
Assume that  $\Gamma (X, \mcO_X)=k$. (For example, this assumption is satisfied if  $X$ is proper over $k$.)
Then, the following assertions hold:
\begin{itemize}
\item[(i)]
The set 
${^\dagger}\DID{\N -1}{, \vartheta}$
 may be identified with the set (not only  the set of ``isomorphism classes") of  dormant indigenous 
 $\mcD_X^{(\N-1)}$-modules 
satisfying the conditions (a) and (b)
  described  before.
That is to say, 
there exists exactly one dormant indigenous 
$\mcD_X^{(\N-1)}$-module  representing each element of  
${^\dagger}\DID{\N -1}{, \vartheta}$.
\item[(ii)]
Suppose further that we are given an $F^\N$-theta characteristic  $\vartheta$ of $X$.
Then, the  set 
 ${^\dagger}\DAID{\N -1}{}$
 may be identified with the set (not only  the set of ``isomorphism classes") of  
  dormant affine-indigenous $\mcD_X^{(\N-1)}$-modules 
satisfying the conditions
  (c) and (d) described  before.
That is to say, 
there exists exactly one
  dormant affine-indigenous $\mcD_X^{(\N-1)}$-module  representing each element of  
 ${^\dagger}\DAID{\N -1}{}$.
\item[(iii)]
When $\N = \infty$, 
the assignment $\mcV^\diamondsuit \mapsto \{ \mcV^\diamondsuit |^{\langle \N \rangle} \}_{\N \in \mbZ_{>0}}$ (resp., ${^\mbA\mcV}^\diamondsuit \mapsto \{ {^\mbA\mcV}^\diamondsuit |^{\langle \N \rangle} \}_{\N \in \mbZ_{>0}}$) defines a bijection of sets
\begin{align} \label{Eh467}
{^\dagger}\DID{\infty}{, \vartheta}
\isom  \varprojlim_{\N \in \mbZ_{>0}}
{^\dagger}\DID{\N -1}{, \vartheta}
\ \left(\text{resp.,} \
{^\dagger}\DAID{\infty}{}
\isom  \varprojlim_{\N \in \mbZ_{>0}}
{^\dagger}\DAID{\N -1}{}
  \right).
\end{align}
\end{itemize}
  \epr
\begin{proof}
Assertions (i) and  (ii) follow directly from Proposition \ref{P3fg6h9}, (ii).
Assertion (iii) follows from assertions (i) and  (ii) together with the identification $\varinjlim_{\N \in \mbZ_{>0}} \mcD_X^{(\N-1)} \cong \mcD_X^{(\infty)}$ (cf. \eqref{E1004}).
\end{proof}

\LSP
\subsection{Correspondence  via  rigidification} \label{SShjk06g7}

In what follows, we prove  that the set of rigidified dormant indigenous (resp.,  affine-indigenous) $\mcD_X^{(\N-1)}$-modules  correspond bijectively to the set of
$F^\N$-projective (resp., $F^\N$-affine) structures including the case of $\N = \infty$ (cf. Theorem \ref{Cfrt0gh}).

\SSP
\bpr \label{P06p0gh} 
\begin{itemize}
\item[(i)]
Suppose that $p \nmid (n+1)$, $\N \in \mbZ_{>0}$,  and $X$ is quasi-projective over $k$.
Then, 
for 
each $F^\N$-theta characteristic $\vartheta := (\varTheta, \nabla_\vartheta)$ of $X$, 
 the composite
 \begin{align} \label{Eh466}
 {\zeta}_{\N, \vartheta}^{\dagger  \Rightarrow \overline{(-)}}
 : 
 {^\dagger}\DID{\N -1}{, \vartheta}
\xrightarrow{\mr{inclusion}}
\DID{\N -1}{}
 \xrightarrow{\mr{quotient}}
 \eDID{\N -1}{}
  \end{align}
 is injective.
 If, moreover, $X$ is projective over $k$, or more generally, satisfies the equality $\Gamma (X, \mcO_X) =k$,
 then this composite is bijective.
 \item[(ii)]
For any $\N \in \mbZ_{>0}$, the composite
\begin{align} \label{Eh3041}
 \xi_{\N}^{\dagger  \Rightarrow \overline{(-)}}
 : 
 {^\dagger}\DAID{\N -1}{}
\xrightarrow{\mr{inclusion}}
\DAID{\N -1}{}
 \xrightarrow{\mr{quotient}}
 \eDAID{\N -1}{}
\end{align}
is bijective.
 \end{itemize}
  \epr
\begin{proof}
First, we shall prove the former assertion of   (i), i.e., the injectivity  of $ {\zeta}_{\N, \vartheta}^{\dagger  \Rightarrow \overline{(-)}}$.
Let $\mcV_i^\diamondsuit := 
(\mcV_i, \nabla_{i}, \mcN_i)$
($i=1,2$) be dormant indigenous $\mcD_X^{(\N-1)}$-modules classified by
${^\dagger}\DID{\N -1}{, \vartheta}$
   (hence $\mcV_i = \mcD_{X, 1}^{(\N-1)}\otimes\varTheta$ and $\mcN_i = \varTheta$) such that $[\mcV_1^{\diamondsuit}] = [\mcV_2^{\diamondsuit}]$.
By 
the definition of {\it being $\mbG_m$-equivalent},
there exists  a collection  $(\mcL, \nabla_\mcL, \eta)$ consisting of  an invertible  $\mcD_X^{(\N-1)}$-module $(\mcL, \nabla_\mcL)$ with    vanishing $p$-$(\N-1)$-curvature  and 
 an isomorphism $\eta : \mcV_1^\diamondsuit \isom (\mcV^\diamondsuit_2)_{\otimes (\mcL, \nabla_\mcL)}$.
This isomorphism $\eta$ restricts to an isomorphism $\mcN_1 \left(=\varTheta\right) \isom \mcN_2 \otimes \mcL \left(= \varTheta \otimes \mcL \right)$, which implies that $\mcL \cong \mcO_X$.
According to the equivalence of categories \eqref{E3002},
$\mcS ol (\nabla_\mcL)$ forms a line bundle on $X^{(\N)}$ with $F_{X/k}^{(\N)*}(\mcS ol (\nabla_\mcL)) \cong \mcO_X$.
It follows from Lemma \ref{L06geh7} below that
\begin{align} \label{Eh600}
\mcS ol (\nabla_\mcL)^{\otimes p^{n\N}} \cong \mcO_{X^{(\N)}}.
\end{align}
On the other hand, 
the determinant of $\eta$ gives an isomorphism
 \begin{align}
 \omega_X^\vee \otimes \varTheta^{\otimes (n+1)}  \left(\stackrel{\eqref{Eh401}}{\cong} \mr{det}(\mcV_1)  \right) \isom  \omega_X^\vee \otimes \varTheta^{\otimes (n+1)} \otimes \mcL^{\otimes (n+1)}\left(\stackrel{\eqref{Eh401}}{\cong} \mr{det}(\mcV_2 \otimes \mcL) \right)
 \end{align}
  compatible with the respective $\mcD_X^{(\N-1)}$-actions $\nabla_\vartheta$
    and $\nabla_\vartheta \otimes \nabla_{\mcL}^{\otimes (n+1)}$, where $\nabla_{\mcL}^{\otimes (n+1)}$ denotes the $(n+1)$-st tensor product of $\nabla_\mcL$.
     This implies that $(\mcL^{\otimes (n+1)}, \nabla_\mcL^{\otimes (n+1)})\cong (\mcO_X, \nabla_{\mcO_X}^{\mr{triv}(\N-1)})$.
  Since the equivalence of categories \eqref{E3002} is compatible with taking tensor products, we see  that
  \begin{align} \label{Eh601}
  \mcS ol (\nabla_\mcL)^{\otimes (n+1)} \cong \mcS ol (\nabla_{\mcL^{\otimes (n+1)}}) \cong \mcS ol (\nabla_{\mcO_X}^{\mr{triv}(\N-1)}) \cong \mcO_{X^{(\N)}}.
  \end{align}
 It follows from \eqref{Eh600},  \eqref{Eh601}, and the assumption $(p^{n\N}, n+1) =1$  that  $\mcS ol (\nabla_\mcL)\cong\mcO_{X^{(\N)}}$, so  $(\mcL, \nabla_\mcL)$ is isomorphic to the trivial $\mcD_X^{(\N-1)}$-module $(\mcO_X, \nabla^{\mr{triv}(\N-1)}_{\mcO_X})$.
Consequently, $\mcV_1^\diamondsuit$ turns out to be isomorphic to $\mcV_2^\diamondsuit$.
This completes the proof of assertion (i).

Next, we shall consider the latter assertion of (i), i.e., the bijectivity of ${\zeta}_{\N, \vartheta}^{\dagger  \Rightarrow \overline{(-)}}$  under the assumption $\Gamma (X, \mcO_X) =k$.
Let $\mcV^\diamondsuit := (\mcV, \nabla, \mcN)$ be  a dormant indigenous $\mcD_X^{(\N-1)}$-module.
Consider  the following sequence of natural isomorphisms:
\begin{align}
\omega_X^\vee \otimes \varTheta^{\otimes (n+1)} \otimes \mr{det} (\mcV)^\vee   
\isom \omega_X^\vee \otimes \varTheta^{\otimes (n+1)} \otimes (\omega^\vee_X \otimes \mcN^{\otimes (n+1)})^\vee 
 \isom \mcL^{\otimes (n+1)},
\end{align}
where $\mcL := \varTheta \otimes \mcN^\vee$ and  the first isomorphism follows from  \eqref{Eh401}.
By passing to the composite of them, we obtain a $\mcD_X^{(\N-1)}$-action $\nabla_{\mcL^{\otimes (n+1)}}$ on 
$\mcL^{\otimes (n+1)}$ corresponding to the tensor product of  $\nabla_\vartheta$ and  the dual of the determinant of $\nabla$.
It follows from Lemma \ref{L06h7} below that
there exists a $\mcD_{X}^{(\N-1)}$-action $\nabla_\mcL$ with vanishing $p$-$(\N-1)$-curvature whose $(n+1)$-st tensor product coincides with $\nabla_{\mcL^{\otimes (n+1)}}$.
Thus, we obtain a dormant indigenous $\mcD_X^{(\N-1)}$-module   $\mcV^\diamondsuit_{\otimes (\mcL, \nabla_\mcL)}$ (cf. \eqref{Ep2}).
If we write  
$(\mcV', \nabla', \mcN') := \mcV^\diamondsuit_{\otimes (\mcL, \nabla_\mcL)}$, then the line bundle 
 $\mcN'  \left( = \mcN \otimes (\varTheta \otimes \mcN^\vee)\right)$  can be identified with $\varTheta$ via the natural isomorphism.
We also identify $\mcV'$ with $\mcD_{X}^{(\N-1)}\otimes \varTheta$ via the isomorphism  $\mr{KS}^\circledast_{\mcV^\diamondsuit_{\otimes (\mcL, \nabla_\mcL)}}$ (cf. \eqref{Eh101}).
Under  these identifications, 
 $\mcV^\diamondsuit_{\otimes (\mcL, \nabla_\mcL)}$ satisfies the condition (a) in Definition \ref{Defgie}.
Moreover, it follows from  the construction of $\nabla_\mcL$ that 
$\mcV^\diamondsuit_{\otimes (\mcL, \nabla_\mcL)}$ also satisfies the condition (b).
Since 
$[\mcV^\diamondsuit_{\otimes (\mcL, \nabla_\mcL)}] = [\mcV^\diamondsuit]$,
the desired surjectivity has been  proved.  
Thus, we finish  the proof of assertion (i).

Assertion (ii) follows from an argument similar to the argument in the proof of Theorem \ref{P01111}, (ii). 
\end{proof}
\SSP

The following two lemmas were used in the above proposition.

\SSP
\ble[cf. ~\cite{EGA2}, Section 6.5] \label{L06geh7} 
Let $f: X \migi X'$ be a finite and faithfully flat morphism between varieties over $k$ and $\mcM$ a line bundle on $X'$.
(In the context of this subsection, we consider the case where $X' = X^{(\N)}$ and $f = F_{X/k}^{(\N)}$ for a positive integer $\N$.)
If $f^*(\mcM) \cong \mcO_X$, then we have $\mcM^{\otimes d} \cong \mcO_{X'}$, where $d$ denotes the degree of $f$.
 \ele
\begin{proof}
Since 
$f$
 is finite and faithfully flat of degree $d$,
 the direct image 
 $f_*(\mcO_X)$
 of $\mcO_X$ forms a rank 
 $d$
 vector bundle on 
 $X'$.
Hence, we have the following sequence of isomorphisms:
\begin{align}
\mr{det} (f_*(\mcO_X)) &\isom \mr{det}(f_*(f^*(\mcM))) 
\isom  \mr{det} ( \mcM \otimes f_*(\mcO_X)) 
\isom  \mcM^{\otimes d} \otimes \mr{det} (f_*(\mcO_X)), 
\end{align}
where the first arrow follows from the assumption, i.e., $f^*(\mcM) \cong \mcO_X$,  and the second arrow follows from the projection formula.
This implies that $\mcM^{\otimes d} \cong \mcO_X$, as desired.
\end{proof}
\SSP

\ble \label{L06h7} 
Suppose that $X$ is projective over $k$, or more generally, satisfies the equality $\Gamma (X, \mcO_X) = k$.
 Let $\N$
 be a positive integer and $l$  a positive integer prime to $p$.
 Also, let $\mcL$ be a line bundle on $X$ and $\nabla_{\mcL^{\otimes l}}$ a  $\mcD_X^{(\N-1)}$-action on $\mcL^{\otimes l}$ with vanishing $p$-$(\N-1)$-curvature.
 Then, there exists   a unique  $\mcD_X^{(\N-1)}$-action $\nabla_\mcL$ on $\mcL$ with vanishing $p$-$(\N-1)$-curvature whose $l$-th tensor product $\nabla_{\mcL}^{\otimes l}$ coincides with $\nabla_{\mcL^{\otimes l}}$.
  \ele
\begin{proof}
Let us consider the Kummer exact  sequence
\begin{align} \label{Eh102}
1 \longmigi \mu_l \xrightarrow{\mr{inclusion}} \mbG_m \xrightarrow{(-)^l} \mbG_m \longmigi 1,
\end{align}
where $\mu_l$ denotes the group of $l$-th roots of unity.
Since the equality $\Gamma (X, \mbG_m) = k^\times$ holds by assumption,
the above sequence remains exact after applying the functor $\Gamma (X, -)$.
Thus, the Kummer sequence induces the following exact sequence:
\begin{align} \label{Eh1809}
0\longmigi  \check{H}^1_{\text{\'{e}t}}(X, \mu_{l}) \longmigi 
\check{H}^1_{\text{\'{e}t}}(X, \mbG_m) \longmigi 
\check{H}^1_{\text{\'{e}t}}(X, \mbG_m) \longmigi 
\check{H}^2_{\text{\'{e}t}}(X, \mu_{l}).
\end{align}
The $\N$-th relative   
Frobenius morphisms of $\mu_l$,  $\mbG_m$ are compatible with the morphisms in
\eqref{Eh102} and hence induce 
the following endomorphism of the sequence \eqref{Eh1809}:
\begin{align} \label{Eh608}
\begin{CD}
0 @>>> \check{H}^1_{\text{\'{e}t}}(X, \mu_{l}) @>>> \check{H}^1_{\text{\'{e}t}}(X, \mbG_m) @>\lambda^{\migi}>> \check{H}^1_{\text{\'{e}t}}(X, \mbG_m)  @>>> \check{H}^2_{\text{\'{e}t}}(X, \mu_{l})
\\
@. @V \wr VV @V {^\downarrow \lambda} VV @VV \lambda^{\downarrow} V @VV \wr V
\\
0 @>>>\check{H}^1_{\text{\'{e}t}}(X, \mu_{l}) @>>> \check{H}^1_{\text{\'{e}t}}(X, \mbG_m) @>> \lambda_{\migi} >  \check{H}^1_{\text{\'{e}t}}(X, \mbG_m) @>>> \check{H}^2_{\text{\'{e}t}}(X, \mu_{l}), 
\end{CD}
\end{align}
where both the leftmost and  rightmost vertical arrows are isomorphisms because of the assumption $(p, l) =1$.
Denote by $[\mcL]$, $[\mcL^{\otimes l}]$, and $[\mcS ol (\nabla_{\mcL^{\otimes l}})]$ the elements of $\check{H}^1_{\text{\'{e}t}}(X, \mbG_m)$ represented by $\mcL$, $\mcL^{\otimes l}$, and $\mcS ol (\nabla_{\mcL^{\otimes l}})$ respectively.
By definition, the equalities  $\lambda_{\migi}([\mcL]) = [\mcL^{\otimes l}]$ and $\lambda^{\downarrow} (\mcS ol (\nabla_{\mcL^{\otimes l}})) = [\mcL^{\otimes l}]$ hold.
One verifies from
 a routine diagram-chasing argument that there exists 
  a line bundle $\mcL_0$, uniquely determined  up to isomorphism,  such that the element $[\mcL_0]$ of $\check{H}^1_{\text{\'{e}t}}(X, \mbG_m)$
  represented by $\mcL_0$ satisfies the equalities 
  $\lambda^{\migi}([\mcL_0]) = [\mcS ol (\nabla_{\mcL^{\otimes l}})]$ 
 and ${^\downarrow \lambda} ([\mcL_0]) = [\mcL]$.
Because of the equality $\Gamma (X, \mbG_m) = k^\times$,
one can find  isomorphisms $\eta_1 : \mcL_0^{\otimes l} \isom \mcS ol (\nabla_{\mcL^{\otimes l}})$ and $\eta_2 : F_{X/k}^{(\N)*}(\mcL_0) \isom \mcL$  making the following square diagram commutative:
\begin{align} \label{Eh700}
\vcenter{\xymatrix@C=56pt@R=36pt{
F^{(\N)*}_{X/k} (\mcL_0^{\otimes l})\ar[r]^-{F^{(\N)*}_{X/k}(\eta_1)}_-{\sim} \ar[d]_-{\wr}& F^{(\N)*}_{X/k} (\mcS ol (\nabla_{\mcL^{\otimes l}}))  \ar[d]^-{\wr}
\\
F^{(\N)*}_{X/k} (\mcL_0)^{\otimes l} \ar[r]^-{\sim}_-{\eta_2^{\otimes l}}& \mcL^{\otimes l},
}}
\end{align}
where the vertical arrows denote the natural isomorphisms.
We shall denote by  $\nabla_{\mcL}$ the $\mcD_X^{(\N-1)}$-action on $\mcL$ corresponding to the $\mcD_X^{(\N-1)}$-action  $\nabla^{\mr{can} (\N-1)}_{\mcL_0}$ on $F_{X/k}^{(\N)*}(\mcL_0)$ (cf. \eqref{Ew9090}) via $\eta_2$.
It follows from  the commutativity of \eqref{Eh700} that $\nabla_\mcL$ satisfies the required conditions.
Moreover, 
 such a $\mcD_{X}^{(\N-1)}$-action  $\nabla_\mcL$ does not depend on the choice of $(\eta_1, \eta_2)$, i.e.,  is uniquely determined; 
 this is because  any automorphism of $\mcL$
  may be given by multiplication by some element of $k^\times  \left(= \Gamma (X, \mcO_X^\times) \right)$, which is compatible with $\nabla_\mcL$.
 This completes the proof of the lemma.  
\end{proof}
\SSP

\bt \label{Cfrt0gh} 
 Let $\N \in \mbZ_{>0}\sqcup \{ \infty \}$.
 Then, the following assertions hold:
\begin{itemize}
\item[(i)]
Suppose that $p \nmid (n+1)$ and $X$ is projective over $k$.
Then, 
 for each $F^\N$-theta characteristic $\vartheta := (\varTheta, \nabla_\vartheta)$ of $X$, 
 the map
   ${^\dagger}\zeta_{\N, \vartheta}^{\diamondsuit \Rightarrow \heartsuit} : 
  {^\dagger}\DID{\N -1}{, \vartheta}
  \migi 
  \Proj{\N}{}$
   is  bijective.
 \item[(ii)]
 Suppose that $X$ is proper over $k$.
Then, the  map
  ${^\dagger}\xi_{\N}^{\diamondsuit \Rightarrow \heartsuit} : 
 {^\dagger}\DAID{\N -1}{}
 \migi 
 \Aff{\N}{}$
  is  bijective.
 \end{itemize}
  \et
\begin{proof}
For $\N \in \mbZ_{>0}$,
both assertions  follow from Theorem \ref{P01111} and Proposition \ref{P06p0gh}.
Notice that the case of $\N = \infty$
 follows from the cases of $\N < \infty$ together with the identifications 
     $\varprojlim_{\N \in \mbZ_{>0}}
     {^\dagger}\DID{\N -1}{, \vartheta}
        =   {^\dagger}\DID{\infty}{, \vartheta}$   and
       $\varprojlim_{\N \in \mbZ_{>0}}
     {^\dagger}\DAID{\N -1}{}
      =   {^\dagger}\DAID{\infty}{}$
      given by 
         \eqref{Eh467}.
In fact, the various  maps
under consideration
  are, by construction, compatible with truncation to lower levels.
\end{proof}
\SSP

By applying  the above theorem, we have the following comparison  between 
$\eDID{\infty}{}$
 and $\varprojlim_{\N \in \mbZ_{>0}} 
 \eDID{\N -1}{}$, as well as its affine version.

\SSP
\bpr \label{C06p0gh}
 The following assertions hold:
\begin{itemize}
\item[(i)]
 Suppose that $p \nmid (n+1)$,  $X$ is projective over $k$, and $X$ admits an $F^\infty$-theta characteristic. Then,
 the natural map 
 \begin{align}
 \eDID{\infty}{}
  \migi \varprojlim_{\N \in \mbZ_{>0}}  
  \eDID{\N -1}{}
 \end{align}
  is bijective.
 In particular, the map
  $\overline{\zeta}_\infty^{\diamondsuit \Rightarrow \heartsuit} : 
 \eDID{\infty}{}
  \migi 
  \Proj{\infty}{}$
 is bijective.
 \item[(ii)]
 Suppose that $X$ is proper over $k$.
 Then, the natural map
 \begin{align}
 \eDAID{\infty}{}
 \migi \varprojlim_{\N \in \mbZ_{>0}}
 \eDAID{\N -1}{}
 \end{align}
 is bijective.
 In particular, the map 
    $\overline{\xi}_\infty^{\diamondsuit \Rightarrow \heartsuit} : 
 \eDAID{\infty}{}
  \migi 
  \Aff{\infty}{}$
 is  bijective.
 \end{itemize}
  \epr
\begin{proof}
We shall prove the former assertion of (i).
Let us fix an $F^\infty$-theta characteristic $\vartheta := (\varTheta, \nabla_\vartheta)$ of $X$ (cf. Remark \ref{Rfkk3401}).
Denote by $\zeta^{\diamondsuit}_{\infty \Rightarrow \mr{lim}}$ the map 
$\eDID{\infty}{} \migi \varprojlim_{\N \in \mbZ_{>0}} \eDID{\N -1}{}$
 in question.
The surjectivity  of 
 $\zeta^{\diamondsuit}_{\infty \Rightarrow \mr{lim}}$ 
 follows immediately from 
the bijectivities of ${^\dagger}\zeta_{\N, \vartheta}^{\diamondsuit \Rightarrow \heartsuit}$'s for various $\N \in \mbZ_{>0}\sqcup \{ \infty \}$ (cf. Theorem \ref{Cfrt0gh}) and the fact that
the equality ${^\dagger}\zeta_{\infty, \vartheta}^{\diamondsuit \Rightarrow \heartsuit} = (\varprojlim_{\N \in \mbZ_{>0}}{^\dagger}\zeta_{\N, \vartheta}^{\diamondsuit \Rightarrow \heartsuit}) \circ \zeta^{\diamondsuit}_{\infty \Rightarrow \mr{lim}}$ holds under 
 the identification $\varprojlim_{\N \in \mbZ_{>0}}
 \mr{Proj}_\N
   = 
   \mr{Proj}_\infty$ (cf. \eqref{Ew010}).

Next, we shall consider the injectivity of  $\zeta^{\diamondsuit}_{\infty \Rightarrow \mr{lim}}$.
Let $\mcV_i^\diamondsuit := (\mcV_i, \nabla_{i}, \mcN_i)$ ($i=1,2$) be dormant indigenous 
$\mcD_X^{(\infty)}$-modules
  such that  $[\mcV_1^{\diamondsuit}|^{\langle \N \rangle}] =  [\mcV_2^{\diamondsuit}|^{\langle \N \rangle}]$ in 
$\eDID{\N -1}{}$
 for every $\N \in \mbZ_{>0}$.
Then, for each $\N \in \mbZ_{>0}$, there exist an invertible $\mcD_X^{(\N-1)}$-module $(\mcL_\N, \nabla_{\mcL, \N})$ with vanishing $p$-$(\N-1)$-curvature
  and  an isomorphism  
\begin{align}
\eta_\N : \mcV_1^\diamondsuit |^{\langle \N \rangle}\isom (\mcV_2^\diamondsuit |^{\langle \N \rangle})_{\otimes (\mcL_\N, \nabla_{\mcL, \N})}.
\end{align}
The isomorphism $\eta_\N$ restricts to an isomorphism $(\mcV_1 \supseteq) \  \mcN_1 \isom \mcN_2 \otimes \mcL_\N \left( \subseteq \mcV_2 \otimes \mcL_\N\right)$, which gives an isomorphism
\begin{align} \label{Eh1100}
\mcL_\N \isom   \mcN_1 \otimes \mcN_2^\vee.
\end{align}
Also, $\eta_\N$ induces, via taking determinants,  an isomorphism  
\begin{align}
\eta_\N^\mr{det} : \mr{det}(\mcV_1) \isom \mr{det} (\mcV_2 \otimes \mcL_\N) \left(\cong \mr{det} (\mcV_2) \otimes \mcL_\N^{\otimes (n+1)} \right).
\end{align}
For each $i=1,2$, 
let  $\nabla^{\mr{det}}_{i}$  denote 
the $\mcD_X^{(\N-1)}$-action on $\mr{det}(\mcV_i)$ induced naturally from $\nabla_{i}$.
Then,   $\eta_\N^\mr{det}$ gives an isomorphism  of $\mcD_X^{(\N-1)}$-modules
\begin{align} \label{Eh1101}
(\mcL_\N^{\otimes (n+1)}, \nabla_{\mcL, \N}^{\otimes (n+1)}) \isom  (\mr{det}(\mcV_1) \otimes  \mr{det}(\mcV_2)^\vee, \nabla_{1}^\mr{det} \otimes (\nabla_{2}^\mr{det})^\vee),
\end{align}
where $\nabla_{\mcL, \N}^{\otimes (n+1)}$ denotes the $(n+1)$-st tensor product of $\nabla_{\mcL, \N}$.
This isomorphism  
 is compatible with  the $(n+1)$-st tensor product of \eqref{Eh1100}  under the composite of natural isomorphisms
\begin{align}
(\mcN_1 \otimes \mcN_2^\vee)^{\otimes (n+1)} &\isom \mcN_1^{\otimes (n+1)} \otimes (\mcN_2^{\otimes (n+1)})^\vee \\
&\isom (\omega_X \otimes \mr{det}(\mcV_1)) \otimes (\omega_X \otimes \mr{det}(\mcV_1))^\vee \notag \\
&\isom \mr{det}(\mcV_1) \otimes \mr{det} (\mcV_2)^\vee, \notag
\end{align}
where the second arrow arises from \eqref{Eh401}.
Hence,
 the uniqueness portion of Lemma \ref{L06h7} implies 
 that, for another positive integer $\N'$ with $\N' > \N$,
the $\N$-th truncation 
of $(\mcL_{\N'}, \nabla_{\mcL, \N'})$ obtained in the same manner as above
is isomorphic to $(\mcL_{\N}, \nabla_{\mcL, \N})$.
It follows that the collection $\left\{(\mcL_{\N}, \nabla_{\mcL, \N}) \right\}_{\N \in \mbZ_{>0}}$ yields
a $\mcD_X^{(\infty)}$-module $(\mcL_\infty, \nabla_{\mcL, \infty})$ such  that, for each $\N \in \mbZ_{>0}$,
the $\N$-th truncations of 
$\mcV_1^\diamondsuit$ and $(\mcV_2^\diamondsuit)_{\otimes (\mcL_\infty, \nabla_{\mcL, \infty})}$ are isomorphic via $\eta_\N$.
By   Proposition \ref{P3fg6h9}, (ii),  
we can obtain an isomorphism  $\mcV_1^\diamondsuit \isom (\mcV_2^\diamondsuit)_{\otimes (\mcL_\infty, \nabla_{\mcL, \infty})}$ determined by  
the collection $\{ \eta_\N \}_{\N \in \mbZ_{>0}}$
 after possibly composing  each $\eta_\N$ with multiplication by some element of $k^\times$.
Thus, $\mcV_1^\diamondsuit$ and $\mcV_2^\diamondsuit$ are $\mbG_m$-equivalent.
This completes the proof of the desired injectivity.

The latter assertion follows directly  from the former assertion and Theorem \ref{P01111}, (i).

Finally,  assertion (ii) follows from \eqref{Ew010},  Proposition \ref{P06p0hgh}, (iii),  Proposition \ref{P06p0gh}, (ii),  and Theorem \ref{P01111}, (ii).
This completes the proof of the proposition.
\end{proof}

\LSP
\subsection{Chern class formula} \label{SS041}

In this subsection, we prove  
 necessary conditions on Chern classes for the existence of  $F^\N$-projective or  $F^\N$-affine structures (cf. Theorem \ref{T016} below). 

For each vector bundle $\mcV$ on $X$ and each $l = 0, 1, \cdots, \mr{rank}(\mcV)$,
we write $c_l (\mcV)$  
for the $l$-th Chern class of $\mcV$, which is defined as an element of   $CH^l (X)$, i.e., the Chow group of codimension-$l$ cycles modulo rational equivalence.
Also, write $c_l^{\mr{crys}} (\mcV)$ for the $l$-th crystalline Chern class of $\mcV$ (cf.  ~\cite{BI} or ~\cite[Chapter III, Theorem 1.1.1]{Gro}), which is defined as 
 an element of  $H^{2l}_\mr{crys}(X/W)$.
 Here,  $H^{2l}_\mr{crys}(X/W)$ denotes  the $2l$-th crystalline cohomology group of $X$ over $W$, where $W$ denotes the ring of Witt vectors over $k$.
As usual, 
we set $c_l (X)  := c_l (\mcT_X)$ and $c_l^{\mr{crys}}(X) := c_l^{\mr{crys}}(\mcT_X)$.

\SSP
\bt \label{T016} 
Let $\N \in \mbZ_{>0} \sqcup \{ \infty \}$ and let $X$ be a smooth quasi-projective  variety over $k$
admitting
an $F^\N$-projective   structure.
Suppose further  that $p \nmid (n+1)$.
Then, the following assertions hold:
\begin{itemize}
\item[(i)]
In the case where
$\N \neq \infty$, 
  the   equality
\begin{align} \label{E45}
c_l (X)= \frac{1}{(n+1)^l} \cdot \binom{n+1}{l} \cdot c_1 (X)^l 
\end{align}
holds  after reduction modulo $p^{\N}$ 
  for each positive integer $l$.
In particular, if $X$ admits an $F^\infty$-projective 
 structure, then the above equality holds after reduction modulo $p^\N$  
  for any positive integers $\N$ and $l$.
Moreover, the same assertions hold for $c_l^{\mr{crys}}(X)$'s.
\item[(ii)]
In the case where $\N = \infty$ and $X$ is projective over $k$,
 the equality \eqref{E45}   in $H^{2l}_{\mr{crys}}(X/W)$, where ``$c_1(X)$'' and ``$c_l(X)$''  are replaced by   $c_1^{\mr{crys}}(X)$ and $c_l^{\mr{crys}}(X)$,  respectively,   holds without any reduction for every  positive integer $l$.
\end{itemize}
  \et
\begin{proof}
We shall prove the first  assertion  of (i) by induction on $l$.
The case of $l=1$ is clear.
Next, suppose that the assertion with $l$ replaced by any $l' <l$
 has been proved.
By assumption, $X$ admits an $F^\N$-projective structure $\mcS^\heartsuit$. 
It follows from 
Theorem \ref{P01111}, (i),
that
there exists a dormant indigenous $\mcD_X^{(\N-1)}$-module  $\mcV^\diamondsuit := (\mcV, \nabla, \mcN)$  with
$\mcV^{\diamondsuit \Rightarrow \heartsuit} = \mcS^\heartsuit$.
Since $(\mcV, \nabla)$ has vanishing $p$-$(\N-1)$-curvature,
we have $F_{X/k}^{(\N)*}(\mcS ol (\nabla)) \cong \mcV$ (cf. \eqref{E3002}).
If $\mcS ol (\nabla)^F$ denotes the vector bundle on $X$ corresponding to $\mcS ol (\nabla)$ via the isomorphism $W_{X/k}: X^{(\N)} \isom X$,
then $F_X^{\N*}(\mcS ol (\nabla)^F) \cong \mcV$, and hence, the equality  $c_l (\mcV) = p^{l\N} \cdot c_l (\mcS ol (\nabla)^F)$ holds (cf. ~\cite[the proof of Lemma 2.1]{G}).
In particular, the equality 
\begin{align} \label{Er1}
c_l (\mcV) = 0
\end{align}
holds after reduction modulo $p^{\N}$.
By Proposition \ref{P36h9}, $\mcV$ may be identified with $\mcD_{X, 1}^{(\N-1)} \otimes \mcN$ via $\mr{KS}_{(\mcV, \nabla, \mcN)}^\circledast$ and we have
\begin{align} \label{Er2}
c_l (\mcV) = c_l (\mcD_{X, 1}^{(\N-1)} \otimes \mcN).
\end{align}
The vector bundle $\mcD_{X, 1}^{(\N-1)} \otimes \mcN$
 fits into the following short exact sequence:
\begin{align}
0 \migi \mcN  \left(= \mcD_{X, 0}^{(\N-1)} \otimes \mcN\right) \migi \mcD_{X, 1}^{(\N-1)} \otimes \mcN \migi \mcT_X \otimes \mcN  \left(= (\mcD_{X, 1}^{(\N-1)} / \mcD_{X, 0}^{(\N-1)} )\otimes \mcN\right) \migi 0.
\end{align} 
It follows that
\begin{align} \label{Er5}
c_l (\mcD_{X, 1}^{(\N-1)} \otimes \mcN) &= \sum_{i+j=l} c_i (\mcN) c_j (\mcT_X \otimes \mcN) \\
&= c_l (\mcT_X \otimes \mcN) + c_1 (\mcN) \cdot c_{l-1} (\mcT_X \otimes \mcN). \notag
\end{align}
Also, consider  the equality
 \begin{align} \label{Er6}
 c_l (\mcT_X \otimes \mcN) = \sum_{i=0}^l \binom{n- i}{l-i} \cdot c_i (X) c_1 (\mcN)^{l-i}
 \end{align}
 resulting from  ~\cite[Chapter 3, Example 3.2.2]{Ful}.
 By combining \eqref{Er1}, \eqref{Er2}, \eqref{Er5}, and \eqref{Er6}, we obtain the following  sequence of equalities in the Chow group $CH^l (X)$ modulo $p^\N$:
\begin{align} \label{Er11}
0 &\stackrel{\eqref{Er1}}{=} c_l (\mcV) \\
&\stackrel{\eqref{Er2}}{=} c_l (\mcD_{X, 1}^{(\N-1)} \otimes \mcN) \notag \\
& \stackrel{\eqref{Er5}}{=}
  c_l (\mcT_X \otimes \mcN) + c_1 (\mcN) c_{l-1} (\mcT_X \otimes \mcN)\notag \\
  & \stackrel{\eqref{Er6}}{=}
  \sum_{i=0}^l \binom{n- i}{l-i} \cdot c_i (X) c_1 (\mcN)^{l-i} 
+ \sum_{i=0}^{l-1} \binom{n- i}{l-1-i} \cdot  c_i (X) c_1 (\mcN)^{l-i} \notag \\
& \, \ = c_l (X) +\sum_{i=0}^{l-1} \binom{n+1-i}{l-i} \cdot c_i (X) \cdot c_1 (\mcN)^{l-i}. \notag
\end{align}
 In the case  $l=1$, this composite equality  reads
\begin{align} \label{Er10}
c_1 (X) = -(n+1) c_1 (\mcN).
\end{align}
Moreover, by \eqref{Er11}, \eqref{Er10}, and the induction assumption,
the following sequence of equalities holds:
\begin{align} \label{Er9}
c_l (X) &\stackrel{\eqref{Er11}}{=} - \sum_{i=0}^{l-1} \binom{n+1-i}{l-i} \cdot c_i (X)  c_1  (\mcN)^{l-i} \\
&\stackrel{\eqref{Er10}}{=}
 - \sum_{i=0}^{l-1} \binom{n+1-i}{l-i} \cdot  \left( \frac{-1}{n+1}\right)^{l-i} \cdot c_i (X) c_1 (X)^{l-i} \notag \\
 & \hspace{-4mm} \stackrel{(\mr{induction})}{=}- \sum_{i=0}^{l-1} \binom{n+1-i}{l-i} \cdot \frac{1}{(n+1)^i} \cdot \binom{n+1}{i}  \cdot \left( \frac{-1}{n+1}\right)^{l-i} \cdot c_1(X)^i  c_1 (X)^{l-i} \notag \\
 & \ = c_1 (X)^l \cdot \frac{1}{(n+1)^l}  \cdot \sum_{i=0}^{l-1} (-1)^{l+1-i} \binom{n+1-i}{l-i} \cdot \binom{n+1}{i} \notag \\
& \ =c_1 (X)^l \cdot \frac{1}{(n+1)^l} \cdot \binom{n+1}{l} \cdot \left(1- \sum_{i=0}^{l} (-1)^{l-i} \binom{l}{l-i} \right)\notag \\
& \ =c_1 (X)^l \cdot \frac{1}{(n+1)^l} \cdot \binom{n+1}{l}. \notag
\end{align}
This completes the proof of the first assertion of (i).

The second  assertion  follows directly from the first one  and their  crystalline versions  
 can be  proved in  similar ways.

Finally,  assertion (ii) follows from assertion (i)  and the fact that $H^{2l}_{\mr{crys}}(X/W)$ is finitely generated  over $W$ (cf. ~\cite[Chap\,VII, Corollaire 1.1.2]{Be0}).
\end{proof}

\SSP
\bt \label{T016a} 
Let $\N \in \mbZ_{>0} \sqcup \{ \infty \}$ and let $X$ be a smooth quasi-projective  variety over $k$
admitting
an  $F^\N$-affine structure.
Then, the following assertions hold:
\begin{itemize}
\item[(i)]
In the case where
$\N \neq \infty$, 
  the   equality
\begin{align} \label{E45aaa}
 c_l (X)=0
\end{align}
holds  after reduction 
 modulo $p^{l \N}$  for each positive integer $l$.
In particular, if $X$ admits an 
$F^\infty$-affine structure, then the above equality holds after reduction 
module $p^{l \N}$ for any positive integers $\N$ and $l$.
Moreover, the same assertions hold for $c_l^{\mr{crys}}(X)$'s.
\item[(ii)]
In the case where $\N = \infty$ and $X$ is projective over $k$,
 the equality \eqref{E45aaa}   in $H^{2l}_{\mr{crys}}(X/W)$, where ``$c_1(X)$'' and ``$c_l(X)$''  are replaced by   $c_1^{\mr{crys}}(X)$ and $c_l^{\mr{crys}}(X)$,  respectively,   holds without any reduction for every  positive integer $l$.
\end{itemize}
  \et
\begin{proof}
 By assumption, $X$ admits an $F^\N$-affine structure $\mcS^\heartsuit$. 
It follows from 
Theorem \ref{P01111}, (ii),  that 
there exists a dormant affine-indigenous $\mcD_X^{(\N-1)}$-module ${^\A \mcV}^\diamondsuit := (\mcV, \nabla, \mcN, \delta)$  with
${^\A}\mcV^{\diamondsuit \Rightarrow \heartsuit} = \mcS^\heartsuit$.
The $\mcD_X^{(\N-1)}$-action $\nabla$ induces $\mcD_X^{(\N-1)}$-actions on $\mr{Ker}(\delta)$ and $\mcN$ via the natural inclusion $\mr{Ker}(\delta) \migiincl  \mcV$ and the surjection $\delta : \mcV \migisurj \mcN$  respectively; we shall
denote these $\mcD_X^{(\N-1)}$-actions  by $\nabla_{\mr{Ker}(\delta)}$ and $\nabla_\mcN$ respectively.
 Let us identity  $\mr{Ker}(\delta)$ with $\mcT_X \otimes \mcN$ via the composite isomorphism \eqref{Eh405}.
 Then, the tensor product of $\nabla_{\mr{Ker}(\delta)}$ and the dual of $\nabla_\mcN$ yields a $\mcD_{X}^{(\N-1)}$-action $\nabla_{\mcT_X}$ on $\mcT_X  \left(\cong (\mcT_X \otimes \mcN) \otimes \mcN^\vee \right)$.
 Since $\nabla$ has vanishing $p$-$(\N-1)$-curvature,
 both $\nabla_{\mr{Ker}(\delta)}$ and $\nabla_\mcN$ has vanishing $p$-$(\N-1)$-curvature, and hence, $\nabla_{\mcT_X}$ has vanishing $p$-$(\N-1)$-curvature.
 It follows that there exists a vector bundle $\mcW$ on $X$ with $F_X^{\N*}(\mcW) \cong \mcT_X$.
 This implies  $c_l (X) \left(:=c_l (\mcT_X) \right) = p^{l\N} \cdot c_l(\mcW)$ (cf.  ~\cite[the proof of Lemma 2.1]{G}).
 In particular, the equality $c_l (X) = 0$ holds  modulo $p^{l \N}$, as desired.
 This completes the proof of the first assertion of  (i).
 
 The second  assertion  follows directly from the first one  and their  crystalline versions  
 can be  proved in  similar ways.

 Finally,  assertion (ii) follows from assertion (i)  and the fact that $H^{2l}_{\mr{crys}}(X/W)$ is finitely generated  over $W$ (cf. ~\cite[Chap\,VII, Corollaire 1.1.2]{Be0}).
\end{proof}

\vspace{10mm}
\section{$F^\N$-projective/affine structure on projective/affine spaces} \label{S0010}\SSP

In this section, we 
study $F^\N$-projective and $F^\N$-affine structures 
 on  projective  and affine spaces, respectively.
Since projective  and affine spaces have global  coordinate charts, there exist the trivial constructions among these structures.
In the case of projective spaces,  we will show  the 
  uniqueness of  $F^\N$-projective structures of prescribed  level $\N$  (cf. Proposition \ref{P0104}). 
  On the other hand,  it will be  verified   that there  exist   many $F^\N$-affine structures on affine spaces; this  is an exotic phenomenon that occurs  in positive characteristic (cf. Proposition \ref{P04ppg0}). 
The main results of this section provide several  characterizations of  projective spaces
in terms of $F^\infty$-structures
  (cf. Theorems \ref{T01022} and \ref{T044}).
  These characterizations can be thought of as positive characteristic analogues of  
  the corresponding results for complex manifolds,  proved in ~\cite{KO1}, ~\cite{Ye},  and ~\cite{JR1}.

\LSP
\subsection{Uniqueness of $F^\N$-projective structures on projective spaces}
 \label{Sdghj1}

Let $\N$ be a positive integer.
To begin with, we construct  the trivial   $F^\N$-projective (resp., $F^\N$-affine)  structure 
on  the projective (resp., affine) space $\mbP^n$ (resp., $\mbA^n$), as well as the corresponding 
 dormant indigenous  $\mcD_{\mbP^n}^{(\N-1)}$-module (resp., dormant affine-indigenous $\mcD_{\mbA^n}^{(\N-1)}$-module). 

We shall set 
\begin{align} \label{Eh1060}
\mcS_{\N, \mr{triv}}^\heartsuit \ \left(\text{resp.,} \  {^\A\mcS}_{\N, \mr{triv}}^\heartsuit \right)
\end{align}
to be the subsheaf of $\mcP^{\text{\'{e}t}}_{\mbP^n}$ (resp., $\mcA^{\text{\'{e}t}}_{\mbA^n}$) 
 which, to any \'{e}tale scheme $U$ over  $\mbP^n$ (resp., $\mbA^n$), assigns the set  
\begin{align}
S_{\N, \mr{triv}}^\heartsuit(U) := 
\left\{ \overline{A} (\phi_U) \in \mcP^{\text{\'{e}t}}_{\mbP^n}(U) \, \Big| \, \overline{A} \in (\mr{PGL}_{n+1})_{\mbP^n}^{(\N)} (U) \right\} \\
\left(\text{resp.,} \  {^\A S}_{\N, \mr{triv}}^\heartsuit(U) := 
\left\{ \overline{A} (\phi_U) \in \mcA^{\text{\'{e}t}}_{\mbA^n}(U) \, \Big| \, \overline{A} \in (\QQ)_{\mbA^n}^{(\N)} (U) \right\} \right), \hspace{-5mm}\notag  
\end{align}
 where $\phi_U$ denotes the natural  \'{e}tale morphism  $U \migi \mbP^n$ (resp., $U \migi \mbA^n$).
That is to say, $\mcS_{\N, \mr{triv}}^\heartsuit$ (resp., ${^\A\mcS}_{\N, \mr{triv}}^\heartsuit$) is obtained as the 
smallest subsheaf of $\mcP^{\text{\'{e}t}}_{\mbP^n}$ (resp.,  $\mcA^{\text{\'{e}t}}_{\mbA^n}$) that is closed under the $(\mr{PGL}_{n+1})_{\mbP^n}^{(\N)}$-action (resp., the $(\QQ)_{\mbA^n}^{(\N)}$-action) and  contains the global section determined by the identity morphism of $\mbP^n$ (resp., $\mbA^n$).
Then, $\mcS_{\N, \mr{triv}}^\heartsuit$ (resp., ${^\A \mcS}_{\N, \mr{triv}}^\heartsuit$) specifies 
 an $F^{\N}$-projective (resp.,  $F^\N$-affine) structure on $\mbP^n$ (resp., $\mbA^n$).
 The formation of $S_{\N, \mr{triv}}^\heartsuit$ (resp.,  ${^\A S}_{\N, \mr{triv}}^\heartsuit$) is compatible with truncation to lower levels, so 
 the collection 
\begin{align} \label{Eh1061}
\mcS_{\infty, \mr{triv}}^\heartsuit := \{ \mcS_\N^\heartsuit\}_{\N \in \mbZ_{>0}} \ \left(\text{resp.,} \ {^\A\mcS}_{\infty, \mr{triv}}^\heartsuit := \{ {^\A\mcS}_\N^\heartsuit\}_{\N \in \mbZ_{>0}} \right)
\end{align}
 forms an $F^\infty$-projective (resp., $F^\infty$-affine) structure on $\mbP^n$ (resp., $\mbA^n$).

Next, we shall observe that there exists a canonical theta characteristic on the projective space $\mbP^{n} := \mr{Proj}(k[x_0, x_1, \cdots, x_n])$ defined as follows.
Let
\begin{align} \label{Effgh}
\eta_0 : \mcO_{\mbP^n}(-1) \migi \mcO_{\mbP^n}^{\oplus (n+1)}
\end{align}
be the $\mcO_{\mbP^n}$-linear injection given by $w \mapsto \sum_{i=0}^n w x_i\cdot  e_i$ for each local section $w \in \mcO_{\mbP^n}(-1)$, where $(e_0, \cdots, e_n)$ is the standard  basis of $\mcO_{\mbP^n}^{\oplus (n+1)}$.
The natural short exact sequence 
\begin{align}
0 \xrightarrow{}  \mcO_{\mbP^n}(-1) \xrightarrow{\eta_0} \mcO_{\mbP^n}^{\oplus (n+1)} \xrightarrow{} \mcO_{\mbP^n}^{\oplus (n+1)}/\mr{Im}(\eta_0) \xrightarrow{} 0
\end{align}
induces a composite  isomorphism
\begin{align} \label{Eh1047}
\mcO_{\mbP^n}(-1) \isom  \mr{det}(\mcO_{\mbP^n}^{\oplus (n+1)}/\mr{Im}(\eta_0))^\vee 
\otimes 
  \mr{det}(\mcO_{\mbP^n}^{\oplus (n+1)}) \isom 
   \mr{det}(\mcO_{\mbP^n}^{\oplus (n+1)}/\mr{Im}(\eta_0))^\vee,
\end{align}
where the second arrow arises from the  isomorphism $\mcO_{\mbP^n} \isom  \mr{det}(\mcO_{\mbP^n}^{\oplus (n+1)})$  given by $1 \mapsto e_0 \wedge \cdots \wedge e_n$.
Also, observe that 
the composite
\begin{align} \label{Eh1052}
\mcO_{\mbP^n}(-1) \xrightarrow{\eta_0}  \mcO_{\mbP^n}^{\oplus (n+1)} \xrightarrow{d^{\oplus (n+1)}}  \Omega_{\mbP^n} \otimes \mcO_{\mbP^n}^{\oplus (n+1)} \migisurj \Omega_{\mbP^n} \otimes (\mcO_{\mbP^n}^{\oplus (n+1)}/\mr{Im}(\eta_0))
\end{align}
is  $\mcO_{\mbP^n}$-linear and 
 induces  an isomorphism of $\mcO_{\mbP^n}$-modules
\begin{align} \label{E5568}
\mcO_{\mbP^n}(-1) \otimes (\mcO_{\mbP^n}^{\oplus (n+1)}/\mr{Im}(\eta_0))^\vee \isom \Omega_{\mbP^n}. 
\end{align}
By \eqref{Eh1047} and \eqref{E5568}, we obtain the following composite isomorphism:
\begin{align}
\mcO_{\mbP^n}(-1)^{\otimes (n+1)} &\isom
\mr{det}(\mcO_{\mbP^n}^{\oplus (n+1)}/\mr{Im}(\eta_0))^\vee \otimes  \mcO_{\mbP^n}(-1)^{\otimes n} \\
&\isom \mr{det}((\mcO_{\mbP^n}^{\oplus (n+1)}/\mr{Im}(\eta_0))^\vee \otimes \mcO_{\mbP^n}(-1)) \notag \\
&\isom \left(\mr{det} (\Omega_{\mbP^n}) =\right)   \omega_{\mbP^n}. \notag
\end{align}
Thus, the line bundle $\mcO_{\mbP^n}(-1)$ together with this composite specifies a theta characteristic  of $\mbP^n$.
Moreover, the restriction $\mcO_{\mbP^n}(-1) |_{\mbA^n}$ of this line bundle to $\mbA^n  \left(\subseteq \mbP^n \right)$ specifies  a theta characteristic of $\mbA^n$.
If  $\delta$ denotes the projection $\mcO_{\mbA^n}^{\oplus (n+1)} \migisurj \mcO_{\mbA^n}$ onto the $1$-st factor, i.e., the morphism given by $\sum_{i=0}^n w_i \cdot e_i \mapsto w_0$, then the composite $\delta \circ (\eta_0 |_{\mbA^n}) : \mcO_{\mbP^n}(-1) |_{\mbA^n} \migi \mcO_{\mbA^n}$ is an isomorphism;  we shall identify $\mcO_{\mbP^n}(-1) |_{\mbA^n}$ with the trivial line bundle $\mcO_X$ via this isomorphism.

Now, let us write
\begin{align} \label{Eh1057}
\mcV_{\N, \mr{triv}}^\diamondsuit := (\mcO^{\oplus (n+1)}_{\mbP^n}, (\nabla^{\mr{triv} (\N-1)}_{\mcO_{\mbP^n}})^{\oplus (n+1)}, \mcO_{\mbP^n}(-1)) \\ 
\left(\text{resp.,} \  {^\A\mcV}_{\N, \mr{triv}}^\diamondsuit := (\mcO^{\oplus (n+1)}_{\mbA^n}, (\nabla^{\mr{triv} (\N-1)}_{\mcO_{\mbA^n}})^{\oplus (n+1)}, \mcO_{\mbP^n}(-1) |_{\mbA^n}, \delta)
\right). \hspace{-12mm}\notag
\end{align}
The bijectivity of \eqref{E5568} shows  that, 
for each $\N \in \mbZ_{>0} \sqcup \{ \infty \}$,
the Kodaira-Spencer map
  $\mr{KS}^\circledast_{\mcV_{\N, \mr{triv}}^\diamondsuit}$ is an isomorphism, and hence,  $\mcV_{\N, \mr{triv}}^\diamondsuit$ (resp., ${^\A\mcV}_{\N, \mr{triv}}^\diamondsuit$)
forms  a dormant indigenous $\mcD_{\mbP^n}^{(\N-1)}$-module (resp., a dormant affine-indigneous $\mcD_{\mbA^n}^{(\N-1)}$-module). 
Moreover, recall the  isomorphism 
\begin{align} \label{Eh1400}
\mr{KS}^\circledast_{\mcV_{\N, \mr{triv}}^\diamondsuit} : \mcD_{\mbP^n, 1}^{(\N-1)} \otimes \mcO_{\mbP^n}(-1) \isom \mcO_{\mbP^n}^{\oplus (n+1)}
\end{align}
induced by $\mr{KS}_{\mcV_{\N, \mr{triv}}^\diamondsuit}$ (cf. \eqref{Eh101}).
By passing through  $\mr{KS}^\circledast_{\mcV_{\N, \mr{triv}}^\diamondsuit}$
(resp.,   $\mr{KS}^\circledast_{\mcV_{\N, \mr{triv}}^\diamondsuit}|_{\mbA^n}$ and $\delta \circ (\eta_0 |_{\mbA^n}) : \mcO_{\mbP^n}(-1) |_{\mbA^n} \isom \mcO_{\mbA^n}$),  
we may   identify the data $\mcO^{\oplus (n+1)}_{\mbP^n}$ (resp., $\mcO^{\oplus (n+1)}_{\mbA^n}$ and $\mcO_{\mbP^n}(-1) |_{\mbA^n}$) in
 \eqref{Eh1057} with $\mcD_{\mbP^n, 1}^{(\N-1)} \otimes \mcO_{\mbP^n}(-1)$ (resp., $\mcD_{\mbA^n, 1}^{(\N-1)}$ and $\mcO_{\mbA^n}$).
 With this in mind, $\mcV_{\N, \mr{triv}}^\diamondsuit$ (resp., ${^\A\mcV}_{\N, \mr{triv}}^\diamondsuit$) may be regarded  as an element of
 ${^\dagger}\DID{\N -1}{, \mcO_{\mbP^n}(-1)}$
  (resp., 
  ${^\dagger}\DAID{}{\N -1, \mbA^n}$).

\SSP
\bpr \label{P0104} 
The following assertions hold:
\begin{itemize}
\item[(i)]
Let $\N \in \mbZ_{>0} \sqcup \{ \infty \}$.
Then, 
the  following equality holds:
\begin{align}
{^\dagger \zeta}_{\N, \mcO_{\mbP^n}(-1)}^{\diamondsuit \Rightarrow \heartsuit} (\mcV_{\N, \mr{triv}}^{\diamondsuit})  =\mcS_{\N, \mr{triv}}^{\heartsuit}.
\end{align} 
Moreover,  if $p \nmid (n+1)$, then we have
\begin{align} \label{Eh1335}
\Proj{\N}{, \mbP^n}= \left\{\mcS_{\N, \mr{triv}}^\heartsuit \right\},
\hspace{10mm}
{^\dagger}\DID{\N -1}{, \mcO_{\mbP^n}(-1)}
 = \left\{ \mcV_{\N, \mr{triv}}^\diamondsuit\right\}. \notag
\end{align}
(The ``affine'' version
of this assertion  will be described in the next subsection.)
 \item[(ii)]
Let $\N$ be a positive integer with 
  $p^\N\nmid (n+1)$.
   Then,  there are no $F^\N$-affine structures on $\mbP^n$.
 In particular,   there are no $F^\infty$-affine structures on $\mbP^n$.
 \end{itemize}
 \epr
\begin{proof}
We shall prove assertion (i).
The former assertion follows directly  from the various definitions involved.
To prove the latter assertion, let us recall that, since $p \nmid (n+1)$, 
 $\overline{\zeta}_\N^{\diamondsuit \Rightarrow \heartsuit}$
 is bijective  (cf.  Theorem \ref{P01111}).
 Thus, the problem is reduced to proving the claim that  
the set  
${^\dagger}\DID{\N -1}{, \mcO_{\mbP^n}(-1)}$
  consists exactly of $\mcV_{\N, \mr{triv}}^\diamondsuit$.
Moreover,  by the identification 
${^\dagger}\DID{\infty}{, \mcO_{\mbP^n}(-1)}
 \cong \varprojlim_{\N \in \mbZ_{>0}}
 {^\dagger}\DID{\N -1}{, \mcO_{\mbP^n}(-1)}$
   resulting from Proposition \ref{P06p0hgh}, (iii),  it suffices to consider the case where
 $\N \neq \infty$.

Let  $\mcV^\diamondsuit := (\mcD_{X, 1}^{(\N-1)}\otimes \mcO_{\mbP^1}(-1), \nabla, \mcO_{\mbP^n}(-1))$ be a dormant indigenous $\mcD_{\mbP^n}^{(\N-1)}$-module  representing an element of 
${^\dagger}\DID{\N -1}{, \mcO_{\mbP^n}(-1)}$.
We identify the underlying vector bundle $\mcD_{X, 1}^{(\N-1)}\otimes \mcO_{\mbP^1}(-1)$  of $\mcV^\diamondsuit$ with $\mcO^{\oplus (n+1)}_{\mbP^n}$ via \eqref{Eh1400}.
Since $\nabla$ has vanishing $p$-$(\N-1)$-curvature, the natural morphism
 $\alpha : F_{\mbP^n/k}^{(\N)*}(\mcS ol (\nabla)) \migi  \mcO_{\mbP^n}^{\oplus (n+1)}$ is an isomorphism.
According to Lemma \ref{Lk0104} described below,
there exists an isomorphism $\beta_0 : \mcS ol (\nabla) \isom \mcO_{\mbP^{n (\N)}}^{\oplus (n+1)}$.
Denote by $\beta : F_{\mbP^n/k}^{(\N)*}(\mcS ol (\nabla)) \isom \mcO_{\mbP^n}^{\oplus (n+1)}$ the pull-back of $\beta_0$ by $F_{\mbP^n/k}^{(\N)}$, where we identify $F_{\mbP^n/k}^{(\N)*} (\mcO_{\mbP^{n (\N)}})$ with $\mcO_{\mbP^n}$  by passing to the natural isomorphism induced from $F_{\mbP^n/k}^{(\N)}$.
Because of the properness of $X$, the automorphism $\alpha \circ \beta^{-1}$ of $\mcO_{\mbP^n}^{\oplus (n+1)}$ may be described as an element of $\mr{GL}_{n+1}(k)$.
 Hence, 
  there exists an automorphism $\gamma_0$ of $\mcO_{\mbP^{n (\N)}}^{\oplus (n+1)}$ inducing $\alpha \circ \beta^{-1}$ via 
 pull-back by $F_{\mbP^n/k}^{(\N)}$.
The pull-back of  $\gamma_0 \circ\beta_0 : \mcS ol (\nabla) \isom \mcO_{\mbP^{n (\N)}}^{\oplus (n+1)}$ by $F_{\mbP^n/k}^{(\N)}$ 
 coincides with 
   $\alpha$.
 It follows from   the equivalence of categories \eqref{E3002} that   the isomorphism $\gamma_0 \circ\beta_0$ implies   the equality  $\nabla = (\nabla_{\mcO_X}^{\mr{triv} (\N-1)})^{\oplus (n+1)}$, i.e.,
  $\mcV^\diamondsuit = \mcV_{\N, \mr{triv}}^\diamondsuit$.
This completes the  proof of assertion (i).

Next,  assertion (ii)  follows from 
 Theorem \ref{T016}, (ii), together with 
the following equalities  in $CH^1 (\mbP^1)  \left(= \mr{Pic}(\mbP^n) \cong \mbZ \cdot c_1  (\mcO_{\mbP^n}(1)) \right)$:
\begin{align}
c_1 (\mbP^n) = c_1 (\omega_{\mbP^n}^\vee) = c_1 (\mcO_{\mbP^n}(n+1)) = (n+1) \cdot c_1  (\mcO_{\mbP^n}(1)).
\end{align}
This completes the proof of the proposition.
\end{proof}
\SSP

The following lemma was used in the proof of the above proposition.

\SSP
\ble \label{Lk0104}
 Let $\N$ be a positive integer and $\mcF$  a vector bundle on $\mbP^{n (\N)}$.
 Then, $\mcF$ is trivial (i.e., isomorphic to a direct sum of finite copies of the trivial line bundle) if and only if  the pull-back $F_{\mbP^n/k}^{(\N)*}(\mcF)$ is trivial.
 \ele
\begin{proof}
It suffices to consider the ``if'' part because the inverse direction is clear.
Suppose that $F_{\mbP^n/k}^{(\N)*}(\mcF)$ is trivial.
Let us take an arbitrary rational curve in $\mbP^{n(\N)}$, which may be 
 obtained  as the $\N$-th Frobenius twist $C^{(\N)}$ for some rational curve $C$ in $\mbP^n$.
Write
 $f : \mbP^1 \migi \mbP^n$  for the normalization of $C$.
 Hence, its base-change $f^{(\N)} : \mbP^{1 (\N)}  \left(= \mbP^1\right) \migi \mbP^{n (\N)}  \left(= \mbP^n \right)$ gives  the normalization of $C^{(\N)}$.
According to ~\cite[Theorem 1.1]{BS}, the problem  is reduced to proving the claim that the pull-back $f^{(\N)*}(\mcF)$ of $\mcF$ via $f^{(\N)}$ is trivial.
(We can apply the result of ~\cite{BS} because  $\mbP^n$ is separably rationally connected, as mentioned in  ~\cite[Chapter IV, Example 3.2.6]{Kol}.)
Recall  the Birkhoff-Grothendieck's theorem, asserting that any vector bundle on $\mbP^1$ is isomorphic to a direct sum of line bundles. Hence, we have a decomposition
\begin{align}
f^{(\N)*}(\mcF) \cong \bigoplus_{j=1}^m \mcO_{\mbP^{1(\N)}}(a_j),
\end{align}
where $m := \mr{rank}(\mcF)$ and $a_1 \leq a_2 \leq \cdots \leq a_m$.
This decomposition implies that 
\begin{align}\label{Ep190}
f^*(F_{X/k}^{(\N)*}(\mcF)) \cong F_{\mbP^1/k}^{(\N)*}(f^{(\N)*}(\mcF)) \cong F_{\mbP^1/k}^{(\N)*}(\bigoplus_{j=1}^m \mcO_{\mbP^{1(\N)}}(a_j)) \cong \bigoplus_{j=1}^m \mcO_{\mbP^1}(p^{\N}a_j).
\end{align}
On the other hand,  the triviality  assumption on $F_{\mbP^n/k}^{(\N)*}(\mcF)$  implies 
\begin{align} \label{Ep191}
f^*(F_{X/k}^{(\N)*}(\mcF)) \cong f^*(\mcO_{\mbP^{n}}^{\oplus m}) \cong \mcO_{\mbP^1}^{\oplus m}.
\end{align}
By \eqref{Ep190} and \eqref{Ep191},  the equalities $a_i =0$ hold for all $i=1,2, \cdots, m$.
Consequently, $f^{(\N)*}(\mcF)$ turns out to be trivial.
This completes the proof of the lemma.
\end{proof}

\LSP
\subsection{Nonuniqueness of $F^\N$-affine  structures on affine spaces} \label{Sffhy1}

Next, we focus on    $F^\N$-affine structures on  affine spaces.
Let  $\N \in \mbZ_{>0} \sqcup \{ \infty \}$.

\SSP
\bde \label{D0gg050} 
 Let $\mcS^\heartsuit$ be an $F^\N$-projective  (resp.,  $F^\N$-affine) structure on $X$.
 If $\N < \infty$,  then we shall say that $\mcS^\heartsuit$ is {\bf global} if 
 $\mcS^\heartsuit$ has a global section, or equivalently, 
 the  $\mr{PGL}_{n+1}$-bundle (resp., the  $\QQ$-bundle)  on $X^{(\N)}$ determined  by $\mcS^\heartsuit$ is trivial.
Also, if $\N = \infty$, then we shall say that  $\mcS^\heartsuit$ is {\bf global} if, for any $\N' \in \mbZ_{>0}$,  the $\N'$-th truncation $\mcS^\heartsuit |^{\langle \N' \rangle}$   is global.
    \ede
\SSP

Regarding global $F^\N$-projective and $F^\N$-affine structures,  we have the following   proposition.

\SSP
\bpr \label{P04pp0} 
The following assertions hold:
\begin{itemize}
\item[(i)]
$X$ admits a global $F^\N$-projective  (resp.,  $F^\N$-affine) structure
 if and only if 
there exists an \'{e}tale $k$-morphism $X \migi \mbP^n$ (resp.,  $X \migi \mbA^n$).
\item[(ii)]
Suppose that $X$ is proper over $k$.
Then, $X$ is isomorphic to $\mbP^n$ if $X$ admits a global $F^\N$-projective structure. 
Moreover, $X$ admits no global $F^\N$-affine structures.
\item[(iii)]
If $X$ admits a global $F^\N$-projective (resp.,  $F^\N$-affine) structure $\mcS^\heartsuit_\N$,
then, for any $\N' \geq \N$, $X$ also admits a global $F^{\N'}$-projective (resp., $F^\N$-affine) structure $\mcS^\heartsuit_{\N'}$ with  $\mcS^\heartsuit_{\N'} |^{\langle \N \rangle} = \mcS_\N^\heartsuit$.
\end{itemize}
  \epr
\begin{proof}
The only if part of assertion (i) follows  from the definition of  a global $F^\N$-projective ($F^\N$-affine) structure.
As for the if part, given an \'{e}tale $k$-morphism $\phi : X \migi \mbP^n$ (resp., $\phi : X \migi \mbA^n$), we obtain a global  $F^\N$-projective (resp.,  $F^\N$-affine) structure $\mcS_\phi^\heartsuit$ on $X$ defined as the subsheaf of $\mcP_X^{\text{\'{e}t}}$ (resp., $\mcA_X^{\text{\'{e}t}}$)  which, to each \'{e}tale  scheme $U$ over  $X$, assigns the set of \'{e}tale morphisms  of the form $\overline{A} (\phi |_U)$ (cf. \eqref{Eg69}) for some $\overline{A} \in (\mr{PGL}_{n+1})_X^{(\N)}(U)$ (resp., $\overline{A} \in (\QQ)_X^{(\N)}(U)$).

The former assertion of  (ii) follows from assertion (i) and  the well-known fact that $\mbP^n$ has no nontrivial \'{e}tale coverings.
In  fact, the properness  on $X$   implies that
an \'{e}tale morphism $X \migi \mbP^n$ defining a global $F^\N$-projective structure
becomes an \'{e}tale covering, and hence, an isomorphism.
Moreover, the latter assertion can be verified as follows.
Suppose that $X$ admits a global $F^\N$-affine structure, which contains an \'{e}tale morphism $\phi : X \migi \mbA^n$.
Since $\mbA^n$ is affine, the image of $\phi$ must be a point.
But, it contradicts the fact that any  \'{e}tale morphism of varieties is open.

Finally, assertion (iii) follows from assertion (i) together with the construction of $\mcS_\phi^\heartsuit$ described above.
\end{proof}
\SSP

Now, let us
discuss
  the affine space $\mbA^n$.
 Consider  $\mbA^{n}  \left(\subseteq \mbP^n \right)$ as the affine scheme associated to the polynomial ring  
$R_0 := k [t_1, \cdots, t_n]$, where $t_i := x_i/x_0$ ($i=1, \cdots, n$).
Also,  for a positive integer $\N$, we identify $\mbA^{n (\N)}$ with 
the affine scheme associated to the subring $R_\N := k [t_1^{p^\N}, \cdots, t_n^{p^\N}]$ of $R_0$ and $F_{\mbA^n/k}^{(\N)}$ with the morphism $\mbA^n \migi \mbA^{n (\N)}$ induced by
the inclusion $R_\N \migiincl R_0$.

We shall set
\begin{align}
\mr{End}^{\text{\'{e}t}}_{\mbA^n} := \left\{(f_i)_{i=1}^n \in R_0^n  \, \Big| \, \mr{det} \left(\partial f_i/\partial t_j \right)_{ij} \in k^\times \right\}.
\end{align}
A $\QQ (R_\N)$-action on 
 $\mr{End}^{\text{\'{e}t}}_{\mbA^n}$ is defined as follows.
Take  $\overline{A} \in \QQ (R_\N)$ and $\vec{f} := (f_i)_i \in \mr{End}^{\text{\'{e}t}}_{\mbA^n}$.
The element $\overline{A}$
may be described uniquely as $\overline{\begin{pmatrix} 1 & {\bf 0} \\ {^t {\bf a}} & A'\end{pmatrix}}$
for some ${\bf a} := (g_i)_{i=1}^n  \in R_\N^n$ and  $A' \in \mr{GL}_n (R_\N)$.
Denote by  $\overline{A} (\vec{f})$ the element of $R_0^n$ defined as
\begin{align} \label{Eh1701}
\overline{A} (\vec{f}) : = (f_1^{\overline{A}}, \cdots, f_n^{\overline{A}}), \  \ 
\text{where} \  f_i^{\overline{A}} := g_i + \sum_{j=1}^n a_{ij} \cdot f_j \hspace{5mm} (i =1, \cdots, n).
\end{align}
Then,  the following equalities hold:
\begin{align}
\mr{det} \left(\partial f_i^{\overline{A}}/\partial t_j \right)_{ij} = 
\mr{det} \left(\sum_{l=1}^na_{il} \cdot (\partial f_l/\partial t_j) \right)_{ij} = \mr{det}(A') \cdot \mr{det} (\partial f_i/\partial t_j)_{ij},
\end{align}
where the rightmost side belongs to $k^\times \left(=R_0^\times \right)$ because
of the assumptions on $A'$ and  $\vec{f}$.
This implies that $\overline{A} (\vec{f})$ belongs to $\mr{End}^{\text{\'{e}t}}_{\mbA^n}$.
The resulting assignment $(\overline{A}, \vec{f}) \mapsto \overline{A}(\vec{f})$ is verified to define a $\QQ (R_\N)$-action on 
 $\mr{End}^{\text{\'{e}t}}_{\mbA^n}$. 
 In particular, we obtain the orbit set
 \begin{align}
 \mr{End}^{\text{\'{e}t}}_{\mbA^n}/\QQ (R_\N)
 \end{align}
of $\mr{End}^{\text{\'{e}t}}_{\mbA^n}$ with respect to this action.

\SSP
\bpr \label{P04ppg0} 
Let $\N$ be a positive integer.
Then, the following assertions hold: 
\begin{itemize}
\item[(i)]
Any $F^\N$-projective  (resp., $F^\N$-affine) structure on $\mbA^n$ is global.
In particular, the map 
\begin{align} \label{Eh1709}
\Proj{\N'}{, \mbA^n}
\migi 
\Proj{\N}{, \mbA^n}
\ 
\left(\text{resp.,} \ 
\Aff{\N'}{, \mbA^n}
\migi
\Aff{\N}{, \mbA^n}
  \right)
\end{align}
 (for any $\N' > \N$) given by truncation is surjective. 
\item[(ii)]
There exists a canonical bijection  of sets
\begin{align} \label{Eh1669}
\mr{End}^{\text{\'{e}t}}_{\mbA^n}/\QQ (R_\N) \isom 
\Aff{\N}{, \mbA^n}.
 \end{align}
\end{itemize}
  \epr
\begin{proof}
First, we shall prove the former assertion of (i).
We only consider the $F^\N$-affine case  because the proof of the $F^\N$-projective case is simpler.
Let $\mcS^\heartsuit$ be an $F^\N$-affine structure on $\mbA^n$.
According to  Theorem \ref{P01111}, (ii),
there exists    a dormant affine-indigenous $\mcD_{\mbA^n}^{(\N-1)}$-module  ${^\A\mcV}^\diamondsuit : = (\mcV, \nabla, \mcN, \delta)$ with 
${^\A}\xi_\N^{\diamondsuit \Rightarrow \heartsuit} ({^\A}\mcV^\diamondsuit) = \mcS^\heartsuit$.
The vector bundle $\mcS ol (\nabla)$ on $\mbA^{n(\N)}$ together with its  subbundle $\mcS ol (\nabla |_{\mr{Ker}(\delta)})  \left(= \mcS ol (\nabla) |_{\mr{Ker}(\delta)} \right)$ corresponds, via projectivization, to    the $\QQ$-bundle  determined by $\mcS^\heartsuit$.
Recall here  that   every  vector bundle on $\mbA^n$ is trivial by the Quillen-Suslin theorem,  and that  every  extension of vector bundles split since $\mbA^n$ is affine.
Hence,  the inclusion $\mcS ol (\nabla |_{\mr{Ker}(\delta)}) \migiincl \mcS ol (\nabla)$ may be identified with the inclusion $\mcO_{\mbA^n}^{\oplus n} \migiincl \mcO_{\mbA^n}^{\oplus (n+1)}$ into the latter  $n$ factors.
This implies that  the  $\QQ$-bundle determined by $\mcS^\heartsuit$ is trivial, i.e., $\mcS^\heartsuit$ is global.
This completes the proof of the former assertion of (i).
The latter assertion of (i) follows from the former assertion and Proposition \ref{P04pp0}, (iii).

Next, we shall consider assertion (ii).
Let us take $\vec{f} := (f_i)_i \in \mr{End}^{\text{\'{e}t}}_{\mbA^n}$.
Since $\mr{det} (\partial f_i/\partial t_j)_{ij} \in k^\times$,  the endomorphism $\phi_{\vec{f}}$ of $\mbA^n$ given by $t_i \mapsto f_i$ ($i=1, \cdots, n$) is \'{e}tale.
In other words, $\phi_{\vec{f}}$ specifies a global section of $\mcA^{\text{\'{e}t}}_{\mbA^n}$.
Denote by $\mcS_{\vec{f}}^\heartsuit$ the global $F^\N$-affine structure on $\mbA^n$ generated by $\phi_{\vec{f}}$.
To be precise, $\mcS_{\vec{f}}^\heartsuit$ denotes the subsheaf of $\mcA^{\text{\'{e}t}}_{\mbA^n}$ determined  by assigning $U \mapsto \left\{ \overline{A} (\phi_{\vec{f}} |_U) \, \Big| \, \overline{A} \in (\QQ)_{\mbA^n}^{(\N)}(U)\right\}$.
The assignment $\vec{f} \mapsto \mcS_{\vec{f}}^\heartsuit$ defines a map of sets
\begin{align} \label{Eh1509}
\mr{End}^{\text{\'{e}t}}_{\mbA^n} \migi 
\Aff{\N}{, \mbA^n}.
\end{align}
By the definition of the $\QQ (R_\N)$-action on $\mr{End}^{\text{\'{e}t}}_{\mbA^n}$, one verifies that
 two elements $\vec{f}_1$, $\vec{f}_2$ of $\mr{End}^{\text{\'{e}t}}_{\mbA^n}$ specify the same $\QQ (R_\N)$-orbit if and only if 
 $\mcS^\heartsuit_{\vec{f}_1} = \mcS^\heartsuit_{\vec{f}_2}$.
 Hence, the map \eqref{Eh1509} factors through the quotient $\mr{End}^{\text{\'{e}t}}_{\mbA^n} \migisurj \mr{End}^{\text{\'{e}t}}_{\mbA^n}/\QQ (R_\N)$ and the resulting map 
 $\mr{End}^{\text{\'{e}t}}_{\mbA^n}/\QQ (R_\N) \migi 
 \Aff{\N}{, \mbA^n}$
  is injective.
 This map is also surjective by the former  assertion of (i).
 Consequently, we have obtained the desired bijection.
 This completes the proof of the proposition.
\end{proof}
\SSP

By applying  the above proposition, we can conclude the nonuniqueness of $F^\N$-affine structures on $\mbA^n$ for $\N>1$, as follows.

\SSP
\bpr \label{P0hjk104}
Let $\N \in \mbZ_{>0} \sqcup \{ \infty \}$.
Then, the following assertions hold:
\begin{itemize}
\item[(i)]
The following equalities hold:
\begin{align}
{^{\dagger} \xi}_\N^{\diamondsuit \Rightarrow \heartsuit} ({^\A\mcV}_{\N, \mr{triv}}^{\diamondsuit}) 
 = {^\A\mcS}^\heartsuit_{\N, \mr{triv}}.
\end{align} 
\item[(ii)]
Suppose   that  $(n, \N) = (1, 1)$.
Then, 
we have
\begin{align} \label{Eh1365}
\Aff{1}{, \mbA^1}
 &= \big\{{^\A\mcS}_{1, \mr{triv}}^\heartsuit \big\},
\hspace{10mm}
{^\dagger}\DAID{0}{, \mbA^1}
= \big\{ {^\A\mcV}_{1, \mr{triv}}^\diamondsuit\big\}.
\end{align}
\item[(iii)]
Suppose that $\N\neq 1$.
Then,  the sets
$\Aff{\N}{, \mbA^n}$
and 
${^\dagger}\DAID{\N -1}{, \mcO_{\mbP^n}(-1)|_{\mbA^n}}$
  are (nonempty and) not singletons.
 \end{itemize}
 \epr
\begin{proof}
Assertion (i) follows from the various definitions involved.

Next, we shall prove the equalities  \eqref{Eh1365} asserted in  (ii).
By  Theorem \ref{P01111}, (ii),  and Proposition \ref{P06p0gh}, (ii),
 it suffices to show that 
the  set 
$\Aff{1}{, \mbA^1}$
consists only of ${^\A\mcS}_{1, \mr{triv}}^\heartsuit$.
It follows from  Proposition \ref{P04ppg0}, (ii), that  
$\Aff{1}{, \mbA^1}$
corresponds bijectively to the set $\mr{End}^{\text{\'{e}t}}_{\mbA^1}/\mr{PGL}_2^\mbA(R_1)$.
This set is nonempty because it has an element corresponding to ${^\A\mcS}_{1, \mr{triv}}^\heartsuit$, which may be represented by the element $t_1 \in \mr{End}^{\text{\'{e}t}}_{\mbA^1} \subseteq  R_0  \left(=k[t_1]\right)$.
Now, let us take an arbitrary element of $\mr{End}^{\text{\'{e}t}}_{\mbA^1}$, given as an element $f \in k [t_1]$.
Since $d f /dt_1 \in k^\times$, we can express $f$ as $g + a \cdot t_1$ for some $a \in k^\times$ and $g \in k[t^p_1]  \left(= R_1 \right)$.
Then, by the definition of the $\mr{PGL}_2^\mbA (R_1)$-action on $\mr{End}^{\text{\'{e}t}}_{\mbA^1}$, 
we can verify  that $f$ and $t_1$ specify the same $\mr{PGL}_2^\mbA (R_1)$-orbit.
Consequently, $\mr{End}^{\text{\'{e}t}}_{\mbA^1}/\mr{PGL}_2^\mbA (R_1)$, as well as  
$\Aff{1}{, \mbA^1}$,
 is a singleton.
The proof of assertion (ii) is therefore complete.

Finally, assertion (iii) for finite  $\N>1$
follows from the fact that 
the $n$-tuple  $(t_1-t^p_1, \cdots,  t_n -t^p_n) \in R_\N^n$ defines an element of $\mr{End}^{\text{\'{e}t}}_{\mbA^n}$ but does not belong to the $\QQ (R_\N)$-orbit of $(t_1, \cdots, t_n) \in \mr{End}^{\text{\'{e}t}}_{\mbA^1}$.
The case of infinite $\N$ follows from the finite cases together with the surjectivity of \eqref{Eh1709}.
This completes the proof of the proposition.
\end{proof}

\LSP
\subsection{Characterization of projective spaces I (Stratified fundamental group)} \label{SS197}

Denote by  $\mr{Str}(X)$  the category of $\mcD_X^{(\infty)}$-modules,  which  is equivalent to the category of $F$-divided sheaves on $X$ (cf. \eqref{E3001}) and  forms an abelian $k$-linear  rigid tensor category.
Once we
choose a $k$-rational point $x : \mr{Spec}(k) \migi X$,
the assignment $(\mcV, \nabla) \mapsto x^*(\mcV)$ defines a functor $\omega_x : \mr{Str}(X) \migi \mr{Vec}_k$, where $\mr{Ver}_k$ denotes the category of finite dimensional $k$-vector spaces.
The category  $\mr{Str}(X)$ together with the functor $\omega_x$  forms a neutral Tannakian category (cf. ~\cite[Section 2.2]{dS}).
In particular, there exists 
a pro-algebraic group $k$-scheme $\pi^\mr{str}_{1}(X, x)$
such that $\omega_x$ induces an equivalence between $\mr{Str}(X)$ and the category $\mr{Rep}_k (\pi^\mr{str}_{1}(X, x))$ of finite dimensional $k$-representations of $\pi^\mr{str}_{1}(X, x)$.
The pro-algebraic group $\pi^\mr{str}_{1}(X, x)$ will be referred to as 
 the {\bf stratified fundamental group} of $X$ (at $x$).
We sometimes use the notation $\pi^\mr{str}_{1}(X)$ to denote $\pi^\mr{str}_{1}(X, x)$ if there is no fear of confusion.
As an example, we know (cf.  ~\cite[Theorem 2.2]{G})  that the stratified fundamental group of $\mbP^n$ is trivial, i.e., 
every $\mcD_{\mbP^n}^{(\infty)}$-module is isomorphic to a direct sum of finite copies of $(\mcO_{\mbP^n}, \nabla_{\mcO_{\mbP^n}}^{\mr{triv}(\infty)})$.

Then,  we prove the  following assertion, which is  
the positive characteristic version of 
 ~\cite[Theorem 4.4]{KO1}.

\SSP
\bt \label{T01022} 
Suppose that $X$ is projective over $k$ and 
admits 
 an $F^\infty$-projective structure.
If $\pi_1^\mr{str}(X)$ is trivial, then
 $X$ is isomorphic to $\mbP^n$.
 \et
\begin{proof}
Suppose that there exists an $F^\infty$-projective structure $\mcS_\infty^\heartsuit := \{ \mcS_\N^\heartsuit \}_{\N \in \mbZ_{>0}}$  on $X$.
For each $\N \in \mbZ_{>0}$,  we shall denote by $\mcE_\N$ the  $\mr{PGL}_{n+1}$-bundle on $X^{(\N)}$ determined by   $\mcS_\N^\heartsuit$.
 If we are given  an arbitrary faithful representation $\mr{PGL}_{n+1} \migi \mr{GL}_l$ (for $l >0$), then 
the $F$-divided sheaf corresponding to the collection $\{ \mcE_\N \times^{\mr{PGL}_{n+1}} \mr{GL}_l \}_{\N \in \mbZ_{>0}}$ is trivial because of the assumption  $\pi_1^\mr{str}(X) =1$.
It follows from   
 ~\cite[Lemma 3.5, (ii)]{BS} that
 the  $\mr{PGL}_{n+1}$-bundle  $\mcE_\N$ is trivial for every  $\N \in \mbZ_{> 0}$.
 (Note that the result in {\it loc.\,cit.} is asserted under the assumption that $X$ is separably rationally connected. But, this assumption was not  used in its proof.)
  Hence,  $\mcS^\heartsuit$ turns out to be global, so
 the assertion follows from Proposition \ref{P04pp0}, (ii).
\end{proof}
\SSP

\begin{rema}[Formulation using the \'{e}tale fundamental group] \label{R1k3401}
We shall recall the relationship between $\pi^\mr{strat}_{1}(X, x)$ and the \'{e}tale fundamental group $\pi_1^{\text{\'{e}t}} (X) := \pi_1^{\text{\'{e}t}} (X, x)$  of $X$ at $x$.
 According to ~\cite[Proposition 13]{dS},
 there exists a natural quotient homomorphism of group schemes
$\pi^\mr{str}_{1}(X) \migisurj \pi_1^{\text{\'{e}t}} (X)$,  where we regard $\pi_1^{\text{\'{e}t}} (X)$ as a constant group scheme.
  In particular, $\pi_1^{\text{\'{e}t}} (X)$ is trivia if $\pi^\mr{str}_{1}(X)$ is trivial.
 Moreover, it follows from ~\cite[Theorem 1.1]{EM} that  the converse is true if $X$ is a projective smooth variety.
  This result may be thought of as a positive characteristic analogue of the Malcev-Grothendieck theorem (cf. ~\cite{Mal}, ~\cite{Gro1}).
  By this fact and the above theorem, we conclude  that {\it any $n$-dimensional smooth projective variety admitting an $F^\infty$-projective structure with trivial $\pi^{\text{\'{e}t}}_1(X)$ is isomorphic to the projective space $\mbP^n$.}
  
 Also, according to 
 ~\cite[Corollaire 3.6]{De2}, 
 any separably rationally connected variety has trivial  $\pi^{\text{\'{e}t}}_1(X)$. Hence,  the above theorem remains true when the assumption ``$\pi^{\text{str}}_1(X)$ is trivial'' is replaced  by  ``$X$ is separably rationally connected''.
 See Section $2$ in {\it loc.\,cit}. for the definition of a  separably rationally connected variety.
 \end{rema}
\SSP

Let us consider   representations of  
the stratified fundamental group
associated to $F^\infty$-projective  and $F^\infty$-affine structures.
For  algebraic groups $G$, $H$ over $k$, then we write
\begin{align} \label{Eu300}
\mr{Out}(G, H) := \mr{Hom}(G, H)/\mr{Inn}(H),
\end{align}
 where $\mr{Hom}(G, H)$ denotes  the set of homomorphisms from $G$ to $H$ and $\mr{Inn}(H)$  denotes the group of inner automorphisms of $H$,  acting on $\mr{Hom}(G, H)$ in the natural manner;
we shall call each element of $\mr{Out}(G, H)$ as an {\bf outer homomorphism} from $G$ to $H$.

Let us fix 
 a $k$-rational point  $x$  of $X$ and 
  an $F^\infty$-projective structure  $\mcS^\heartsuit$ 
 on $X$.
  Suppose that $p \nmid (n+1)$.
  Since the map $\overline{\zeta}_\infty^{\diamondsuit \Rightarrow \heartsuit} : 
\eDID{\infty}{}
\migi 
\Proj{\infty}{}$
  is bijective (cf. Proposition \ref{C06p0gh}, (i)),
there exists  a dormant indigenous $\mcD_X^{(\infty)}$-module  $\mcV^\diamondsuit := (\mcV, \nabla, \mcN)$  with
$\zeta_{\infty}^{\diamondsuit \Rightarrow \heartsuit} (\mcV^\diamondsuit) = \mcS^\heartsuit$.
The $\mcD_X^{(\infty)}$-module $(\mcV, \nabla)$ determines 
the $k$-vector space $x^*(\mcV)$ endowed with a $\pi^\mr{str}_1 (X, x)$-action.
By choosing a basis of $x^*(\mcV)$, we obtain 
a  representation $\pi_1^{\mr{str}}(X, x) \migi \mr{GL}_{n+1}$.
The class in $\mr{Out}(\pi_1^{\mr{str}}(X, x), \mr{PGL}_{n+1})$ represented by 
the composite $\pi_1^{\mr{str}}(X, x) \migi \mr{GL}_{n+1} \migisurj \mr{PGL}_{n+1}$ does not depend on the choices of both  $\mcV^\diamondsuit$ and the basis of $x^*(\mcV)$ (i.e., depend only on $\mcS^\heartsuit$).
Thus, we obtain a well-defined  outer homomorphism
\begin{align} \label{Eh9010}
\rho_{\mcS^\heartsuit} 
 \in \mr{Out}(\pi_1^\mr{str} (X, x), \mr{PGL}_{n+1}).
\end{align}

Next, 
let ${\mcS}^\heartsuit$ be an $F^\infty$-affine structure on $X$.
According to 
Proposition \ref{C06p0gh}, (ii),
there exists a dormant affine-indigenous $\mcD_X^{(\infty)}$-module ${^\A}\mcV^\diamondsuit := (\mcV, \nabla, \mcN, \delta)$
with ${^\A}\xi_{\infty}^{\diamondsuit \Rightarrow \heartsuit} ({^\A}\mcV^\diamondsuit) = {\mcS}^\heartsuit$.
Just as in the above discussion, the $\mcD_{X}^{(\infty)}$-module 
$(\mcV, \nabla)$ induces, up to conjugation, a homomorphism
$\pi_1^\mr{str} (X, x) \migi \mr{GL}_{n+1}$.
Moreover, the $\mcD_X^{(\infty)}$-submodule 
$(\mr{Ker}(\delta), \nabla |_{\mr{Ker}(\delta)})$ determines 
a $\QQ$-reduction of this homomorphism, i.e., a homomorphism
$\pi_1^\mr{str} (X, x) \migi \mr{GL}^\A_{n+1} := \mr{GL}_{n+1} \times_{\mr{PGL}_{n+1}}\QQ$.
By composing it with  the quotient $\mr{GL}_{n+1}^\A \migisurj \QQ$, we obtain
a well-defined outer homomorphism
\begin{align} \label{Eh9011}
{^\A}\rho_{\mcS^\heartsuit} 
\in \mr{Out}(\pi_1^\mr{str} (X, x), \QQ).
\end{align}

We shall refer to $\rho_{\mcS^\heartsuit}$  (resp., ${^\A}\rho_{\mcS^\heartsuit}$) constructed above as the {\bf monodromy representation} associated to the $F^\infty$-projective  (resp., $F^\infty$-affine) structure  $\mcS^\heartsuit$.

\LSP
\subsection{Characterization of projective spaces II (Embedded rational curves)} \label{S0331}

Next, let us prove another  characterization of  projective spaces using the existence of an $F^\infty$-projective structure.
According to the main theorem in ~\cite{Ye},
a complex  Fano variety  admitting  a projective structure is isomorphic to the projective space of the same dimension. 
Also, we can see its generalization in ~\cite[Theorem 4.1]{JR1}, saying that
if a compact  K\"{a}her manifolds admitting a projective structure (or more generally, a projective connection) contains a rational curve, then it is isomorphic to the projective space.
The crucial step for this result is to show that $f^*(\mcT_X)$ is ample for an arbitrary nontrivial morphism $f: \mbP^1_\mbC \migi X$. Once this is proved, we can use  a result  of S. Mori (cf. ~\cite[Theorem 8]{M1}) to conclude the assertion.
A similar argument gives us  the positive characteristic analogue of this assertion, as follows.

\SSP
\bt \label{T044} 
Let $X$ be a smooth projective variety over $k$ of dimension $n$ admitting an $F^\infty$-projective structure.
If $X$ contains a rational curve, then $X$ is isomorphic to $\mbP^n$.
In particular, $\mbP^n$ is the only Fano variety of dimension $n$ admitting an  $F^\infty$-projective structure.
  \et
\begin{proof}
Since  
any Fano variety contains a rational curve (cf. ~\cite[Chapter V, Theorem 1.6, (1.6.1)]{Kol}), it suffices to prove the former assertion.

Let $C$ be any rational curve in $X$ and $f : \mbP^1 \migi X$ its normalization.
According to Mori's theorem (cf. ~\cite{M1},  or ~\cite[Theorem 4.2]{MP}),
the problem is reduced to proving the claim that $f^*(\mcT_X)$ is ample.
Since any vector bundle on $\mbP^1$ decomposes into a finite direct sum of line bundles, we can write
\begin{align} \label{Ep79}
f^*(\mcT_X) \cong \bigoplus_{i=1}^n \mcO_{\mbP^1} (a_i),
\end{align}
where $a_1 \leq a_2 \leq \cdots \leq a_n$.
Since the differential  $\mcT_C \left(\cong \mcO_{\mbP^1}(2) \right)  \migi f^*(\mcT_X)$ of $f$ is injective, the inequality 
\begin{align} \label{Ep8000}
a_n \geq 2
\end{align}
holds.
Now, 
let us take   an $F^\infty$-projective structure $\mcS^\heartsuit$ on $X$.
If  $\mcE_\infty := \{ \mcE_\N \}_{\N \in \mbZ_{> 0}}$ denotes 
the $F$-divided  $\mr{PGL}_{n+1}$-torsor  induced from $\mcS^\heartsuit$,
then  $\mcE_0$ admits   a $\PP$-reduction $\mcE_\mr{red}$ in a natural manner (cf. the proof of Theorem \ref{P01111}, (i)).
By pulling-back this data via $f$, we obtain  an $F$-divided  $\mr{PGL}_{n+1}$-torsor  $f^*(\mcE_\infty)$ on $\mbP^1$ and a $\PP$-reduction $f^*(\mcE_\mr{red})$ on $f^*(\mcE_0)$.
The equality  $\pi_1^\mr{str}(\mbP^1) = \{1 \}$  implies that $f^*(\mcE_\infty)$ is trivial, i.e.,
obtained from the rank $(n+1)$ trivial $F$-divided sheaf on $\mbP^1$, or equivalently,  from the rank $(n+1)$ trivial $\mcD_{\mbP^1}^{(\infty)}$-module $(\mcO_{\mbP^1}^{\oplus (n+1)}, (\nabla_{\mcO_{\mbP^1}}^{\mr{triv}(\infty)})^{\oplus (n+1)})$ (cf. \eqref{E3001}).
Denote by $\mcN$ the line subbundle of $\mcO_{\mbP^1}^{\oplus (n+1)}$ determined by $\mcE_\mr{red}$.
Then, we have the following sequence of isomorphisms:
\begin{align}
f^*(\mcT_X) \xrightarrow{\sim} f^*(\mcA d_{\mcE_0})/f^*(\mcA d_{\mcE_\mr{red}}) \xrightarrow{\sim}
\mcA d_{f^*(\mcE_0)}/\mcA d_{f^*(\mcE_\mr{red})} \xrightarrow{\sim} \mcN^\vee \otimes (\mcO^{\oplus (n+1)}_{\mbP^1}/\mcN),
\end{align}
where the first arrow denotes the pull-back of the Kodaira-Spencer map 
$\mr{KS}_{(\mcE_0, \nabla^{\mr{can}(0)}_{\mcE_1}, \mcE_\mr{red})}$ in the sense of ~\cite[Definition 4.4.1]{Wak12}
(composed with the inverse to  the isomorphism displayed preceding ~\cite[Example 4.1.3]{Wak12}) and 
the last arrow follows from the discussion in Remark \ref{R13401}.
The composite of these isomorphisms induces the following exact sequence: 
\begin{align} \label{Eh2031}
0 \longmigi \mcN \longmigi \mcO_{\mbP^1}^{\oplus (n+1)} \longmigi f^*(\mcT_X) \otimes \mcN \left(\cong \bigoplus_{i=1}^n \mcN \otimes \mcO_{\mbP^1}(a_i) \right) \longmigi 0.
 \end{align}
The injectivity of the second arrow in \eqref{Eh2031} implies the inequality $l := \mr{deg}(\mcN) \leq 0$.
On the other hand, the surjectivity of the third arrow implies the ineqaulity 
\begin{align} \label{Ep78}
a_i + l \geq 0
\end{align}
 for any $i =1, \cdots,  n$.
By the decomposition  \eqref{Ep79} and the exactness of \eqref{Eh2031}, we have
\begin{align} \label{Eh2032}
0  \ = \mr{deg}(\mcO_{\mbP^1}^{\oplus (n+1)}) 
 \stackrel{\eqref{Eh2031}}{=} \mr{deg}(\mcN) + \mr{deg} (f^*(\mcT_X) \otimes \mcN)
 \stackrel{\eqref{Ep79}}{=} l(n+1) + \sum_{i=1}^n a_i. 
\end{align}
Because of the inequality \eqref{Ep8000} and \eqref{Eh2032},
the integer $l$ must be negative.
Hence,  the inequality   \eqref{Ep78}  implies 
 the  inequality  $a_i >0$,  so  $f^* (\mcT_X)$ is ample.
 This completes the proof of the assertion.
 \end{proof}

\LSP
\subsection*{Acknowledgements} \label{S90fffff25}

We would  like to thank
Professor Michel Gros
  for 
his  helpful 
 comments on the theory of  differential operators   of finite  level.
We are also grateful to 
Professor  Yuichiro Hoshi 
 for carefully reading  the preprint of this paper and for his valuable suggestions.
We wish to express our deepest gratitude to the anonymous referees for  their very careful  reading of this paper  and for many insightful comments and suggestions, which greatly improved the quality and clarity of this paper.
 We are  also thankful  for the many inspiring conversations we have had with 
algebraic varieties equipped with  Frobenius-projective or Frobenius-affine structures,  that continue to live and breathe  in the world of mathematics!
Finally, this  work was   partially  supported by JSPS KAKENHI Grant Numbers 21K13770 and  25K06933.

\LSP

\end{document}